\documentclass[hidelinks,onefignum,onetabnum]{siamart250211}

\usepackage[T1]{fontenc}
\usepackage{geometry}
\usepackage{graphicx} 
\usepackage{pgfplots,tikz,float,tabularx}
\pgfplotsset{compat=1.18}
\usetikzlibrary{plotmarks}
\usetikzlibrary{external}
\usepackage{amsmath,amsfonts, amssymb, makecell,xcolor,enumitem,hyperref,bm}
\usepackage{algorithmic}

\usepackage{chngcntr}
\counterwithout{algorithm}{section}
\usepackage{setspace}
\allowdisplaybreaks
\newtheorem{assumption}{Assumption}
\definecolor{tured}{RGB}{190,30,60}
\definecolor{tugreen}{RGB}{109,131,0}
\definecolor{vtred}{rgb}{0.52549,0.12157,0.25490}%
\definecolor{vtorange}{rgb}{0.89804,0.45882,0.12157}%
\definecolor{vtgrey}{rgb}{0.45882,0.47059,0.48235}%
\definecolor{tublue}{rgb}{0.00000,0.43922,0.60784}%
\definecolor{vtturquoise}{rgb}{0.17255,0.83529,0.76863}%
\definecolor{tudred}{rgb}{0.74510,0.11765,0.23529}%
\definecolor{tulgreen}{rgb}{0.54118,0.61176,0.20000}%

\newcommand{\dprod}{~\Diamond~}
\newcommand{\RR}{\mathbb{R}}

\newcommand{\cN}{\mathcal{N}}
\newcommand{\cX}{\mathcal{X}}
\newcommand{\orth}[2]{#1^{\!\perp\!{#2}}}
\newcommand{\argmin}{\operatornamewithlimits{argmin}}

\newcommand{\nc}{$\mathrm{not~converged}$}
\DeclareMathOperator{\ran}{\mathcal{R}}
\newcommand{\spann}{\mathrm{span}}

\newcommand{\myplotsvert}[1]{
    \centering
    \begin{minipage}[t]{0.15\textwidth}
        \vspace{0pt}
        \centering
        \includegraphics[scale=.75]{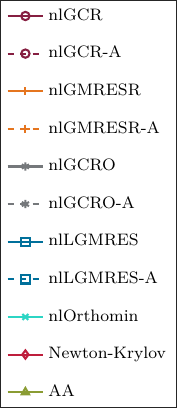}
    \end{minipage}
    \begin{minipage}[t]{0.41\textwidth}
    \vspace{0pt}
        \centering
        \includegraphics[scale=.78]{newplots/#1_iter.pdf}
    \end{minipage}
    \begin{minipage}[t]{0.41\textwidth}
    \vspace{0pt}
        \centering
        \includegraphics[scale=.78]{newplots/#1_feval.pdf}
    \end{minipage}
}

\newlist{feat}{enumerate}{1}            
\setlist[feat,1]{label=\textbf{F\arabic*}}        
\usepackage[nameinlink]{cleveref}       
\crefname{feati}{}{}     

\newsiamremark{remark}{Remark}
\newsiamremark{hypothesis}{Hypothesis}
\crefname{hypothesis}{Hypothesis}{Hypotheses}
\newsiamthm{claim}{Claim}
\newsiamremark{fact}{Fact}
\crefname{fact}{Fact}{Facts}

\headers{\emph{nlKrylov}}{T. Werner, N. Wan, and A. Mi\k{e}dlar}

\title{\emph{nlKrylov}: A Unified Framework for Nonlinear GCR-type Krylov Subspace Methods
\thanks{Submitted to the editors \today.
\funding{Work of N. Wan and A. Mi\k{e}dlar has been supported by the NSF awards DMS \#2144181 and \#2324958.}}}

\author{Tom Werner\thanks{Institute for Numerical Analysis, TU Braunschweig, Germany, \email{tom.werner@tu-braunschweig.de}, corresponding author.}
\and Ning Wan\thanks{Department of Mathematics, Virginia Tech, Blacksburg, VA, USA, \email{wning@vt.edu}, \email{amiedlar@vt.edu}.}
\and Agnieszka Mi\k{e}dlar\footnotemark[3]}

\usepackage{amsopn}

\usepackage[nopatch=footnote]{microtype}    
\usepackage[font=small]{caption}
\ifpdf
\hypersetup{
  pdftitle={\textit{nlKrylov}: A Unified Framework for Nonlinear GCR-type Krylov Subspace Methods},
  pdfauthor={T. Werner, N. Wan, and A. Mi\k{e}dlar}
}
\fi

\newcommand{\dnc}[1]{--}
\newcommand{\sep}{$\,\vert\,$}
\allowdisplaybreaks
\begin{document}

\maketitle
\begin{abstract}
In this paper, we introduce a unified framework for nonlinear Krylov subspace methods (\textit{nlKrylov}) to solve systems of nonlinear equations. Building on classical GCR-like/type linear Krylov solvers such as GMRESR, we generalize these approaches to nonlinear problems via nested algorithmic structures. We present rigorous convergence results for problems, relying on relaxed assumptions that avoid the need for exact line searches.
The framework is further extended to matrix-valued root finding problems using global nonlinear Krylov approaches. Extensive numerical experiments validate the theoretical insights and demonstrate the robustness and efficiency of our proposed algorithms.
\end{abstract}


%
%
\section{Introduction}

A classical problem in numerical analysis, now increasingly important in data science and machine learning, is the root finding problem 
\begin{equation}
    \text{Find } x\in\RR^n \ \text{ that satisfies } \ f(x)=0,\label{eq:rootfinding}
\end{equation}
where $f:\RR^n\rightarrow\RR^n$ is continuously differentiable and possibly nonlinear. Such nonlinear systems commonly arise in ODEs/PDEs solvers and unconstrained optimization problems 
\begin{equation}
\min_{x\in\RR^n} \phi(x), \ \phi:\RR^n\rightarrow\RR \ \text{ twice continuously differentiable,} \label{eq:optimization}    
\end{equation}
using the optimality condition $\nabla \phi(x)=:f(x)=0$.
Since \eqref{eq:rootfinding} frequently arises in practical applications, numerous methods have been developed and refined since the mid-20th century, most of which are based on fixed-point iteration,
\begin{equation}
    x_{j+1}=g(x_j), \quad g:\RR^n\rightarrow\RR^n. \label{eq:fixpoint}
\end{equation}
\noindent
Choosing $g(x)=x+\beta f(x)$ with $\beta\in\RR$ yields fixed points identical to the roots of $f(x)=0$. However, convergence of the fixed-point scheme \eqref{eq:fixpoint} is generally not guaranteed and can be very slow, motivating acceleration schemes such as Anderson Acceleration (AA) \cite{And1965,And2019} and Pulay's mixing method \cite{Pul1980,Pul1982}, owing to their simplicity and effectiveness. Both methods rely on sampling consecutive differences of iterates and function-related quantities, i.e., $\Delta x_j=x_{j+1}-x_j$ and $\Delta f_j=f(x_{j+1})-f(x_j)$. A second family of fixed-point schemes \eqref{eq:fixpoint} for \eqref{eq:rootfinding} is the family of Newton-type methods 
$$x_{j+1}=x_j + \Delta x_j, \ \ \text{where} \ \ B_j\Delta x_j=-f(x_j)+t_j.$$ In this setting, $B_j\approx J_f(x_j)\in\mathbb{R}^{n,n}$ is an approximation to the Jacobian $J_f(x_j)$ of $f$ at $x_j$ and $t_j\in\mathbb{R}^n$ is an error term. The classical Newton method \cite{Kelley2003} uses $B_j=J_f(x_j)$ and $t_j=0$, resulting in 
$$g(x) = x - J_f(x)^{-1}f(x).$$
In practice, the update equation is usually solved approximately with $t_j\neq 0$, leading to inexact Newton~\cite{DemES1982} or Newton--Krylov~\cite{Kelley2003} methods. Quasi-Newton variants replace $J_f(x_j)$ with iterative approximations $B_j$~\cite{Dennis1977}. Recently, nonlinear extensions of GMRES and GCR have been explored as Krylov-based accelerators~\cite{Hans2021,OosW2000,WasO1997,He2024}, which coincide with AA on linear problems~\cite{WalN2011} and connect to quasi-Newton and multi-secant updates~\cite{Eyert1996,FanS2009}.

\paragraph{Contributions and Outline}
This paper develops a unified framework for nonlinear Krylov (\textit{nlKrylov}) methods, building on the nonlinear truncated GCR (nlTGCR) method \cite{He2024}, which itself extended GCR \cite{EisES1983,HesS1952} to nonlinear systems. Central to this framework is the subroutine $\mathcal{SR}_j$, see \Cref{sec:nlkrylovframework}, which drives the construction of nonlinear search directions and nested acceleration schemes.
We emphasize that \Cref{sec:background} is purely expository, serving only to collect background material on Krylov and Newton-type methods. All original methodological developments, theoretical results, and new connections are contained in Sections \ref{sec:nlKrylov}--\ref{sec:convergence} and are as follows. A general \textit{nlKrylov} framework based on nlGCR provides a unified description of Krylov-type methods for nonlinear systems, within which extensions of classical linear solvers are derived, including GMRESR \cite{Van1994}, GCRO/GCROT \cite{DeS1996,Des1999}, and LGMRES \cite{Baker2005}, and their embedding in the proposed structure is demonstrated. Links are established to existing approaches such as quasi-Newton/multisecant methods and subspace projected Newton methods \cite{ShoT1994}. In addition, we clarify the relationship between nlTGCR and nlOrthomin, showing that the resulting \textit{nlKrylov} methods can be viewed as flexible preconditioned variants of nlOrthomin. Convergence results are proven for problems with nonsingular Jacobians under relaxed assumptions (e.g., without requiring exact line search) and for singular Jacobians via the subspace projected Newton framework \cite{ShoT1994}. The methodology further extends naturally to matrix-valued root finding problems \cite{AbsMS2007,Higham2008,ZhaYSY2020}, as illustrated in the numerical experiments.
The paper is organized as follows. \Cref{sec:background} collects background material on Krylov and Newton-type methods (expository only; no original contributions). \Cref{sec:nlKrylov} introduces the \textit{nlKrylov} framework and the main algorithmic developments. \Cref{sec:connection} establishes new connections to quasi-Newton, nlOrthomin, and projection-based methods. \Cref{sec:convergence} presents the convergence theory. Sections~\ref{sec:implementation}--\ref{sec:experiments} cover implementation details, extensions, and numerical experiments.
\paragraph{Notation}
Throughout the paper, vectors are denoted by lowercase Roman letters and matrices by uppercase Roman letters, unless stated otherwise. For a (nonlinear) function $f:\RR^n\rightarrow\RR^n$, $J_f(x)\in\RR^{n,n}$ denotes its Jacobian evaluated at a point $x\in\RR^n$. In certain contexts, we use the notation $\texttt{GMRES}(A,b,m)$ to indicate the application of the Generalized Minimal Residual method to the linear system $Ax=b$ for $m$ steps, returning the approximate solution $x_m$. Unless a better initial guess is available, GMRES is assumed to be initialized with the zero vector. The Euclidean scalar product is denoted by $\langle\cdot,\cdot\rangle$, and the associated $2$-norm by $\lVert\cdot\rVert$. By $\mathcal{R}(V)$ and $\mathcal{N}(V)$, we denote the range and Null space of a rectangular matrix $V\in\RR^{n,j}$, $j\leq n$, respectively. 
%
%
\section{Background}\label{sec:background}

\subsection{Krylov subspace methods for linear equations} \label{sec:linearkrylov}

When solving large linear systems of the form 
\begin{equation}
    Ax=b, \label{eq:linearsystem}
\end{equation} 
where the matrix $A\in\RR^{n,n}$ is large and sparse and $b \in \RR^{n}$ is the right-hand side, Krylov subspace methods are among the most widely used iterative solvers in practice~\cite{Saad2003}. These methods generate a sequence of approximations within the Krylov subspace of size $j > 0$
\begin{equation}
     K_j(A,r_0):=\spann\lbrace r_0,Ar_0,A^2r_0,\dots,A^{j-1}r_0\rbrace\subseteq\RR^n ,\label{eq:krylov} 
\end{equation}
where $r_0=b-Ax_0\in\RR^n$ is the residual for the initial guess $x_0\in\RR^n$, and each iterate satisfies $\displaystyle x_{j+1} \in x_0+K_{j+1}(A,r_0).$
It is well-known~\cite[\S 6.2]{Saad2003} that for a nonsingular $A\in\RR^{n,n}$, the Krylov subspace $K_n(A,r_0)$ of size $n$ contains an exact solution to \eqref{eq:linearsystem}. In practice, however, one typically seeks an accurate approximation $x_j$ with $j\ll n$. Among Krylov methods, this work focuses on the minimal residual~\cite[\S 2.5.5]{LieS2012} and conjugate residual methods~\cite{EisES1983}. Minimal residual methods minimize the residual norm $\|r_{j+1}\|=\|b-Ax_{j+1}\|$ over $K_{j+1}(A,r_0)$, i.e.,
$$x_{j+1} = x_0 + p_j, \quad p_j=\argmin_{p\in K_{j+1}(A,r_0)}\|b-A(x_0+p)\|.$$
The correction $p_j$ is typically determined via the associated projected problem
$$x_{j+1}=x_0+V_jy_j, \quad y_j=\argmin_{y\in\RR^{j+1}}\|b-A(x_0+V_jy)\|,$$
where $V_j=[v_0,v_1,\dots,v_j]\in\RR^{n,(j+1)}$ is an orthonormal basis for $K_{j+1}(A,r_0)$ \cite{Saad2003}. Prominent examples are MINRES~\cite{PaiS1975} for symmetric and GMRES~\cite{SaaS1986} for general matrices. A second class of Krylov subspace methods comprises variational iterative approaches derived from {CG}~\cite{HesS1952}, including GCR~\cite{EisES1983,Stiefel1955} and Orthomin($k$) \cite{Vinsome1976}. The two mentioned classes of methods fall within the general framework of (Truncated) Petrov Galerkin Krylov ((T)PGK) methods, where well-known equivalences hold in exact arithmetic, such as the equivalence of GMRES and GCR, and of Orthomin with truncated GCR \cite{SaaS1985}. Moreover, GMRES is connected to rank-one acceleration~\cite{EirN1989} and can be viewed as a Broyden-type scheme~\cite{YanG1995}. Preconditioning plays a crucial role in enhancing the efficiency of Krylov subspace methods. While constant preconditioners based on incomplete factorizations are commonly used~\cite[\S 10]{Saad2003}, variable preconditioners offer an adaptive and flexible alternative that can be adjusted dynamically during the iteration.
Examples include Flexible GMRES~\cite{Saad1993}, GMRESR~\cite{Van1994}, GCRO/GCROT~\cite{DeS1996,Des1999}, and LGMRES~\cite{Baker2005,HicZ2010}. \\
Modern applications often require solving linear systems with multiple right-hand sides, i.e.,
\begin{equation*}
    Ax_i=b_i,~i=1,\dots,p \quad \Leftrightarrow \quad AX=B,~X=[x_1,\dots,x_p]\in\RR^{n,p},~B=[b_1,\dots,b_p]\in\RR^{n,p}.
\end{equation*}
To address this, block variants of Krylov subspace methods construct iterates in the block-subspace
$$ \mathcal{K}^\square_j(A,R_0)=\mathrm{blockspan}\lbrace R_0,AR_0,A^2R_0,\dots,A^{j-1}R_0\rbrace,$$
 using a block-Arnoldi or -Lanczos procedures \cite{Ruhe1979,Saad2003,Sadkane1993}. For more general linear operator equations 
\begin{equation*}
  \mathcal{A}(X)=B,\quad X,B\in\RR^{n,p},\quad\mathcal{A}:\RR^{n,p}\rightarrow\RR^{n,p}\quad\text{linear}, 
\end{equation*}
global Krylov subspace methods were introduced in \cite{JbiMS1999} and later extended to GCRO, Flexible GMRES, QMR, LSQR, and related algorithms~\cite{MenGLF2021,TouK2006,WanG2007,ZadTW2019}. These approaches have proven effective, particularly when only operator evaluations are available, such as in Newton--Krylov frameworks~\cite{Wer2024}.

\subsection{Inexact Newton, Quasi--Newton and Newton--Krylov Methods} \label{sec:inexactnewton}

A standard approach for solving \eqref{eq:rootfinding} is Newton's method, which updates iterates via
\begin{equation}
\label{eq:Newton}
x_{j+1} = x_j + \Delta x_j, \ \mbox{ where } \ J_f(x_j)\Delta x_j = -f(x_j). 
\end{equation}
While Newton’s method converges rapidly for a good (near a solution) initial guess, computing the Jacobian
$J_f(x_j)$ can be costly or infeasible. Quasi-Newton methods improve efficiency by approximating the Jacobian or its inverse using past iterates~\cite{Dennis1977}. A prominent example is Broyden's method~\cite{Bro1965}, which replaces $J_f(x_j)$ with an approximation $B_j$ and updates it via
\begin{equation}
\label{eq:qNewton-matrix}
B_{j+1} = B_j + \frac{(\Delta f_j-B_j\Delta x_j)\Delta x_j^T}{\Delta x_j^T\Delta x_j},
\end{equation}
corresponding to a minimal Frobenius-norm update satisfying $B\Delta x_j=\Delta f_j$. Other quasi-Newton schemes, including Powell Symmetric Broyden (PSB)~\cite{Powell1970}, Broyden-Fletcher-Goldfarb-Shanno (BFGS)~\cite{Broyden1970,Fletcher1970,Goldfarb1970,Shanno1970}, and Davidon-Fletcher-Powell (DFP) method \cite{Davidon1991}, use similar updates and are locally q-superlinearly convergent~\cite{Dennis1974}.
Multisecant methods generalize these approaches by incorporating multiple recent iterates and function evaluations~\cite{GayS1978,Sch1983,FanS2009}, enabling limited-memory implementations without explicit Jacobian computations. Since quasi-Newton steps may not always provide descent directions, inexact Newton methods introduce a controlled residual error $t_j$, i.e.,
\begin{equation}
\label{eq:inexactNewton}
J_f(x_j)\Delta x_j  = -f(x_j) + t_j,  \quad \mbox{with} \quad \frac{\|t_j\|}{\|f(x_j)\|} \leq \eta_j,
\end{equation}
where $\lbrace\eta_j\rbrace_j$ governs inner accuracy~\cite{DemES1982} and $t_j\in\RR^n$ represents the resulting inexactness. Krylov-based solvers efficiently update iterates without forming $J_f(x_j)$ explicitly, leading to Jacobian--free Newton--Krylov methods~\cite{KnoK2004}.

\subsection{The Nonlinear Truncated Generalized Conjugate Residual (nlTGCR)} \label{sec:nlgcr}
In \cite{He2024}, the Nonlinear Truncated GCR (nlTGCR) method extended the linear GCR method \cite{EisES1983} to nonlinear systems. This approach closely relates to Anderson Acceleration and inexact Newton methods, offering a unified view of nonlinear iterative acceleration techniques.
Recall that GCR solves the linear system \eqref{eq:linearsystem} using the Krylov subspace \eqref{eq:krylov}, where a basis $P_j=[p_0,~p_1,\dots,~p_{j}]$ of $K_{j+1}(A,r_0)$ is constructed to satisfy $P_j^TA^TAP_j=I_{j+1}$, i.e., $P_j$ is $A^TA$-orthogonal, or equivalently, $V_j=AP_j$ is orthogonal. Orthogonality of $V_j$ is maintained using a modified Gram-Schmidt step after constructing a new direction $v_{j+1}$.
A truncated variant, GCR($k$)~\cite{EisES1983} (see \Cref{alg:GCR}) limits orthogonalization to the most recent $k$ vectors, recovering full GCR for $k=\infty$. Following~\cite{He2024}, let
\begin{equation}
j_k:=\max\lbrace j-k+1,0\rbrace,\quad
    P_j:=\left[p_{j_k},p_{j_k+1},\dots,p_j\right],\quad V_j:=\left[v_{j_k},v_{j_k+1},\dots,v_j\right], \label{eq:defV} 
\end{equation}
such that $P_j,V_j\in\RR^{n,k}$ when $j\geq k$ (truncated), and $P_j,V_j\in\RR^{n,(j+1)}$ otherwise. Since all methods considered in this paper use truncation, we will drop the "T" in nlTGCR($k$) and define $n_j:=j-j_k+1$ as the number of columns in $P_j$ (resp. $V_j$).\\
\noindent
While extending the GCR method to nonlinear systems,~\cite{He2024} specified four key features:
\begin{feat}
    \item The nonlinear version recovers the linear GCR in exact arithmetic.~\label{feat:1}
    \item The algorithm can fit into the inexact Newton or multi-secant framework.~\label{feat:2}
    \item The algorithm should exploit a more accurate linear model than Newton or Quasi-Newton methods at the cost of potential extra function evaluations. \label{feat:3}
    \item The algorithm can handle ``fuzzy'' functions in stochastic or machine learning applications. \label{feat:4}
\end{feat}
To achieve \Cref{feat:1}--\Cref{feat:4}, nlGCR($k$) assumes matrices $P_j, V_j\in\RR^{n,n_j}$ with $J_f(x_j)P_j\approx V_j$ that allow the local representation of the negative nonlinear residual $-r_{j+1}:=f(x_{j+1})$ as 
\begin{equation}
 f(x_{j+1})=f(x_j + P_jy_j)\approx f(x_j)+J_f(x_j)P_jy_j\approx f(x_j) + V_jy_j,    \label{eq:localmodel}
\end{equation}
with $y_j\in\RR^{n_j}$ selected such that the residual norm $\|r_{j+1}\|$ is minimized and
$$f(x_j)+V_jy_j \ \perp \ \ran(V_j).$$
In the case of equality in \eqref{eq:localmodel}, $y_j$ is directly determined in terms of the normal equations as 
\begin{equation}
    y_j=\argmin_{y\in\RR^{n_j}} \|f(x_j)+V_jy\|=-V_j^+f(x_j)\stackrel{(\star)}{=}V_j^Tr_j, \label{eq:gcrnormal}    
\end{equation}
where we assume the columns of $V_j$ to be orthonormal for $(\star)$.
Once $y_j$ is computed, the update $x_{j+1}=x_j+P_jy_j$ is performed and $r_{j+1}=-f(x_{j+1})$ is set. The matrices $P_j$ and $V_j$ are then extended by setting $p_{j+1}=r_{j+1}$ and $v_{j+1}=J_f(x_{j+1})p_{j+1}$ as well as orthogonalizing $v_{j+1}$ against all columns of $V_j$ and modifying $p_{j+1}$ accordingly. For linear $f(x)=Ax-b$, $J_f(x)=A$ and so
\begin{equation}
V_j=AP_j=J_f(x_j)P_j, \quad P_jy_j=P_jV_j^Tr_j=P_j(v_j^Tr_j)e_{n_j}=(v_j^Tr_j)p_j\equiv \alpha_j p_j, \label{eq:gcrPalphap} 
\end{equation}
due to the orthogonality condition $r_j\perp\ran(V_{j-1})$ in GCR, where $e_{n_j}\in\RR^{n_j}$ is the $n_j$-th unit vector. However, equations \eqref{eq:gcrnormal} and \eqref{eq:gcrPalphap} are not true in the nonlinear version since, generally, $J_f(x_{j+1})\neq J_f(x_j)$ and, as such, 
$$V_j = \left[J_f(x_{j_k})p_{j_k},J_f(x_{j_k+1})p_{j_k+1},\dots,J_f(x_j)p_j\right]\neq J_f(x_j)\left[p_{j_k},p_{j_k+1},\dots,p_j\right]=J_f(x_j)P_j.$$
A direct comparison between GCR for linear equations and its nonlinear counterpart is presented in \Cref{alg:GCR} and \Cref{alg:nlGCR}, where the notation from \eqref{eq:defV} is used.\\
\begin{minipage}{.495\textwidth}
\begin{algorithm}[H]
    \caption{GCR($k$) for $Ax=b$  \ \cite{DeS1996,Van1994}}    
    \label{alg:GCR}
    \begin{algorithmic}[1]
    \setstretch{1.05}
    \REQUIRE $A\in\mathbb{R}^{n,n}$, $k\in\mathbb{N}$, $x_0,b\in\mathbb{R}^{n}$
    \ENSURE $x^*$ approximate solution to $Ax=b$
    \STATE $\widehat{p}=r_0=b-Ax_0$,\quad $\widehat{v}=A\widehat{p}$
    \STATE $p_0=\frac{\widehat{p}}{\lVert \widehat{v}\rVert}$,\quad $v_0=\frac{\widehat{v}}{\lVert \widehat{v}\rVert}$
    \STATE $j=0$
    \WHILE{\nc}
        \STATE $\alpha_j=\langle v_j,r_j\rangle$  \label{line:inner1}
        \STATE $x_{j+1} = x_j + \alpha_j p_j$,\quad $r_{j+1} = r_j - \alpha_j v_j$
        \STATE $\widehat{p} = r_{j+1}$, $\widehat{v}=A\widehat{p}$
        \FOR{$i=j_k:j$}
            \STATE $\beta_i=\langle \widehat{v},v_i\rangle$ \label{line:inner2}
            \STATE $\widehat{p} = \widehat{p} - \beta_ip_i$,\quad $\widehat{v}=\widehat{v} - \beta_iv_i$
        \ENDFOR
        \STATE $p_{j+1}=\frac{\widehat{p}}{\lVert \widehat{v}\rVert}$,\quad $v_{j+1}=\frac{\widehat{v}}{\lVert \widehat{v}\rVert}$
        \STATE $j=j+1$
    \ENDWHILE
    \RETURN $x^*=x_j$
    \end{algorithmic}
\end{algorithm}
\end{minipage}
\begin{minipage}{.495\textwidth}
\begin{algorithm}[H]
    \caption{nlGCR($k$) for $f(x)=0$  \ \cite{He2024}}   
    \label{alg:nlGCR}
    \begin{algorithmic}[1]
    \setstretch{1.05}
    \REQUIRE $x_0\in\mathbb{R}^{n}$, $k\in\mathbb{N}$, $f,J_f$
    \ENSURE{$x^*$ approximate solution to $f(x)=0$}
    \STATE $\widehat{p}=r_0=-f(x_0)$,\quad $\widehat{v}=J_f(x_0)\widehat{p}$
    \STATE $p_0=\frac{\widehat{p}}{\lVert \widehat{v}\rVert}$,\quad $v_0=\frac{\widehat{v}}{\lVert \widehat{v}\rVert}$
    \STATE $j=0$
    \WHILE{\nc}
        \STATE $y_j={V_j}^Tr_j$ \label{line:VTr}
        \STATE $x_{j+1} = x_j +  P_jy_j$,\quad $r_{j+1} = -f(x_{j+1})$ \label{line:nlupdate}
        \STATE $\widehat{p} = r_{j+1}$, $\widehat{v}=J_f(x_{j+1})\widehat{p}$ \label{line:startgs}
        \FOR{$i=j_k:j$}
            \STATE $\beta_i=\langle \widehat{v},v_i\rangle$
            \STATE $\widehat{p} = \widehat{p} - \beta_ip_i$,\quad $\widehat{v}=\widehat{v} - \beta_iv_i$
        \ENDFOR
        \STATE $p_{j+1}=\frac{\widehat{p}}{\lVert \widehat{v}\rVert}$,\quad $v_{j+1}=\frac{\widehat{v}}{\lVert \widehat{v}\rVert}$ \label{line:endgs}
        \STATE $j=j+1$
    \ENDWHILE
    \RETURN $x^*=x_j$
\end{algorithmic}
\end{algorithm}
\end{minipage}
%
%
\section{From Linear to Nonlinear Krylov methods} \label{sec:nlKrylov}
In this section, we will derive a unified framework for generalizing linear nested Krylov methods based on GCR to nonlinear equations and present three notable members of the resulting class of methods. Following
the terminology of \cite{He2024}, which introduces nlGCR as the nonlinear analogue of GCR, we refer to methods that fit within our framework as \emph{nonlinear Krylov} (\textit{nlKrylov}) methods.
We first revisit the ideas from the original nlTGCR paper presented in \autoref{sec:nlgcr} to develop a more general framework for this family of methods in \autoref{sec:nlkrylovframework}. Afterward, we will specify how to apply the introduced framework to three popular nested GCR-based algorithms, namely GMRESR, GCRO, and LGMRES in \autoref{sec:nestednlkrylov}, obtaining nonlinear extensions of these methods. 
\subsection{The general nlKrylov framework} \label{sec:nlkrylovframework}
In establishing the framework for \textit{nlKrylov} methods, we recall that at step $j$ of the nlGCR algorithm, the matrices $P_j$ and $V_j$ are employed as sets of directions for progressing the iteration, as shown in lines \ref{line:VTr} and \ref{line:nlupdate} of \Cref{alg:nlGCR}. Since in the derivation of nlGCR, the basis matrices $P_j,V_j$ seem rather arbitrary, one might also think of improving the local linear model \eqref{eq:localmodel} at the cost of some additional function evaluations, as indicated by \Cref{feat:3}. Inspired by the idea of nested GCR-type methods \cite{DeS1996}, our aim is to compute $y_j$, $x_j$, and $r_j$ as in nlGCR, but with $p_j$ defined via a specified subroutine $\mathcal{SR}_j(r_j,J_f(x_j))$, which will involve the current residual $r_j$ as well as applications of the Jacobian $J_f(x_j)$. In the simplest case, where $\mathcal{SR}_j(r_j,J_f(x_j))=r_j$, the standard nlGCR algorithm is recovered.
Building on the concept of nested Krylov methods, $\mathcal{SR}_j$ represents a specific algorithm for solving the linear equation
\begin{equation}
    J_f(x_j)\widehat{p}=r_j,\label{eq:refinedlocal}    
\end{equation}
which requires only $r_j$ as well as applications of $J_f(x_j)$ in the process. We will refer to every algorithm that uses nlGCR as an outer iteration and is equipped with a subroutine $\mathcal{SR}_j$ for \eqref{eq:refinedlocal} as a \emph{nonlinear Krylov method}. Observe that performing a few iterations of an iterative inner solver for \eqref{eq:refinedlocal} effectively incorporates the inexact Newton update direction $\Delta x_j$ 
from \eqref{eq:inexactNewton} into the search space $P_j$, thereby improving the local linear model. However, as a large number of applications of $J_f(x_j)$ within $\mathcal{SR}_j$ may introduce substantial computational overhead, their use should therefore be carefully balanced to enhance the local linear model \eqref{eq:localmodel} as effectively as possible without dominating the overall complexity.
In \cite{DeS1996}, classical inner methods include Krylov subspace methods such as GMRES or the Bi-Conjugate Gradient Stabilized method (BiCGStab, \cite{Van1992}) if $J_f(x_j)$ is non-symmetric and CG for symmetric Jacobians. However, our framework is not restricted to Krylov subspace methods to obtain $p_j$. Different techniques such as classical iterative methods, randomized projection approaches or whatever algorithm is viable for subproblem \eqref{eq:refinedlocal} can be used. Nevertheless, in this paper, we will focus on three particular Krylov subspace approaches that will lead to nonlinear extensions of GMRESR($m,k$) , GCRO($m,k$) and LGMRES($m,k$). Before getting into details, the general form of \textit{nlKrylov} methods is displayed in \Cref{alg:nlKrylovSR}, where again, notation from \eqref{eq:defV} is used.
\begin{algorithm}[H]    
    \caption{A template for nlKrylov($k$) methods to solve $f(x)=0$}
    \label{alg:nlKrylovSR}
    \begin{algorithmic}[1]
    \REQUIRE{$x_0\in\mathbb{R}^{n}$, $k\in\mathbb{N}$, $f:\RR^n\rightarrow\RR^n$, $J_f:\RR^n\rightarrow\RR^{n,n}$, $\mathcal{SR}:\RR^n\times\RR^{n,n}\rightarrow\RR^n$}
    \ENSURE{$x^*$ approximate solution to $f(x)=0$}
    \STATE $r_0=-f(x_0)$
    \STATE $\widehat{p}=\mathcal{SR}_0(r_0,J_f(x_0))$,\quad $\widehat{v}=J_f(x_0)\widehat{p}$
    \STATE $p_0=\frac{\widehat{p}}{\lVert \widehat{v}\rVert}$,\quad $v_0=\frac{\widehat{v}}{\lVert \widehat{v}\rVert}$,\quad $j=0$
    \WHILE{\nc}
        \STATE $y_j = V_j^Tr_j$ 
        \STATE $x_{j+1} = x_j + P_jy_j$,\quad $r_{j+1} = -f(x_{j+1})$
        \STATE $\widehat{p} = \mathcal{SR}_{j+1}(r_{j+1},J_f(x_{j+1}))$ \label{line:SRsolve}
        \STATE $\widehat{v}=J_f(x_{j+1})\widehat{p}$ \label{line:SRupdate}
        \FOR{$i=j_k:j$}
            \STATE $\beta_i=\langle \widehat{v},v_i\rangle$,\quad $\widehat{p} = \widehat{p} - \beta_ip_i$,\quad $\widehat{v}=\widehat{v} - \beta_iv_i$ \label{line:nlKrylovortho}
        \ENDFOR
        \STATE $p_{j+1}=\frac{\widehat{p}}{\lVert \widehat{v}\rVert}$,\quad $v_{j+1}=\frac{\widehat{v}}{\lVert \widehat{v}\rVert}$,\quad $j=j+1$
    \ENDWHILE
    \RETURN $x^*=x_j$
\end{algorithmic}
\end{algorithm}
\subsection{Obtaining nonlinear extensions of established nested GCR-type methods}\label{sec:nestednlkrylov} 
In \cite{Van1994,DeS1996} and later \cite{Baker2005}, various techniques have been suggested to improve GCR for linear problems, leading to GMRESR, GCRO and LGMRES. Here, we want to demonstrate how these concepts apply to nonlinear problems. 
\subsubsection{Solving by GMRES: nlGMRESR} \label{sec:nlgmresr}
As mentioned before, an intuitive way to choose $\mathcal{SR}_j$ would be to apply $m$ steps of an iterative solver to \eqref{eq:refinedlocal}. For general, potentially large-scale problems, we want to avoid an explicit computation of the Jacobian, and its transpose is generally not available if we use finite difference approximations. In \cite{Van1994}, GMRESR($m$) was originally introduced as a nested GCR algorithm, where $m$ steps of GMRES for the linear system $A\widehat{p}=r_j$ are used as an inner method. This approach easily extends to the nonlinear setting by choosing 
\begin{equation}
\mathcal{SR}_{j+1}(r_{j+1},J_f(x_{j+1}))=\texttt{GMRES}(J_f(x_{j+1}),r_{j+1},m)  \label{eq:SRgmresr}
\end{equation}
as a subroutine in line \ref{line:SRsolve} of \Cref{alg:nlKrylovSR}. We will call the resulting \textit{nlKrylov} method nlGMRESR($m,k$) and expect it to require fewer iterations than nlGCR($k$) due to its enhanced local linear model, while the additional function evaluations introduced by the inner solve \eqref{eq:SRgmresr} can be controlled through an appropriate choice of $m$. For linear systems, a detailed discussion on the choice of $m$ can be found in \cite[\S 5.2]{Van1994}. A comprehensive analysis for nlGMRESR($m,k$) is still an open question. In the numerical experiments presented in \Cref{sec:experiments}, we select $m$ heuristically, typically between $m=4$ and $m=30$ depending on the problem; a brief intuition for this choice is given in Section \ref{sec:config}.
Note that, by choosing a zero initial guess in GMRES, we can recover $\widehat{p}$ in line \ref{line:SRsolve} and $\widehat{v}$ in line \ref{line:SRupdate} of \Cref{alg:nlKrylovSR} from the Arnoldi relation as a low-rank update without the need of additional function evaluations. Thus, one step of nlGMRESR($m,k$) requires $(m-1)$ additional function evaluations for \eqref{eq:SRgmresr} compared to nlGCR($k$).
\subsubsection{Taking the nonlinear basis into account: nlGCRO and nlLGMRES}
For linear problems, it has been noted in \cite{DeS1996} and later \cite{Baker2005} that a major limitation of GMRESR($m$) is, that it treats the subproblem \eqref{eq:refinedlocal} in isolation, without leveraging information from the outer GCR loop to improve the GMRES-solve. This is a common drawback in nested and flexible iterative methods, which can be overcome by techniques like augmentation, deflation, and recycling, which are all closely related for linear problems \cite{SoodDK2020}. As we will see in a moment and in numerical experiments in \autoref{sec:experiments}, the behavior of deflated and augmented nonlinear methods can be very different.\\
The first method to take into account the outer basis $V_j$ in the inner solve is the GCRO($m$) algorithm \cite{DeS1996}, which aims at deflating information contained in $V_j$ by solving the projected problem 
\begin{equation}
    (I-V_jV_j^T)A\widehat{p}=:\orth{A}{V_j}\widehat{p}=\widetilde{r}_{j+1}:=(I-V_jV_j^T)r_{j+1}, \label{eq:gcrosolve}
\end{equation} 
by GMRES instead. Note that for linear GCR, we have $r_{j+1}\perp \ran(V_j)$ and, as such, $\widetilde{r}_{j+1}=r_{j+1}$. Also, the application of $(I-V_jV_j^T)$ to $\widehat{v}=A\widehat{p}$ orthogonalizes $\widehat{v}$ against the outer basis $V_j$, which corresponds to the Gram-Schmidt process in line \ref{line:nlKrylovortho} of \Cref{alg:nlKrylovSR}, i.e., GCRO($m$) is respecting the progress of the outer method in the inner solve. Additionally, assuming that a zero initial guess is used and 
\begin{equation}
\orth{A}{V_j}Q_m=Q_{m+1}\underline{H}_m,\qquad \gamma_m=\argmin_{\gamma\in\RR^m} \Big\|\|r_{j+1}\|e_1-\underline{H}_m \gamma\Big\| \label{eq:gcroarnoldi}    
\end{equation}
hold inside GMRES, we get 
\begin{align}
\|v_{j+1}\|v_{j+1}&=(I-V_jV_j^T)AQ_m\gamma_m = \orth{A}{V_j}Q_m\gamma_m= Q_{m+1}(\underline{H}_m\gamma_m),\label{eq:gcrov}\\ 
\|v_{j+1}\|p_{j+1}&=A^{-1}\orth{A}{V_j}Q_m\gamma_m=(I-A^{-1}V_jV_j^TA)Q_m\gamma_m=(Q_m-P_jB_m)\gamma_m, \label{eq:gcrop}
\end{align}
where $B_m = V_j^TAQ_m\in\RR^{n_j,m}$, and $p_{j+1}$ and $v_{j+1}$ can be recovered from a low-rank update without an additional function evaluation and without the need for orthogonalization in line \ref{line:nlKrylovortho} of \Cref{alg:nlKrylovSR}.
Moving to the nonlinear case, we generally do not have $J_f(x_j) P_j = V_j$, which means that the information in $V_j$ is not solely associated with $J_f(x_j)$ but rather reflects contributions from all previously used Jacobians. As such, including a projection onto $\ran(V_j)^\perp$ in the linear solve is only beneficial when the Jacobian is changing slowly, i.e., when $J_f(x_j)\approx J_f(x_{j-1})$. If the subspace $\ran(V_j)$ becomes outdated relative to the subproblem \eqref{eq:refinedlocal}, projecting $J_f(x_{j+1})$ onto $\ran(V_j)^\perp$ may fail to improve the local linear model and can even hinder convergence. Additionally, we generally have $\widetilde{r}_{j+1}\neq r_{j+1}$ in \eqref{eq:gcrosolve} since $r_{j+1} \not\perp \ran(V_j)$. Nevertheless, using 
\begin{equation}
\mathcal{SR}_{j+1}(r_{j+1},J_f(x_{j+1}))=\texttt{GMRES}(\orth{J_f(x_{j+1})}{V_j},\widetilde{r}_{j+1},m)  \label{eq:SRgcro}
\end{equation}
as subroutine allows to solve the deflated problem \eqref{eq:gcroarnoldi} and update $p_{j+1}$ and $v_{j+1}$ through \eqref{eq:gcrop} and \eqref{eq:gcrov}, where we have to use $B_m = V_j^TJ_f(x_{j+1})Q_m\in\RR^{n_j,m}$. The additional cost of nlGCRO($m,k$) turns out to be $m-1$ function evaluations as well as $m\cdot n_j$ inner products compared to nlGCR($k$), if we skip the outer orthogonalization, which might decrease stability of the algorithm.\\
A second way to incorporate the outer basis in the inner solve is to augment the linear Krylov subspace used by GMRES. This technique was used in the Loose GMRES (LGMRES) algorithm \cite{Baker2005,HicZ2010} to accelerate the convergence of GMRES. That is, the augmented solver minimizes the residual over the space $\mathcal{K}_{m_j}(A,r_{j+1})\oplus \ran(P_j)$ through the Arnoldi-like relation 
\begin{equation}
    A [Q_{m_j},P_j] = Z_{m+k}\underline{H}_{m+k}, \label{eq:augment} 
\end{equation}
where $m_j=m+k-n_j$ such that $[Q_{m_j},P_j]$ always has $m+k$ columns. For linear GCR($k$), we have $AP_j=V_j$ and therefore $Z_{m+k}=[Q_{m_j+1},V_j]$, which leads to $m_j$ matrix-vector products required for the inner solve. For nonlinear problems, we generally have $J_f(x_{j+1})P_j\neq V_j$ and \eqref{eq:augment} should be used, leading to $m+k$ additional function evaluations compared to nlGCR($k$). The augmented GMRES (\texttt{AGMRES}) algorithm to solve \eqref{eq:augment} is stated in the supplementary material. Using 
\begin{equation}
\mathcal{SR}_{j+1}(r_{j+1},J_f(x_{j+1}))=\texttt{AGMRES}(J_f(x_{j+1}),r_{j+1},m,k,P_j)  \label{eq:SRlgmres}
\end{equation}
in \Cref{alg:nlKrylovSR}, we get a nonlinear extension of LGMRES, which we will call nlLGMRES($m,k$). Note that the original LGMRES algorithm introduced in \cite{Baker2005} does not enforce orthogonality of $V_j$ nor $P_j$. As such, the nlLGMRES($m,k$) algorithm obtained by using \eqref{eq:SRlgmres} in \Cref{alg:nlKrylovSR} violates \ref{feat:1} in the sense that the nonlinear version does not reduce to the original algorithm when applied to linear problems. However, the idea integrates naturally into the \textit{nlKrylov} framework and is able to improve the performance of nlGMRESR($m,k$) in some of the numerical experiments.
To conclude this section, we want to emphasize the eminent need for truncation in those algorithms that consider the outer space in the subroutine $\mathcal{SR}_j$. If no truncation is employed, the additional computational cost due to orthogonalization, in the case of nlGCRO, or function evaluations, in the case of nlLGMRES, will dominate the total cost of the linear solve with minor contributions to the overall convergence of the nonlinear scheme. Nevertheless, we also observed that, in the cases where taking the outer space into account improved the convergence, choosing $k$ around ten lead to faster convergence than choosing $k$ around one or two.  
\subsubsection{Comparing nlGCR, nlGMRESR, nlGCRO and nlLGMRES}
In this section, we want to compare the three methods introduced in \Cref{sec:nlKrylov} to nlGCR concerning the underlying spaces used to update $\widehat{p}$ in line \ref{line:SRsolve} of \Cref{alg:nlKrylovSR} as well as the number of additional function evaluations and inner products required for this refinement of the local linear model. A summary is presented in \Cref{tab:nlkrylov_update}, where we assume that matrix-free products in the inner solve are evaluated via the finite difference formula
\begin{equation}
    J_f(x_j,z)=\frac{f(x_j+\varepsilon z)-f(x_j)}{\varepsilon}=\frac{f(x_j+\varepsilon z)+r_j}{\varepsilon},\label{eq:matrixfreederivative}
\end{equation}
for a small step-size $\varepsilon>0$, at the cost of one extra function evaluation and $n_j\equiv k$. We can see that the total number of function evaluations for nlGMRESR($m,k$) and nlGCRO($m,k$) is the same, while an additional $k\cdot m$ inner products are needed for the projection onto $V_j^\perp$. nlLGMRES($m,k$) requires $k$ additional function evaluations compared to nlGMRESR($m,k$) and nlGCRO($m,k$) due to the explicit augmentation. This may be avoided if the linearization $V_j=J_f(x_{j+1})P_j$ is used. Due to the larger subspace of dimension $m+k$, the number of inner products in nlLGMRES($m,k$) is the same as in nlGMRESR($m+k,k$). We want to point out again that the subspaces used by nlGCRO($m,k$) and nlLGMRES($m,k$) are generally not the same for nonlinear problems.    
\begin{table}[H]
    \centering
    \begin{tabular}{|l|l|l|l|}
        \hline
        Method & Space for $\widehat{p}$ in step $j$ & $\#$ fevals & $\#$ inner products\\
        \hline\hline
         nlGCR($k$) & $\spann\lbrace r_{j+1}\rbrace$ & $1$ & $0$ \\ 
         \hline
         nlGMRESR($m,k$)& $\mathcal{K}_m((J_f(x_{j+1}),r_{j+1})$ & $m-1$ & $\frac{m(m+1)}{2}+1$\\
         \hline 
         nlGCRO($m,k$) & $\ran(P_j)\oplus \mathcal{K}_m(\orth{J_f(x_{j+1})}{V_j},\widetilde{r}_{j+1})$ & $m-1$& $\frac{(m+1)(m+2k)}{2}+1$\\
         \hline
        nlLGMRES($m,k$) & $\ran(P_j) \oplus\mathcal{K}_m(J_f(x_{j+1}),r_{j+1})$& $m+k-1$& $\frac{(m+k)(m+k+1)}{2}+1$\\
        \hline
    \end{tabular}
    \caption{A comparison of selected \textit{nlKrylov}($k$) methods.}
    \label{tab:nlkrylov_update}
\end{table}
%
%
%
\section{Connection of nlKrylov methods to existing nonlinear methods} \label{sec:connection}
\subsection{Local updating and connection to quasi-Newton methods}
Recall that in the quasi-Newton framework, the classical Newton-update \eqref{eq:Newton} is replaced by its quasi-Newton counterpart, i.e.,
$$x_{j+1}=x_j-G_j f(x_j),$$
where $G_j\in\RR^{n,n}$ is an approximation of the inverse Jacobian $G_j\approx J_f(x_j)^{-1}$. Among such methods, the notable Broyden’s scheme \cite{Bro1965} updates $G_j$ iteratively via a low-rank correction similar to \eqref{eq:qNewton-matrix}. Within the family of  \textit{nlKrylov} algorithms, at step  $j\geq 1$,  all methods perform the update
\begin{equation}
    x_{j+1}=x_j+P_j{V_j}^Tr_j=x_j+P_j{V_j}^T(-f(x_j))= x_j-P_j{V_j}^Tf_j, \qquad r_{j+1}=-f(x_{j+1}),   \label{eq:nlKrylov_update}
\end{equation}
with the sole difference between the methods arising from how the $P_j$ and $V_j$ are selected via the corresponding local linear model \eqref{eq:localmodel} and the subroutine $\mathcal{SR}_j$. From \eqref{eq:nlKrylov_update}, we can immediately see that any \textit{nlKrylov} method can be viewed as a quasi-Newton method with $G_j\approx P_jV_j^T$ being used as a rank-$n_j$-approximation to the inverse Jacobian $J_f(x_j)^{-1}$. Although this has already been observed for nlGCR in \cite{He2024}, this result in fact extends to all \textit{nlKrylov} methods.
Note that in \cite{He2024}, a connection between nlTGCR($k$) and Anderson Acceleration with truncation (AA($k$), \cite{And1965}) has been established through the framework of multi-secant methods. Thus, the \textit{nlKrylov} family can be regarded as nonlinear acceleration methods for fixed-point schemes.

\subsection{Connection to nonlinear Orthomin} \label{sec:orthomin}
It is well known that, for linear problems, the truncated version of the GCR method, denoted GCR($k$), is equivalent to Orthomin($k$) \cite{EisES1983}. Extensions of
the classical Orthomin \cite{Vinsome1976} to nonlinear problems have been first proposed in the form of nonlinear Orthomin(1) \cite{Chr1992} and later generalized to nonlinear Orthomin($k$) in \cite{CheC2001} for $k\geq 2$.  In \Cref{alg:nlorthomin}, we recall the 
nonlinear Orthomin($k$) algorithm, with notation consistent with the conventions of this paper.
\begin{algorithm}[ht]
    \caption{nlOrthomin($k$) for $f(x)=0$ \quad \cite[Alg. 3]{CheC2001}}
    \label{alg:nlorthomin}
\begin{algorithmic}[1]
    \REQUIRE{$x_0\in\RR^n$, $k\in\mathbb{N}$, $f,J_f$}
    \ENSURE{$x^*$ approximate solution to $f(x)=0$}
    \STATE $p_0=r_0=-f(x_0)$
    \STATE $j=0$
    \WHILE{\nc}
    \STATE Solve $y_j= \mathrm{argmin}_{y\in\RR^{n_j}}\|f(x_j+P_jy)\|$ (stagnation for $y_j=0$) \label{line:linesearchorthomin}
    \STATE $x_{j+1}=x_j+P_jy_j$
    \STATE $r_{j+1}=-f(x_{j+1})$
    \STATE $p_{j+1}=r_{j+1}-\sum_{i=j_k}^{j} \beta_i^{(j)} p_i$, ~where~ $\beta_i^{(j)}=\frac{\left(J_f(x_{j+1})r_{j+1}\right)^T \left(J_f(x_{j+1})p_i)\right)}{\|J_f(x_{j+1})p_i\|^2}$ \label{line:gramschmidt}
    \STATE $j=j+1$
    \ENDWHILE
    \RETURN $x^*=x_{j}$   
\end{algorithmic}
\end{algorithm}
Two significant differences emerge when comparing \Cref{alg:nlorthomin} and \Cref{alg:nlGCR}. Firstly, in line \ref{line:linesearchorthomin} of nlOrthomin($k$), the vector $y_j$ is the exact minimizer of $\|f(x_j+d_j)\|$ along $d_j \in \ran(P_j)$. In nlGCR($k$), $y_j=V_{j}^Tr_{j}$ is the minimizer of the linearized problem \eqref{eq:localmodel}, which involves a first order Taylor approximation of $f(x_j+P_jy_j)$ as well as the local approximation $J_f(x_j)P_j\approx V_j$. If these approximations are \emph{sufficiently} accurate, the nlGCR update is essentially equivalent to a Gauss-Newton step \cite{Deuflhard2011} for the exact minimization problem
\begin{equation}
 \min_y \|s_j(y)\|,\quad \text{where}\quad s_j:\RR^{n_j}\rightarrow\RR^n, \ y\mapsto f(x_j+P_jy).\label{eq:nlmin}
\end{equation}
Secondly, in nlGCR, the matrix $V_j$ is stored as an additional set of vectors and is successively orthogonalized as each new pair ($p_j,v_j$) is computed, ensuring $V_j^TV_j=I_{n_j}$ and $V_j\approx J_f(x_j)P_j$. In contrast, in nlOrthomin($k$), $V_j=J_f(x_j)P_j$ holds exactly and is recomputed at each iteration using the current Jacobian applied to $P_j$. The orthogonalization of $V_j$ is performed implicitly using the coefficients $\beta_i^{(j)}$ before $V_j$ is discarded, which corresponds to the classical Gram-Schmidt process in line \ref{line:gramschmidt} of \Cref{alg:nlorthomin}. The similarity between
nlGCR($k$) and nlOrthomin($k$) primarly depends on the quality of the local linearization of $f(x_j+P_jy_j)$, i.e., on how slowly the Jacobian varies during the iteration, which-- following~\cite{He2024}--can be measured by the error matrix
\begin{equation}
W_j=J_f(x_j)P_j-V_j.  \label{eq:errormatrix}  
\end{equation}
Hence, if $\|W_j\|$ is small, nlGCR($k$) and nlOrthomin($k$) should exhibit comparable behavior.
In \cite{CheC2001}, it is also shown that preconditioning nlOrthomin($k$) with a constant preconditioner \textbf{M}, i.e.,
$$0=g(z)=f(\textbf{M}z),\quad z=\textbf{M}^{-1}x,$$
can significantly improve convergence. In the context of \textit{nlKrylov} methods, the subroutine $\mathcal{SR}_j$ can serve as a flexible preconditioner $\textbf{M}_j$ within nlOrthomin($k$), thus linking flexible preconditioned nlOrthomin($k$) to the family of \textit{nlKrylov} methods.
Upon comparing \textit{nlKrylov} methods to nlOrthomin, one main difference in behavior was observed. The convergence, measured by the number of iterations, was comparable between nlOrthomin and nlGCR, despite the fact that nlOrthomin employed an exact line-search, whereas nlGCR relied on the linear model. However, because the exact line-search involves an inner optimization, the number of function evaluations of $f$ was significantly higher for nlOrthomin than for most other \textit{nlKrylov} methods.

\subsection{Connection to subspace projection methods}\label{sec:paa}
In \cite{ShoT1994}, a \emph{Preconditioned Subspace Projection} approach is proposed to accelerate the convergence of Newton's method, which can be summarized as the simple three step~\Cref{alg:ppnm}.
\begin{algorithm}[ht]
    \caption{Preconditioned Projection Newton method}    
    \label{alg:ppnm}  
\begin{algorithmic}[1]
    \REQUIRE{$x_0\in\RR^n$, $k\in\mathbb{N}$, $f,J_f$}
    \ENSURE{$x^*$ approximate solution to $f(x)=0$}
    \STATE $j=0$
    \WHILE{\nc}
       \STATE Choose a \emph{Preconditioner} $Y_j\in\mathbb{R}^{n,k}$ and form the matrix $J_j = J_f(x_j)Y_j$
       \STATE Solve $J_jy_j=-f(x_j)$ \ for \ $y_j\in\mathbb{R}^k$, e.g. using $y_j=-J_j^+f(x_j)$ \label{line:paay}
       \STATE Update $x_{j+1}=x_j+Y_jy_j$
       \STATE $j=j+1$
    \ENDWHILE
    \RETURN $x^*=x_j$
\end{algorithmic}
\end{algorithm}
Note that the term \emph{Preconditioner} in the context of \Cref{alg:ppnm} should not be confused with the preconditioner \textbf{M} used in \autoref{sec:orthomin} for Krylov subspace methods. In \Cref{alg:ppnm}, $Y_j$ can be understood as a set of carefully selected directions that ``sketch'' the action of the current Jacobian $J_f(x_j)$, enabling the Newton update equation to be solved within this subspace. If we assume that locally, the Jacobian $J_f(x_j)$ does not change \emph{too much} (see \Cref{ass:upperbound} for details) and we choose $Y_j=P_j$ in iteration $j$ (assuming $n_j=k$), then
$$J_j=J_f(x_j)P_j\approx V_j.$$ 
Hence, assuming $V_j$ is orthogonal and $r_j=-f(x_j)$ is the nonlinear residual, the optimal mixing parameter $y_j\in\RR^k$ in line \ref{line:paay} of \Cref{alg:ppnm} can be computed as
$$y_j=V_j^+(-f(x_j))=V_j^T r_j.$$ Note that this is exactly the update \eqref{eq:gcrnormal}
within the \textit{nlKrylov} framework.
In \cite{ShoT1994}, \Cref{alg:ppnm} is called a \emph{Preconditioned Angle Algorithm} (PAA), if at each step, some column $z\in\RR^n$ of $J_j$ satisfies the angle criterion 
\begin{equation}
  \cos^2(\theta)=\frac{(r_j^Tz)^2}{\|r_j\|^2\cdot\|z\|^2}>\tau,\quad \tau\in(0,1), \label{eq:paa}  
\end{equation}
with angle $\theta=\angle(z,r_j)$ and tolerance $\tau$. Local convergence of PAA is ensured by 
\cite[Theorem 5]{ShoT1994}, under the assumption that $J_f(x^*)$ is invertible at a solution $x^*\in\RR^n$ to $f(x)=0$. The proof frames PAA as an inexact Newton method, measuring the inexactness by means of $Y_j$ and $J_j$ as well as \eqref{eq:paa}. In the following section, we present a comprehensive convergence analysis of \textit{nlKrylov} methods, combining both theoretical insights and practical implications, and employing a similar approach based on $P_j$ and $V_j$.
%
%
\section{Convergence analysis of nlKrylov methods} \label{sec:convergence}
The goal of this section is to examine the convergence behavior of \textit{nlKrylov} methods, beginning with problems that have nonsingular Jacobians at the solution and then addressing those with singular Jacobians.

\subsection{Problems with nonsingular Jacobian}
\label{subsec:convnonsingular}
Suppose that at iteration $j$, the spaces $P_j,V_j\in\RR^{n,n_j}$ have been constructed, with $V_j^TV_j=I_{n_j}$. Let $W_j\in\RR^{n,n_j}$ denote the error matrix defined in \eqref{eq:errormatrix}, such that
\begin{equation}
    W_j=J_f(x_j)P_j-V_j~\Leftrightarrow~V_j=J_f(x_j)P_j-W_j. \label{eq:errormatrix2}
\end{equation}
Assuming $J_f(x_j)$ is locally invertible in the neighborhood of $x^*$--and in particular invertible at $x_j$-- we set 
$$Y_j=P_j - J_f(x_j)^{-1}W_j\in\RR^{n,n_j},$$ such that 
$$J_j = J_f(x_j)Y_j=J_f(x_j)P_j - W_j = V_j\in\RR^{n,n_j}.$$
Thus, we can compute $y_j=-J_j^+f(x_j)=-V_j^Tf(x_j),$ which is 
precisely how it is computed in the \textit{nlKrylov} framework. This relation provides a natural interpretation of \textit{nlKrylov} methods through the lens of preconditioned projection methods, where the dependence on the underlying subspaces is explicit, and these subspaces are dictated by the choice of the subroutine $\mathcal{SR}_j$. We use this perspective to further analyze the convergence behavior of \textit{nlKrylov} methods.
In \cite[Thm. 4.5]{He2024}, the convergence of nlGCR was established under the assumption of an exact line search, smoothness of $f(x)$, the existence of uniform bounds $\eta,\mu\in[0,1)$ for the local linear model 
\begin{equation}
    \|W_jy_j\|\leq \mu\|f(x_j)\|, \label{eq:asslinres}
\end{equation}
and a sufficient decrease of the linearized residual 
\begin{equation}
    \|f(x_j)+V_jy_j\|\leq\eta\|f(x_j)\|. \label{eq:assdecrease}
\end{equation}
As we will see, in our analysis, we do not require the first two assumptions and can relax the uniform bounds \eqref{eq:asslinres} and \eqref{eq:assdecrease} as follows.
\begin{assumption}\label{ass:upperbound}
    Let $\displaystyle\mu_j:=\frac{\|W_jy_j\|}{\|f(x_j)\|}$ and $\displaystyle\eta_j:=\frac{\|f(x_j)+V_jy_j\|}{\|f(x_j)\|}$. Then, for all $j$ there exists a constant $c>0$ such that 
    \begin{equation}
        \mu_j+\eta_j=:c_j\leq c<1, \label{eq:upperbound}  
    \end{equation}
\end{assumption}
Note that, within the PAA framework, $\eta_j$ can be interpreted as the angle between $f(x_j)$ and $\ran(V_j)$; see \eqref{eq:paa}. The convergence of the \textit{nlKrylov} methods can now be established by leveraging \cite[Theorem 5]{ShoT1994}.
\begin{remark}
\label{rem:Rem_5_1}
The condition $\mu_j + \eta_j \le c < 1$ should be interpreted as an abstract assumption on the quality of the constructed subspaces rather than a property that can be verified \emph{a priori} for specific algorithms. In this sense, the result provides a connection to inexact Newton theory rather than a solver-specific convergence guarantee.
\end{remark}
\begin{theorem} \label{thm:nonsingular}
    Let $f(x)$ be continuously differentiable in a neighborhood $B_1(x^*)$ of a solution $x^*$ to $f(x)=0$, with $J_f(x^*)$ nonsingular. Suppose that \Cref{ass:upperbound} holds. Then there exists a neighborhood $B_0(x^*)$ of $x^*$ such that, for any $x_0\in B_0(x^*)$, the \textit{nlKrylov} iterates $\{x_j\}_j$ converge to $x^*$.
\end{theorem}

\begin{proof}
Following the proof of \cite[Theorem 5]{ShoT1994}, it is sufficient to show that \textit{nlKrylov} methods are in fact inexact Newton methods, with the ratio $\tfrac{\|t_j\|}{\|f(x_j)\|}$ uniformly bounded by a constant less than $1$, where $t_j$ denotes the deviation from the exact update, see \eqref{eq:inexactNewton}.
Then, local convergence follows from \cite[Theorem 2.3]{DemES1982}.
Within our \textit{nlKrylov} formulation, we have
\[
  \Delta x_j = x_{j+1}-x_{j}  = -P_jV_j^Tf(x_j) \quad \mbox{ and } \quad 
  t_j = J_f(x_j)\Delta x_j+f(x_j) = (I-J_f(x_j)P_jV_j^T)f(x_j).
\]
Since \eqref{eq:errormatrix2} yields $J_f(x_j)P_j=W_j+V_j$, we obtain
\begin{equation}
\label{eq:tj}
t_j=J_f(x_j)\Delta x_j+f(x_j)=(I-V_jV_j^T-W_jV_j^T)f(x_j).
\end{equation}
Then, by \eqref{eq:gcrnormal} and \Cref{ass:upperbound}, we get
$$\|W_jV_j^Tf(x_j)\| = \|W_jy_j\|=\mu_j\|f(x_j)\|.$$
Thus, using \eqref{eq:tj} and denoting by $\theta_j\in[0,\pi/2]$ the angle between $f(x_j)$ and $\ran(V_j)$, gives
\begin{align*}
\label{eq:inexact}
    \|t_j\|^2=& \;f(x_j)^T(I-V_jV_j^T)^T(I-V_jV_j^T)f(x_j)-f(x_j)^T(I-V_jV_j^T)^TW_jV_j^Tf(x_j)\\
    &-f(x_j)^TV_jW_j^T(I-V_jV_j^T)f(x_j)+f(x_j)^TV_jW_j^TW_jV_j^Tf(x_j)\\
    \le& \; \|f(x_j)\|^2\sin^2\theta_j+2\|f(x_j)\|^2\mu_j\sin\theta_j +\|f(x_j)\|^2\mu_j^2\\
    =&\;\|f(x_j)\|^2(\sin\theta_j+\mu_j)^2.
\end{align*}
By \Cref{ass:upperbound}, 
$$\eta_j \|f(x_j)\|=\|f(x_j)+V_jy_j\|=\|f(x_j)-V_jV_j^Tf(x_j)\|=\|(I-V_jV_j^T)f(x_j)\|=\|f(x_j)\|\sin(\theta_j),$$ and hence, $\sin(\theta_j)=\eta_j$.
Finally, local convergence is guaranteed by \eqref{eq:upperbound} and \cite[Theorem 2.3]{DemES1982} as 
\begin{align*}
\frac{\|t_j\|}{\|f(x_j)\|}\le\eta_j+\mu_j\le c_j\le c<1.  
\end{align*}
\end{proof}
\begin{remark}\label{rem:slowchange}
Thus far, we have implicitly assumed that the Jacobian varies slowly during iterations. However, this assumption has not been explicitly quantified. In our analysis, it is captured by \eqref{eq:errormatrix} and \Cref{ass:upperbound}, i.e., 
the alignment between $\ran(V_j)$ and $\ran(J_f(x_j)P_j)$, and consequently the small norm of $W_j$ reflects a slowly changing linear model. Imposing stronger conditions, such as Lipschitz continuity of $J_f(x)$, would lead to assumptions similar to
\cite[Assumption 2]{He2024}, which require a uniform bound on the quadratic term of the Taylor expansion of $f(x)$ -- a condition that is too restrictive for our framework. Our result, \Cref{thm:nonsingular}, is formulated in terms of the subspaces $\ran(P_j)$ and $\ran(V_j)$, allowing more flexibility in selecting the subroutine $\mathcal{SR}_j$.     
\end{remark}
\subsection{Problems with singular Jacobian}
\label{convsingular}
Extensive results exist on the convergence of Newton's method for singular problems, where the derivative at the root is singular. \cite{Red1978, Red1979} considered $C^3$ functions with a one-dimensional null space of the Jacobian. These results were extended in \cite{DecK1980, DecK1982, DecKK1983, KelS1983}, which provided detailed convergence rate analyses under various scenarios. \cite{Gri1980,Gri1981} showed a starlike domain for which Newton's method converges linearly to the root. For inexact Newton methods applied to singular problems, \cite{Kelley1993} provided convergence results analogous to \cite[Theorem 2.3]{DemES1982},
while \cite{ShoT1994} showed that subspace projected inexact Newton methods converge whenever the corresponding exact Newton method does. \cite{Arg1999} extended the results of \cite{DecKK1983} to inexact Newton schemes. Convergence of quasi-Newton methods, such as Broyden’s method, for singular problems has also been studied in \cite{DecK1985,Man2022}. For singular systems, including underdetermined and overdetermined cases, several results \cite{Chen2008, FerG2011, Zhou2014An, Zhou2014On, ArgG2015, AnaA2016} establish convergence of inexact Newton methods under weak Lipschitz or majorant conditions.
Although \textit{nlKrylov} methods can be interpreted within the inexact Newton framework, the forcing sequence $c_j = \mu_j + \eta_j$  is not explicitly user-controlled; rather, it is implicitly determined by the Krylov subspace construction and the associated projections at each iteration. Consequently, the aforementioned convergence results do not apply directly. Nevertheless, under additional assumptions on $t_j$, generalized convergence results for singular systems can be derived by leveraging the frameworks of \cite{DecKK1983} and \cite{Kelley1993}.
\begin{assumption}\label{ass:singular}
   Let $J_f(x^*)$ have a one dimensional null space $\cN$ spanned by $\phi\in \mathbb{R}^n$ and closed range $\cX$ such that $\mathbb{R}^n = \cN\oplus \cX$. In addition, let $f$ be twice Lipschitz continuously differentiable in a neighborhood of $x^*$, and for any projection $P_\cN$ onto $\cN$ parallel to $\cX$, we have
    $$P_\cN H_f(x^*)(\phi,\phi)\neq0,$$
    where $H_f(x)$ denotes the Hessian of
   a (nonlinear) function $f:\RR^n\rightarrow\RR^n$.
\end{assumption}
Let $\widetilde x=x-x^*$ for $x\in \RR^n$, and define $P_\cX=I-P_\cN$. For singularities satisfying \Cref{ass:singular}, consider the neighborhood of $x^*$ (in terms of $\widetilde x$), i.e.,
$$W(\rho,\gamma)=\Big\{x\in \mathbb{R}^n\mid 0<\|\widetilde x\|\le\rho, \|P_\cX\widetilde x\|\le \gamma\|P_\cN \widetilde x\|\Big\},$$
where $\rho > 0$ and $\gamma > 0$ are sufficiently small.
Let 
\begin{align}
    A(x)&=P_\cX J_f(x)P_\cX, &A_1(x)(\cdot) &= P_\cX H_f(x^*)(\widetilde x,P_\cX(\cdot)),\label{eq:defA}\\
    B(x)&=P_\cX J_f(x)P_\cN, &B_1(x)(\cdot) &= P_\cX H_f(x^*)(\widetilde x,P_\cN(\cdot)),\label{eq:defB}\\
    C(x)&=P_\cN J_f(x)P_\cX, &C_1(x)(\cdot) &= P_\cN H_f(x^*)(\widetilde x,P_\cX(\cdot)),\label{eq:defC}\\
    D(x)&=P_\cN J_f(x)P_\cN, &D_1(x)(\cdot) &= P_\cN H_f(x^*)(\widetilde x,P_\cN(\cdot)), \label{eq:defD}
\end{align}
and let $\bar D_1(x)(\cdot) := P_\cN H_f(x^*)(P_\cN\widetilde x,P_\cN(\cdot))$ be a linear map from $\cN$ to $\cN$. Note that, by \cref{eq:defA}--\cref{eq:defD}, we can write
\begin{equation}
    H_f(x^*)(\widetilde{x},\widetilde{x})=\Big(A_1(x)+B_1(x)+C_1(x)+D_1(x)\Big)\widetilde{x}. \label{eq:hessiansum}    
\end{equation}

\begin{figure}[H]
        \centering
        \includegraphics[scale=0.53]{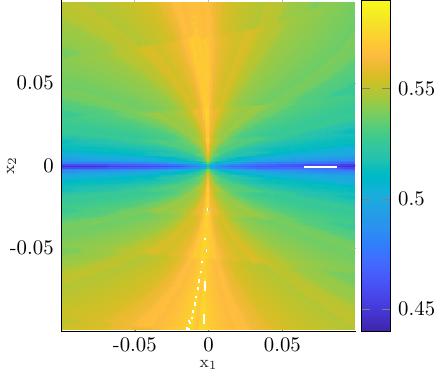}
        \includegraphics[scale=0.53]{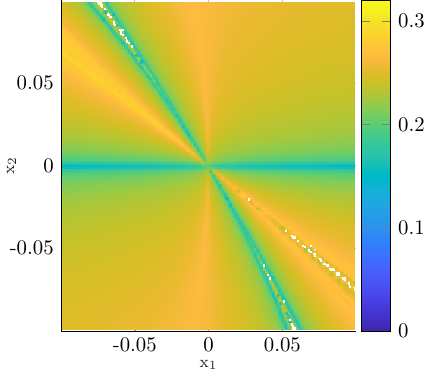}
        \includegraphics[scale=0.53]{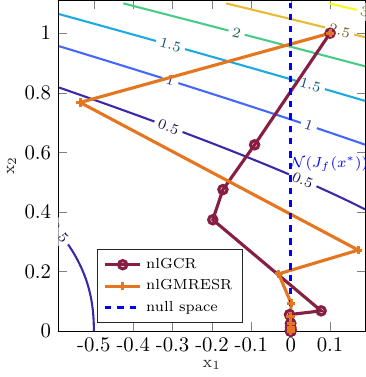}
        \includegraphics[scale=0.53]{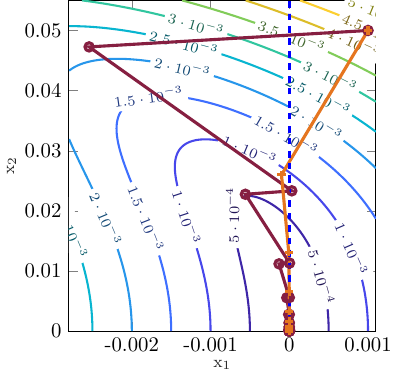}
    \caption{Convergence maps of nlGCR(2) and nlGMRESR(2,2) for \eqref{eq:smallsingular} on $\Omega=[-0.1,0.1]^2$ as well as trajectories for $x^{(1)}_0=[0.1,1]^T$  and $x^{(2)}_0=[0.001,0.05]^T$ along the null space (left to right)}
    \label{fig:visualize}
\end{figure}

\begin{lemma}\label{lem:Jinv}
If $\bar D_1(x)(\cdot)$ is invertible whenever $P_\cN\widetilde x\neq0$, there exist constants $\rho>0$ and $\gamma>0$ such that $J_f(x)$ is invertible in the region $W(\rho,\gamma)$, and
$$J_f(x)^{-1}=P_\cN D_1(x)^{-1}P_\cN+\mathcal{O}(1).$$
\end{lemma}
\begin{remark}
The invertibility of $J_f(x)$ is discussed in \cite{DecKK1983,DecK1980,KelS1983,Red1978}.
\end{remark}
Based on \cite[Theorem 5.1]{DecKK1983}, and assumption $\|t_j\|\le c\|f(x_j)\|^2$, we have the following convergence result for the cases covered under \Cref{ass:singular}.
\begin{theorem}\label{thm:singular}
Let \Cref{ass:singular} hold and define $\widetilde{x} = x - x^*$. Assume for all $P_\cN\widetilde x\neq0$, that $\bar D_1$ is nonsingular as a map on $\cN$. 
If $\|t_j\|\le c\|f(x_j)\|^2$ for some $c$,
then for $\rho >0$ and $\gamma>0$ sufficiently small, $J_f(x_0)^{-1}$ exists for all $x_0\in W(\rho,\gamma)$, and the sequence $\{x_j\}_j$ of iterates generated by a \textit{nlKrylov} method converges to $x^*$ with 
\begin{equation}
\label{eq:ConvQuadTerm}
\|P_\cX(x_{j+1}-x^*)\|\le K_1\|x_{j}-x^*\|^2,
\end{equation}
for some $K_1>0$, and
\begin{equation}
 \lim_{j\to\infty}\frac{\|P_\cN(x_{j+1}-x^*)\|}{\|P_\cN(x_{j}-x^*)\|}=\frac{1}{2}.  \label{eq:limitnull} 
\end{equation}
\end{theorem}
\begin{remark}
The condition $\|t_j\| \le c \|f(x_j)\|^2$ should be interpreted in the same way as the condition in \Cref{rem:Rem_5_1}, i.e., it serves as an abstract analytical assumption rather than a practically verifiable criterion for a specific algorithm.
%
\end{remark}
\begin{proof}
Let us denote $\widetilde{x}_j=x_j-x^*$. For $x_0\in W(\rho,\gamma)$, using \eqref{eq:hessiansum} in the expansion of $f$ and $J_f$, respectively, around $x^*$ gives
\begin{align}
 f(x_0)&=J_f(x^*)\widetilde x_0+\frac12\Big(A_1(x_0)+B_1(x_0)+C_1(x_0)+D_1(x_0)\Big)\widetilde x_0+\mathcal{O}(\|\widetilde x_0\|^3), \label{eq:expandf}\\   
 J_f(x_0)&=J_f(x^*)+\Big(A_1(x_0)+B_1(x_0)+C_1(x_0)+D_1(x_0)\Big)+\mathcal{O}(\|\widetilde x_0\|^2), \label{eq:expandJf}  
\end{align}
and \eqref{eq:expandJf} yields
$$J_f(x_0)\widetilde x_0=J_f(x^*)\widetilde x_0+\Big(A_1(x_0)+B_1(x_0)+C_1(x_0)+D_1(x_0)\Big)\widetilde x_0+\mathcal{O}(\|\widetilde x_0\|^3).$$
Therefore,
\begin{equation}
f(x_0)=J_f(x_0)\widetilde x_0-\frac12\Big(A_1(x_0)+B_1(x_0)+C_1(x_0)+D_1(x_0)\Big)\widetilde x_0+\mathcal{O}(\|\widetilde x_0\|^3). \label{eq:fx0}
\end{equation}
Also, from \eqref{eq:expandf} and $J_f(x^*)P_\cN x\equiv0$, we have $$f(x_0)=J_f(x^*)\widetilde x_0+\mathcal{O}(\|\widetilde x_0\|^2)=J_f(x^*)(P_\cX + P_\cN)\widetilde x_0+\mathcal{O}(\|\widetilde x_0\|^2)=J_f(x^*)P_\cX \widetilde{x_0}+\mathcal{O}(\|\widetilde x_0\|^2),$$ which shows that there is a constant $K > 0$ such that 
$$\|f(x_0)\|\;\le \;K\|P_\cX\widetilde x_0\|\; \le \;\gamma K\|P_\cN\widetilde x_0\|\; \le \; \gamma K\|\widetilde x_0\|.$$ 
Thus, $\|t_0\|\;\le \; c\|f(x_0)\|^2\; \le \; c  \gamma^2K^2\|\widetilde x_0\|^2$, i.e., $\|t_0\|= \gamma^2\mathcal{O}(\|\widetilde{x}_0\|^2)$. It follows from \eqref{eq:fx0} and \Cref{lem:Jinv} that
\begin{align*}
\widetilde x_1&=\widetilde x_0 - J_f(x_0)^{-1}(f(x_0)-t_0)\\
&=J_f(x_0)^{-1}\bigg(\frac12\Big(A_1(x_0)+B_1(x_0)+C_1(x_0)+D_1(x_0)\Big)\widetilde x_0+t_0 +\mathcal{O}(\|\widetilde x_0\|^3)\bigg)\\
&=\frac12\big(P_\cN D_1(x_0)^{-1}P_\cN+\mathcal{O}(1)\big)\bigg(\Big(A_1(x_0)+B_1(x_0)+C_1(x_0)+D_1(x_0)\Big)\widetilde x_0+t_0 +\mathcal{O}(\|\widetilde x_0\|^3)\bigg).
\end{align*}
Moreover, equations \eqref{eq:defA}--\eqref{eq:defD} yield $P_\cN A_1(x_0)\!=\!P_\cN B_1(x_0)\!=\!0$,~$C_1(x_0)\widetilde x_0\!=\!\gamma\mathcal{O}(\|\widetilde x_0\|^2)$ and $P_\cN D_1(x_0)\!=\!D_1(x_0)$.
Thus,
\begin{equation}\label{eq:singularupdate}
\widetilde x_1=\frac{1}{2}P_\cN\widetilde x_0+\gamma P_\cN\mathcal{O}(\|\widetilde x_0\|)+ \gamma^2P_\cN\mathcal{O}(\|\widetilde x_0\|)+\mathcal{O}(\|\widetilde x_0\|^2)+\mathcal{O}(\|\widetilde x_0\|^3).
\end{equation}
Applying $P_\cX$ to \eqref{eq:singularupdate} allows us to show that there exists a constant $K_1 > 0$ such that
$$\|P_\cX\widetilde x_1\|\le K_1\|\widetilde x_0\|^2,$$
proving \eqref{eq:ConvQuadTerm}. Similarly, applying $P_\cN$ to \eqref{eq:singularupdate} shows that for small enough $\gamma > 0$ there exists a constant $K_0 > 0$ with $K_0\gamma<\frac12$ such that $$\big(\frac12-K_0\gamma\big)\|P_\cN\widetilde x_0\|\le \|P_\cN\widetilde x_1\| \le\big(\frac12+K_0\gamma\big)\|P_\cN\widetilde x_0\|.$$ 
Now, we prove the convergence by induction starting with $j=0$. To show that $x_1\in W(\rho,\gamma)$, we set $\rho_0=\|x_0-x^*\|$, $\gamma_0=\gamma$, and note that $\rho_0\le\|P_\cN\widetilde x_0\|+\|P_\cX\widetilde x_0\|\le(1+\gamma_0)\|P_\cN\widetilde x_0\|$. Then
$$\|P_\cX\widetilde x_1\|\le K_1\rho_0^2\le K_1\rho_0(1+\gamma_0)\|P_\cN\widetilde x_0\|\le K_1\rho_0(1+\gamma_0)(\frac12-K_0\gamma_0)^{-1}\|P_\cN\widetilde x_1\|.$$
By letting $\rho_j=\|x_j-x^*\|$, $\gamma_{j}=K_1\rho_{j-1}(1+\gamma_{j-1})(\frac12-K_0\gamma_{j-1})^{-1}$ for $j>0$, we have $\|P_\cX\widetilde x_1\|\le\gamma_1\|P_\cN\widetilde x_1\|$ and $x_1\in W(\rho_1,\gamma_1)$. 
Also,
$$\rho_1\le \|P_\cN\widetilde x_1\|+\|P_\cX\widetilde x_1\|\le (\frac12+K_0\gamma)\|P_\cN\widetilde x_0\|+K_1\rho_0^2\le \left((\frac12+K_0\gamma)(1+\gamma_0)^{-1}+K_1\rho_0\right)\rho_0.$$
For $\rho_0>0$, $\gamma_0>0$ small enough, there exists $\tau\in(\frac12,1)$ such that $\rho_1<\tau\rho_0$, $\gamma_1<\tau\gamma_0$, and therefore $x_1\in W(\rho_1,\gamma_1)\subset W(\rho,\gamma)$.
By induction, repeating the same argument at each step, we obtain
$\rho_{j+1}<\tau\rho_j$, $\gamma_{j+1}<\tau\gamma_j$, and  $x_j\in W(\rho_j,\gamma_j)$ for all $j$.  Thus $x_j\to x^*$ with a convergence rate less than $\tau$, and 
\begin{align*}
    \lim_{j\to\infty}\frac{\|P_\cN(x_{j+1}-x^*)\|}{\|P_\cN(x_{j}-x^*)\|}=\frac{1}{2}.
\end{align*}
\end{proof}
\noindent
Since \Cref{thm:singular} and its proof are rather technical, we illustrate the underlying result using a simple example. Consider the problem
\begin{equation}
f(x)=\begin{bmatrix}x_1+x_2^2\\\frac32x_1x_2+x_2^2+x_2^3\end{bmatrix}, \label{eq:smallsingular}
\end{equation}
discussed in \cite{DecKK1983}. It is straightforward to verify that $x^*=[0,0]^T$ is a root of $f$, with the null space of $J_f(x^*)$ given by $\cN=\spann\{\phi\}$, with $\phi=[0,1]^T$, so that $\dim(\cN) = 1$. Moreover, $P_\cN H_f(x^*)(\phi, \phi) \neq 0$, showing that \eqref{eq:smallsingular} satisfies \Cref{ass:singular}.
We solve \eqref{eq:smallsingular} using nlGCR(2) and nlGMRESR(2,2) under two different configurations. First, we discretize the domain $\Omega=[-0.1,0.1]^2$ using $N=200$ grid points per dimension, and the convergence rates of both methods are evaluated for various initial values. Second, we select two initial points $x^{(1)}_0=[0.1,1]^T$ and $x^{(2)}_0=[0.001,0.05]^T$, and show the trajectories of the iterates generated by the two methods, overlaid with contour lines of $\|f(x)\|$. The results of both experiments are displayed in Figure \ref{fig:visualize}.
\noindent In the first experiment, we observe that for both methods, convergence along the null space (the $x_2$-axis) is considerably slower than in most other regions of the plane. In contrast, convergence orthogonal to the null space, along the $x_1$-axis, is the fastest. For nlGMRESR, an additional diagonal region connecting $p_1 = [-0.06, 0.1]^T$ and $p_2 = [0.06, -0.1]^T$ also exhibits rapid convergence. Overall, nlGMRESR demonstrates superior convergence behavior compared to nlGCR. In both plots, white points refer to initial values that did not lead to convergence within 100 steps, either due to very slow convergence rates or stagnation.
To study convergence along the null space $\mathcal{N}(J_f(x^*))$ (displayed by a dashed blue line), we use the two initial guesses introduced above and visualize their corresponding trajectories. We can see that both methods follow the direction of the null space as they approach the root, with nlGMRESR reaching the null space in fewer steps and remaining closer to it over successive iterations. Upon convergence, the algorithm achieved an approximate twofold reduction in $x_2$ per step, indicating that \eqref{eq:limitnull} holds. 
%

%
%
\section{Practical implementation} \label{sec:implementation}
\subsection{Linear and nonlinear updates} \label{sec:adaptive}
In \cite{He2024}, an adaptive version of nlTGCR($k$) was developed to reduce the number of function evaluations by adapting the residual update from the original (linear) GCR algorithm. 
Let us denote by 
\begin{equation}
    r^{nl}_{j+1}=-f(x_{j+1}),\qquad r^{lin}_{j+1}=r_j-V_jy_j, \label{eq:res}
\end{equation}
the nonlinear and linearized residual, respectively. Assume that in step $j$ of an \textit{nlKrylov} method, we have used the nonlinear residual so far, i.e., $r_j\equiv r_j^{nl}$ in \eqref{eq:res}. Then, if the cosine angle condition holds, i.e.,
\begin{equation}
    \theta_{j+1}:=1-\frac{(r_{j+1}^{nl})^Tr_{j+1}^{lin}}{\|r_{j+1}^{nl}\|\|r_{j+1}^{lin}\|}<\theta \label{eq:adaptiveswitch}   
\end{equation}
for some predefined threshold $\theta$ (see \eqref{eq:paa}), we adaptively switch to the linear version. Once switched, the nonlinear residual is periodically recomputed to ensure the local linear model remains valid, i.e.,
checking \eqref{eq:adaptiveswitch}. 
If the angle becomes too large, we revert to the nonlinear update and discard $P_j$ and $V_j$, effectively restarting the algorithm using $x_0 = x_j$ to avoid inaccurate directions that could hinder convergence. 
\subsection{Automatic restarts} \label{sec:restart}
In the standard GCR-algorithm, we have $V_j^TV_j=I_{n_j}$. Moreover, since $V_j=A P_j$,\; $P_j$ is orthogonal with respect to the inner product induced by $A^TA$, i.e.,
$$\delta_{i,j}=\langle v_i,v_j\rangle = \langle Ap_i,Ap_j\rangle = p_i^T(A^TA)p_j:=\langle p_i,p_j\rangle_{A^T\!A}.$$
If $A$ is \emph{well-conditioned}, orthogonality of $V_j$ induces $A^T\!A$-orthogonality of $P_j$ and allows for a numerically stable procedure. However, in the nonlinear case, we generally do not have $V_j=AP_j$ since $A_j=J_f(x_j)$ changes in every step. To avoid ill-conditioning of $P_j$ caused by the loss of $A_j^T\!A_j$-orthogonality, we employ a slight modification of the automatic restart strategy presented in \cite{TangXHSX2024}.
Consider the Gram-Schmidt process in the $j$-th step of an \textit{nlKrylov} method, i.e., line \ref{line:nlKrylovortho} of \Cref{alg:nlKrylovSR}. With $\beta_i=\widehat{v}^Tv_i$, we expand $P_j$ by 
$$p_{j+1}=\frac{1}{\lVert \widehat{v}\rVert}\left(\widehat{p}-\sum_{i=j_k}^j\beta_i p_i\right).$$
If we assume rounding errors on $p_i$, i.e., the numerically stored basis matrices are $\widetilde{p}_i=p_i+\varepsilon_i$, $i=j_k,\dots,j$, and an error originating from the computation of $p_{j+1}$ denoted by $\delta_{j+1}$, which is assumed to be bounded by $\lVert \delta_{j+1}\rVert_{\infty}\leq \frac{C}{\lVert \widehat{v}\rVert}\lVert \widehat{p} \rVert_\infty$, we get an estimate of the error $\varepsilon_{j+1}$ in the stored vector $p_{j+1}$ by \cite[Eq. (3.6)]{TangXHSX2024} 
$$\lvert \varepsilon_{j+1}\rvert \leq \frac{1}{\lVert \widehat{v}\rVert}(C\cdot\lVert \widehat{p}\rVert_\infty + \sum_{i=j_k}^j\lvert\beta_i\rvert \lVert\varepsilon_i\rVert_\infty),$$
meaning we can keep track of the level of ill-conditioning by introducing the scalar sequence
\begin{equation}
    w_{j+1}=\frac{1}{\lVert \widehat{v}\rVert}(C\cdot\lVert \widehat{p}\rVert_\infty + \sum_{i=j_k}^j\lvert\beta_i\rvert w_i),\quad w_0=0, \label{eq:restart}
\end{equation}
and restarting whenever $w_{j+1}>\tau$, where $\tau>0$ is a user-defined tolerance.  If ill-conditioning is detected and the algorithm is restarted, we set $P_j=V_j=\emptyset$ and restart using 
$$P_{j+1}=\frac{\widehat{p}}{\|J_f(x_{j+1})\widehat{p}\|},\quad V_{j+1}=\frac{J_f(x_{j+1})\widehat{p}}{\|J_f(x_{j+1})\widehat{p}\|},\quad w_{j+1}=C\frac{\|\widehat{p}\|_\infty}{\|J_f(x_{j+1})\widehat{p}\|}.$$
In our implementation, we use $C=1$ and $\tau=10^3$ as default parameters and refer the interested reader to \cite{TangXHSX2024} for a discussion of different parameter choices. 
\subsection{Damping and choosing descent directions}
Within the quasi-Newton framework,  $d_j=P_jy_j=P_jV_j^Tr_j$ can be interpreted as an inexact update direction obtained using the approximation $J_f(x_j)^{-1}\approx P_jV_j^T$. Hence, we cannot guarantee that $d_j$ is a descent direction for $f(x)$~\cite{AriG2023}. Therefore, if $d_j$ is not a descent direction, its negation $-d_j$ should be used to ensure a decrease of $f(x_j+\alpha_jd_j)$ along $d_j$. This condition can be verified by computing $\zeta_j:=\langle r_j,J_f(x_j)d_j\rangle.$
If $\zeta_j>0$, then $d_j$ is a descent direction, otherwise, $-\zeta_j=\langle r_j,J_f(x_j)(-d_j)\rangle$ is likely to be positive and we use $-d_j$ instead. Additionally, we want to perform a simple line search to ensure that $\alpha_j\in(0,1]$ satisfies the Armijo-Goldstein condition \cite{Arm1966} 
\begin{equation}
\lVert f(x_j + \alpha_j d_j)\rVert^2\leq \|f(x_j)\|^2 - c_1\alpha_j \langle r_j,J_f(x_j)d_j\rangle=\|r_j\|^2 - c_1\alpha_j \zeta_j.\label{eq:armijogoldstein}
\end{equation}
Note that, within the backtracking procedure, $\zeta_j$ can be inexpensively approximated using a finite difference approximation, i.e., 
$$\zeta_j=\langle r_j,J_f(x_j)d_j\rangle\approx\frac{1}{\varepsilon}\langle r_j,f(x_j+\varepsilon d_j)-f(x_j)\rangle = \frac{1}{\varepsilon}\langle r_j,f(x_j+\varepsilon d_j)+r_j\rangle,$$
where $\varepsilon=\alpha_j^{(0)}$ is used before the backtracking loop to check if $d_j$ is a descent direction and the estimate is refined inside the loop using $\varepsilon=\alpha_j^{(\ell)}=\frac{\alpha_j^{(0)}}{2^\ell}$. In our implementation, we use $c_1=10^{-3}$ and adapt the initial step-size heuristic from \cite{He2024}, where it is suggested to use $\alpha_0^{(0)}=1$ and
$$\alpha_{j+1}^{(0)}=\begin{cases} \min\lbrace 1,2\alpha_j^{(0)}\rbrace, & \text{if $\alpha_j^{(0)}$ was accepted in iteration $j$} \\\frac{\alpha_j^{(0)}}{2},&\text{else}\end{cases}$$
to avoid small steps while simultaneously reducing the number of line-search steps. The interested reader is referred to the supplementary material for a detailed statement of the algorithm.
Note that, if the linear update is performed, the damped linear residual 
$$r_{j+1}^{lin}=r_j-\alpha_j V_jy_j$$
is used on the left-hand side of \eqref{eq:armijogoldstein} instead of $f(x_j+\alpha_jd_j)$.
\subsection{Solving matrix valued equations}\label{sec:changematrix}
Drawing on the fact that global Krylov subspace methods (see e.g. \cite{JbiMS1999,WanG2007,ZadTW2019}) implicitly vectorize the matrix-valued equations, we can easily implement slight modifications of \textit{nlKrylov} methods to solve matrix-valued problems of the form 
$$ F(X)=0,\quad X\in\RR^{n,p},~F:\RR^{n,p}\rightarrow\RR^{n,p}.$$
Although this extension of the algorithmic framework yields a powerful, memory-efficient tool for solving general nonlinear matrix equations, the idea follows naturally and is deferred to keep the discussion focused. Detailed descriptions of the algorithms can be found in the supplementary material; we will just briefly comment on the modifications and their implementation. First, we have to replace the inner products and norms in every method by their Frobenius equivalents. Note that in {\sc MATLAB}, the inner product $\langle X,Y\rangle_F:=\mathrm{tr}(X^TY)$ can be computed efficiently in $\mathcal{O}(np)$ flops using $$\langle X,Y\rangle _F= \sum_{i,j} x_{ij}y_{ij} = \texttt{sum(X.*Y,'all')}.$$ 
Secondly, every matrix-vector multiplication involving the Jacobian $J_f(x)z$ has to be replaced by the Fr\'echet derivative $L_F(X,Z)$, which can be computed efficiently using a finite difference approximation similar to \eqref{eq:matrixfreederivative} if no explicit formula is available. If there is a need for a high-accuracy approximation, the complex step approximation 
\begin{equation}
L_F(X,Z)\approx\frac{\Im(F(X+\imath \varepsilon Z))}{\varepsilon},\quad\imath=\sqrt{-1},~\Im(a+\imath b)=b,~a,b\in\RR, \label{eq:cs}   
\end{equation}
can be used~\cite{AlmH2010}.
Matrix-vector products have to be performed in terms of the $\dprod$-product
$$\mathcal{U}_j\dprod\gamma:=\mathcal{U}_j(\gamma\otimes I_p)=\sum_{i=1}^j \gamma_iU_i,$$
where $\gamma=(\gamma_i)_i\in\RR^j$ is a vector and $\mathcal{U}_j=[U_1,U_2,\dots,U_j]\in\RR^{n,(j\cdot p)}$ is a block matrix of $(n\times p)$-matrices $U_i$. Finally, one may apply the automatic restart strategy from \autoref{sec:restart} when the maximum matrix norm is used in \eqref{eq:restart}.

\subsection{Subroutine configuration} \label{sec:config}
To conclude the discussion of implementation details, we provide brief remarks on the three specific subroutines employed in the numerical experiments: GMRES as used within nlGMRESR, AGMRES as used within nlLGMRES, and the projected GMRES for $\orth{J_f(x_{j+1})}{V_j}$ as used within nlGCRO. Across all solvers and their global matrix case counterparts, a fixed number of $m$ steps (or $m+k$ steps in the case of AGMRES) is performed without residual monitoring or early breakdown detection, starting from a zero initial guess. While this straightforward approach proves adequate for the purposes of this work, performance could be further improved using standard techniques from inexact Newton methods, such as adaptive tolerance selection. Modified Gram-Schmidt orthogonalization is used throughout, and both the solution $\widehat{p}$ and the quantity $\widehat{v}=J_f(x_{j+1})\widehat{p}$ are recovered via low-rank products involving the Krylov basis matrices, thereby avoiding any additional function evaluations.

When the system matrix is symmetric, the Hessenberg matrix $\underline{H}_m$ arising in GMRES and projected GMRES reduces to a symmetric tridiagonal matrix, while in the case of AGMRES, the augmented matrix $\underline{H}_{m+k}$ takes the form
$$\underline{H}_{m+k}=\begin{bmatrix} \underline{T_m} & \underline{B}_k \\ 0 & R_k\end{bmatrix}\in\RR^{(m+k+1),(m+k)},$$
where $\underline{T}_m\in\RR^{m+1,m}$ is tridiagonal with an extra row, $\underline{B}_k\in\RR^{m+1,k}$ is dense and $R_k\in\RR^{k,k}$ is upper triangular. Symmetry can be detected automatically by checking whether the $(1,3)$-element of $\underline{H}_m$ (respectively $\underline{H}_{m+k}$ in the case of AGMRES) is numerically zero, after which this structure can be exploited in subsequent steps to yield short recurrences within the Gram-Schmidt orthogonalization. As noted in Section \ref{sec:nlgmresr}, the selection of optimal parameters $m$ and $k$ remains an open problem in the context of \textit{nlKrylov} methods, and is discussed extensively for the linear setting in \cite{Van1994}.

In the numerical experiments presented here, parameters are chosen so as to allow as many solvers as possible to converge within competitive runtimes. Broadly, nlGCRO and nlLGMRES perform best when the nonlinearity of the problem is moderate and $k$ is not excessively large. Values of around ten strike a good balance between leveraging the information contained in the outer basis and keeping the computational overhead of the subroutines manageable. In contrast, nlGCR tends to benefit from smaller values of $k$, typically one or two, particularly for problems with symmetric or nearly symmetric Jacobians, as noted in \cite{He2024}. Among all methods considered, nlGMRESR proved to be the most robust and broadly competitive across a range of parameter choices. The selection of an optimal $m$ similarly involves balancing the quality of the local linear model against the overall cost of the subroutine. Moderate values in the range of four to ten performed well across a variety of experiments, while larger values such as $m=20$ may be necessary for strongly nonlinear problems. Although larger values of $m$ and $k$ can accelerate convergence, memory constraints in certain applications may necessitate the use of smaller parameter values.
%
%
%
%
\section{Numerical experiments} \label{sec:experiments}

In this section, we evaluate the performance of the \textit{nlKrylov} methods developed in this work on a variety of nonlinear problems. We compare them against related approaches discussed in \autoref{sec:connection}, including inexact Newton methods, particularly the Jacobian-free Newton-Krylov approach nlOrthomin($k$), and Anderson Acceleration (AA). Unless otherwise specified, we use GMRES as inner solver for all methods. In nlOrthomin, we use a Gauss-Newton-procedure to compute the optimal $y_j=\argmin_y\|f(x_j+P_jy)\|$ at each iteration, where we set $V_j\equiv J_{s_j}(y)$ and restrict the Gauss-Newton iterations to a maximum of $20$ steps. For Newton-Krylov methods, we follow the Eisenstat-Walker heuristic \cite{EisW1996} to adaptively control the linear solver tolerance at each nonlinear iteration, i.e.,
$$\eta_j=\left(\frac{\|r_j\|}{\|r_{j-1}\|}\right)^\alpha,\quad \eta_0=\frac{1}{3},\quad\alpha=\frac{1+\sqrt{5}}{2}.$$ In AA, updates are damped using a fixed parameter $\beta\in\RR$, such that  $x_{j+1}=x_j+\beta f_j$, which corresponds to a root finding problem for the scaled function $f_\beta(x):=\beta f(x)$. This straightforward acceleration scheme was similarly employed by \cite{He2024} as a baseline algorithm. While the authors acknowledge that more powerful fixed-point iterations may exist for specific problem classes, the present approach is deliberately selected for its simplicity and generality, with only the parameters  $\beta$ and $k_{AA}$ requiring careful tuning. Solver results are color-coded in all experiments as follows: nlGCR results are displayed with dark red circles (\textcolor{vtred}{$\bm{\circ}$}), nlGMRESR with orange pluses (\textcolor{vtorange}{$\bm{+}$}), nlGCRO with grey asterisks (\textcolor{vtgrey}{$\bm{*}$}) and nlLGMRES with blue squares (\textcolor{tublue}{$\bm{\square}$}). Solid lines represent the full nonlinear methods, while dashed lines indicate their 
linear-nonlinear variants (see Section \ref{sec:adaptive}), e.g., nlGCR-A. For comparison, nlOrthomin is shown with turquoise crosses (\textcolor{vtturquoise}{$\bm{\times}$}), Newton-Krylov with red diamonds (\textcolor{tured}{$\bm{\Diamond}$}) and AA with green triangles (\textcolor{tulgreen}{$\bm{\triangle}$}). 
We compare solvers in terms of (outer) iterations as well as total function evaluations, where each evaluation of $f(x)$ or matrix-free evaluation of $J_f(x)y$ is counted as a single function evaluation. This measure reflects the extra evaluations needed in inner iterations of nested methods to enhance the local linear model, illustrating the potential gains of linear-update approaches. 
\Cref{tab:nlkrylov_update} presents estimates of the expected number of function evaluations per iteration for the \textit{nlKrylov} methods, excluding those required for damping or adaptive switching as described in \Cref{sec:implementation}. The present AA implementation requires one evaluation of $f(x)$ per iteration. The number of function evaluations required by nlOrthomin($k$) and Newton-Krylov is difficult to estimate a priori, owing to the exact line search employed in nlOrthomin and the adaptive forcing strategy in Newton-Krylov, respectively. Note that the total computational cost also includes operations like inner products, which are not counted here and may significantly affect runtime depending on the cost of function evaluations. The total runtime of all algorithms, under the experimental setup described in the subsequent sections, is reported in \Cref{tab:runtime}. It should be noted that the present implementations have not been fully optimized with respect to numerical performance. All experiments are conducted using {\sc MATLAB} 2024b on a notebook with AMD Ryzen 7 PRO 4750U CPU and integrated Radeon RX Vega 7 GPU.
\begin{table}[H]
    \centering
    \resizebox{\linewidth}{!}{
    \begin{tabular}{l|c|c|c|c|c|c}
         Solver & LJ & \makecell[c]{H-Equation\\$\omega=(0.99\,\vert\,1)$} & \makecell[c]{Bratu\\$\lambda=(0.5\,\vert\,6)$} & NARE & NEP & NEPv ($F\,\vert\,\widetilde{F}$)\\
         \hline 
         nlGCR & $0.8893$& ($1.7706$\sep\dnc{37.429}) & ($0.7181$\sep\dnc{12.027}) & $64.966$ & \dnc{0.3021} & (\dnc{2.5899}\sep\dnc{})  \\
         \hline
         nlGCR-A & \textcolor{tured}{$0.5653$}& ($1.7416$\sep\dnc{19.995}) & ($5.2733$\sep\dnc{11.753}) & $58.699$ & \dnc{0.1963} & (\dnc{0.6822}\sep\dnc{})\\
         \hline 
         nlGMRESR & $0.6426$& ($2.5119$\sep$8.7678$) & ($0.1351$\sep\textcolor{tured}{$0.5296$}) & $32.975$& \textcolor{tured}{$0.1958$} & (\dnc{11.47}\sep$6.0326$)\\
         \hline
         nlGMRESR-A & $0.6255$& ($2.5453$\sep$8.5156$) &($0.1381$\sep$0.555$) & \textcolor{tured}{$32.89$}& $0.3827$ & ($6.9538$\sep$5.9278$)\\
         \hline 
         nlGCRO & $1.055$ &(\dnc{83.435}\sep\dnc{85.967}) & (\textcolor{tured}{$0.1299$}\sep$3.9505$) & $138.94$& $0.266$ & (\dnc{11.89}\sep\dnc{})\\
         \hline
         nlGCRO-A & $1.2576$ &($10.825$\sep$38.906$) & ($0.2495$\sep$0.8819$) & $139.43$& $0.4195$ & (\dnc{10.004}\sep$7.4485$)\\
         \hline 
         nlLGMRES & $0.8581$&($8.2075$\sep$27.714$) & ($0.2312$\sep$0.7934$)& $71.327$ & $0.3326$ & (\dnc{15.266}\sep$5.316$)\\ 
        \hline
        nlLGMRES-A & $0.75941$ & ($8.2468$\sep$27.295$) & ($0.2561$\sep$2.0077$) & $71.636$ & $0.5053$ & ($10.184$\sep$5.2347$)\\
        \hline
        AA & $1.1186$& (\textcolor{tured}{$1.069$}\sep\textcolor{tured}{$6.3679$})& ($4.7845$\sep\dnc{196.89})& $323.39$& $2.5412$ & (\textcolor{tured}{$6.3783$}\sep$6.7799$)\\
        \hline 
        JFNK & $0.6821$& ($2.848$\sep$10.633$) & ($0.1548$\sep$0.6237$)& $ 39.199$ & $0.9649$ & ($8.3591$\sep\textcolor{tured}{$3.3295$})\\
        \hline 
        nlOrthomin & $2.4641$ &($6.166$\sep\dnc{$12.126$}) &($1.6299$\sep\dnc{564.78}) & $154.22$ & \dnc{1.0804} & (\dnc{7.4621}\sep\dnc{})
    \end{tabular}
    }
    \caption{Runtime comparison of nonlinear solvers for all problems. Timings indicated by \dnc{} either did not converge or did not find the desired solution, e.g., some other than the smallest magnitude eigenvalue. The fastest converged method for every problem is marked in \textcolor{tured}{red}.}
    \label{tab:runtime}
\end{table}

\subsection{The Lennard-Jones problem}
\begin{figure}[H]
    \myplotsvert{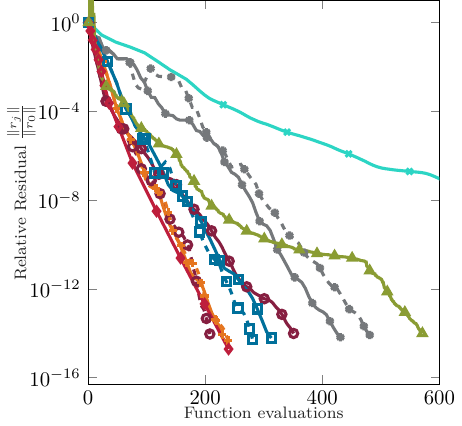}
    \caption{Convergence results for Lennard-Jones problem.}
    \label{fig:lj}
\end{figure}

\begin{figure}[H]
\centering
\includegraphics[width=.8\textwidth]{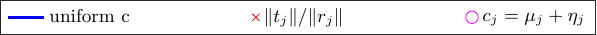}
\begin{minipage}{0.375\textwidth}
\centering
\includegraphics[scale=.6]{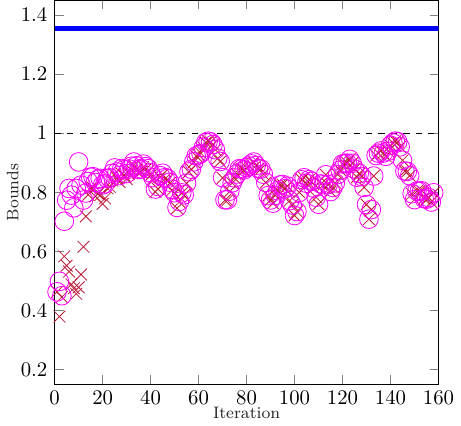}
\end{minipage}
\begin{minipage}{0.375\textwidth}
\centering
\includegraphics[scale=.6]{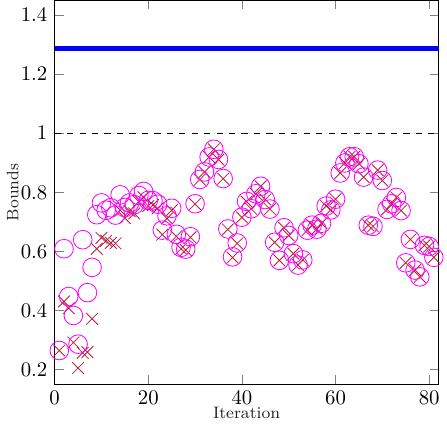}
\end{minipage}
\caption{Theoretical and observed convergence bounds for nlGCR (left) and nlGMRESR (right) applied to the Lennard-Jones problem.}
\label{fig:ljbounds}
\end{figure}

We first consider the Lennard-Jones problem, also covered in \cite{He2024}, a molecular optimization problem that seeks to determine atomic positions that minimize the total potential energy of a molecule described by the Lennard-Jones-Potential
\begin{equation}
    E(Y) = 4\sum_{i=1}^N\sum_{j=1}^{i-1}\left(\frac{1}{\| y_i-y_j\|^{12}}-\frac{1}{\|y_i-y_j\|^6}\right), \label{eq:LJenergy}
\end{equation}
where the matrix $Y=[y_1,\dots,y_N]\in\RR^{3,N}$ contains the 3-dimensional coordinate vectors $y_i\in\RR^3$ for atoms $i=1,\dots,N$. A local minimizer of \eqref{eq:LJenergy} can be computed by finding roots of 
$$f(x)=\nabla E(\mathrm{vec}(Y)),$$
where $x\in\RR^{3\cdot N}$ is the vectorization of $Y$. The derivative $J_f(x)\Delta{x}$ is applied to high accuracy by a complex-step approximation similar to \eqref{eq:cs} using $\varepsilon=10^{-10}$. Note that $J_f(x)$ is symmetric as the Hessian of $E$, so we expect the CS approximation to yield a numerically self-adjoint linear operator. Following \cite{He2024}, an Argon cluster is initialized from a perturbed Face-Centered-Cubic (FCC) structure. Using $3$ cells per direction gives $3^3=27$ unit cells, each with 4 atoms, resulting in $N=4 \cdot 27 = 108$ and a root finding problem of size $n=3\cdot N=324$. Setup files and initial state are freely available at \href{https://github.com/Data-driven-numerical-methods/Nonlinear-Truncated-Conjugate-Residual}{GitHub} \cite{Git2024}.\\[0.01in]
\emph{Experiment 1:} In our first experiment, we compare the convergence properties of various methods using the Lennard-Jones problem as a test case. We use a truncation window of $k=2$ for \textit{nlKrylov} methods and $m=5$ for all nested methods. Adaptive methods use $\theta=10^{-3}$ as the switching tolerance. For AA, we use a truncation window of $k_{AA}=10$ and damping parameter $\beta=-3\times10^{-4}$. Newton-Krylov uses a maximum of $40$ steps of MINRES for the inner solve. The iteration is stopped after $250$ steps or when a relative tolerance of $10^{-14}$ is reached. The convergence results are displayed in \Cref{fig:lj}, where we restrict the view to at most $200$ iterations and $600$ function evaluations.
The left panel of \Cref{fig:lj} illustrates that the nested \textit{nlKrylov} methods (nlGMRES, nlGCRO, and nlLGMRES) converge more rapidly than nlGCR, owing to the improved local linear model provided by the inner solve. Among these, nlGMRESR and nlLGMRES exhibit comparable iteration counts, while nlGCRO requires few more iterations. The nlOrthomin method behaves similarly to nlGCR in terms of iteration count. However, when considering the number of function evaluations, the exact line search and full construction of $V_j=J_f(x_j)P_j$ in nlOrthomin result in more than three times as many evaluations as nlGCR. Comparing the full nonlinear methods, both nlGMRESR and nlLGMRES require fewer function evaluations than nlGCR. Among the adaptive variants, all \textit{nlKrylov} methods except nlGCRO benefit from switching to the linear-update version, with the effect being most pronounced for nlGCR-A, which surpasses nlGMRESR in efficiency upon convergence. Because the inner solve incurs additional function evaluations, the relative gain is smaller for the nested variants. In this example, nlGCR-A, nlGMRESR, and nlGMRESR-A perform comparably to Newton-MINRES in terms of total function evaluations. This is further reflected in the computational timings reported in \Cref{tab:runtime}, where nlGMRESR and nlGMRESR-A demonstrate improved performance over Newton-Krylov, trailing only nlGCR--A, which achieves the best overall runtime. The relatively strong nonlinearity of the problem limits the benefit of reusing the outer basis in the inner iteration, as done in nlGCRO and nlLGMRES.
\noindent\emph{Experiment 2:} In the second experiment, we use the Lennard-Jones problem to compare the theoretical and observed convergence bounds for nlGCR and nlGMRESR (\Cref{thm:nonsingular}) with those from~\cite{He2024}. Hence, we use the same setup as before but now we choose $m=2$ and $k=5$. In \Cref{fig:ljbounds}, we display the uniform bound $c=\mu+\eta$ required in \cite[Eq. (4.33) and (4.34)]{He2024} as well as $c_j=\mu_j+\eta_j$ from \eqref{eq:upperbound} and the observed relative inexactness $\vartheta_j:=\frac{\|t_j\|}{\|r_j\|}$. Recall that \Cref{thm:nonsingular} guarantees the convergence of \textit{nlKrylov} methods provided $\vartheta_j\leq c_j<1$. We can see in \Cref{fig:ljbounds}, that the uniform bound $c>1$ does not satisfy the assumption and, as such, cannot provide a theoretical foundation for convergence. However, the same figure shows that we always have $c_j<1$ and $\vartheta_j\leq c_j$ in the experiments and for some $j_0\in\mathbb{N}$, $\vartheta_j=c_j$ for all $j\geq j_0$. For nlGCR, $j_0 \approx 40$ while for nlGMRESR, $j_0 \approx 25$. Thus, we can see that the bound presented in \eqref{eq:upperbound} is sufficient for convergence of the methods as it is always an upper bound for $\vartheta_j$ and, at convergence, this upper bound is sharp. 
\noindent 
\subsection{The Chandrasekhar H-equation}

\begin{figure}[H]
    \myplotsvert{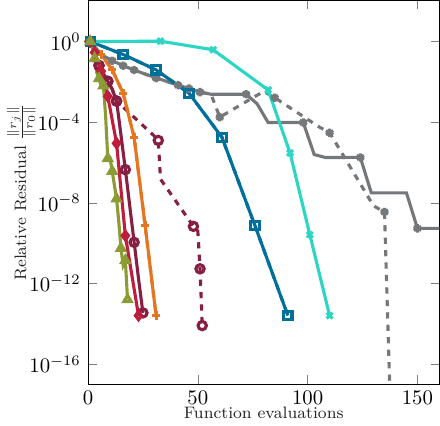}
    \caption{Convergence result for H-equation with nonsingular Jacobian ($\omega=0.99$).}
    \label{fig:ch099}
\end{figure}

For our second example, we consider the integral equation 
\begin{equation}
F(H)(\mu)=H(\mu)-\left(1-\frac{\omega}2\int_0^1\frac{\mu}{\mu+\nu}H(\nu)d\nu\right)^{-1}=0,  \label{eq:hequation}   
\end{equation}
where a continuously differentiable function $H\in\mathcal{C}([0,1])$ satisfying \eqref{eq:hequation} is to be determined. Equation \eqref{eq:hequation}, known as the Chandrasekhar H-equation~\cite{Chandrasekhar1960}, originates from radiative transfer problems. This equation has been examined in several studies focusing on the convergence properties of Newton’s method for singular problems
\cite{DecK1982,Kelley2018}. In particular, \cite{DecK1982} shows that \eqref{eq:hequation} has a singular Jacobian when $\omega=1$, with $\dim(\cN)=1$. As described in \cite{Kelley2018}, using $n$ equally spaced grid points $\mu_i=(i-1/2)/n$, $i=1,\dots,n$, and setting $h_i=H(\mu_i)$, \eqref{eq:hequation} can be discretized via the midpoint rule as 
\begin{equation*}
 \left[f(h)\right]_i=h_i-\left(1-\frac{\omega}{2n}\sum_{j=1}^n\frac{ih_j}{i+j-1}\right)^{-1},\quad i=1,\dots,n,  
\end{equation*}\noindent
where the summation term can be written in terms of a product $K\!\cdot\! h$ involving a Hankel matrix $K\in\RR^{n,n}$ and evaluated efficiently using the Fast Fourier Transform (FFT). The resulting formulation is an $n$-dimensional root finding problem $f(h)=0$. In our experiments, we use $n=10^5$ and consider both the nonsingular case ($\omega=0.99$) and the singular case ($\omega=1$).\\[0.01in]
\emph{Experiment 1:} In the first experiment, we consider the nonsingular case $\omega=0.99$. We use a truncation window of $k=10$ for the \textit{nlKrylov} methods, nlOrthomin and Anderson Acceleration. For nested variants, the inner linear systems are solved using $m=4$ steps of GMRES. The adaptive methods use $\theta=10^{-2}$ and Anderson Acceleration is executed with $\beta=-1$. The Newton-Krylov method is limited to a maximum of $100$ GMRES iterations per Newton step. All iterations are \makebox[\linewidth][s]{terminated after $30$ steps or once the relative residual falls below $10^{-12}$. As shown in \Cref{fig:ch099},} \par 
\noindent the adaptive variants of nlGMRESR and nlLGMRES did not switch to their linear counterparts during the iterations, resulting in overlapping convergence curves. In contrast, both nlGCR-A and nlGCRO-A activated the linear updates. Notably, nlGCRO--A was able to overcome the stagnation observed in its nonlinear counterpart while nlGCR--A did not benefit from the linear update. In terms of iteration counts, nlGMRESR and nlLGMRES converged comparably to Newton-GMRES. When considering the total number of function evaluations, nlGMRESR and nlGCR remained close to Newton-GMRES in efficiency, whereas nlLGMRES required more than twice as many evaluations. Although nlOrthomin reached the desired tolerance in the same number of iterations as Newton-GMRES, it did so at the cost of significantly more function evaluations. With respect to computational timings, AA(10) outperformed all other methods in this test case, with nlGCR and its adaptive variant serving as the closest competing approaches.
\vspace{-.5cm}
\begin{figure}[H]
    \myplotsvert{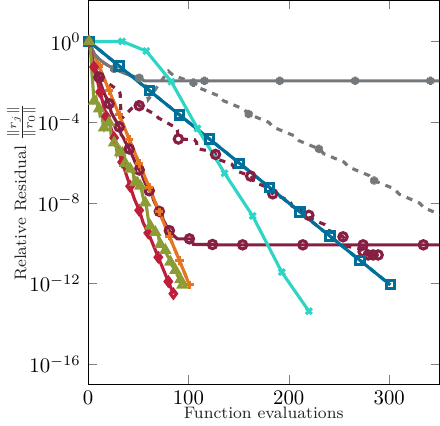}
    \caption{Convergence result for H-equation with singular Jacobian ($\omega=1$).}
    \label{fig:ch100}
\end{figure}
\noindent\emph{Experiment 2:} Our second experiment considers the singular case $\omega=1$. We use the same setup as in Experiment 1, except that the maximum number of (outer) iterations is increased to $100$. As shown in~\Cref{fig:ch100} all methods require substantially more iterations compared to the nonsingular case. Once again, the adaptive variants of nlGMRESR and nlLGMRES did not switch to their linear counterparts. We observe that nlGCR and nlGCRO failed to converge for this problem, whereas their adaptive variants significantly improved accuracy, i.e., nlGCR-A improved from approximately $10^{-10}$ to around $2\times10^{-11}$ before stagnating and nlGCRO-A managed to overcome the stagnation at around $10^{-2}$ and converged after around $70$ iterations. Upon convergence, nlGMRESR, nlLGMRES and Newton-GMRES required the same number of iterations, with nlGMRESR and Newton-Krylov showing comparable numbers of function evaluations. In contrast, nlOrthomin achieved convergence in fewer iterations but incurred more than twice the number of function evaluations. However, the computed solution did not agree with the reference solution. Again, AA(10) was the most efficient among all methods at around 6 seconds, with nlGMRESR and its adaptive version being the most competitive at around 8.5 seconds of total runtime. These results indicate that the adaptive variants can enhance convergence relative to the fully nonlinear methods and, in some cases, overcome the stagnation observed in the original nonlinear iterations. Moreover, most methods successfully solved the singular problem up to the desired accuracy. 
\subsection{The symmetric Bratu problem}\label{sec:bratu}
As a third benchmark problem, we study the well-known symmetric Bratu problem \cite{HajJB2018}, arising from the discretization of the nonlinear PDE 
\begin{eqnarray*}
    \Delta u + \lambda e^u &=&0, \quad(y,z)\in\Omega\\
    u &=&0, \quad(y,z)\in\partial\Omega,
\end{eqnarray*}
with $\Omega=(0,1)\times(0,1)$ and $\lambda=0.5$. Using an equally spaced grid  $y_i=z_i=i\cdot h$, $h=\frac{1}{N+1}$, $i=1,\dots,N$ in both spatial directions, the problem reduces to the following nonlinear problem
\begin{equation}
    f(x)=Lx - h^2\lambda\exp(x)=0 \label{eq:bratu},
\end{equation}
where $L\in\RR^{n,n}$ is a 2D-Laplacian ($L = L_N\otimes I_N + I_N\otimes L_N,~ L_N=\mathrm{tridiag}(-1,2,-1)\in\RR^{N,N}$), $n=N^2$ and $\exp(\cdot)$ is applied element-wise. 
From \eqref{eq:bratu}, it follows immediately that 
$$J_f(x)\Delta{x}=L\Delta{x} - h^2\lambda\exp(x)\odot \Delta{x}, $$
where $x\odot\Delta{x}$ is the Hadamard product. In our experiments, we follow the setup of \cite{He2024} with $N=100$, giving a problem size $n=10,000$, and initialize with the all-ones vector $x_0=\mathbf{1}_n$.

\begin{figure}[H]
    \myplotsvert{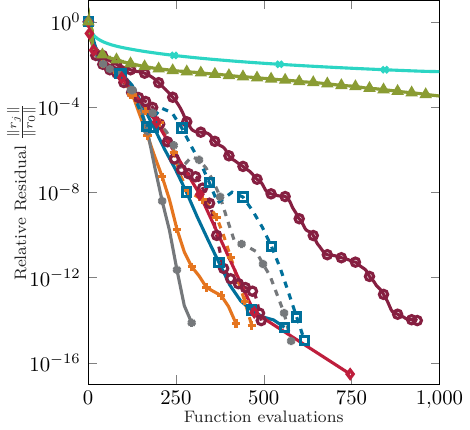}
    \caption{Convergence results for Bratu problem}
    \label{fig:bratu}
\end{figure}

\smallskip
\noindent
\emph{Experiment 1:} 
As before, we analyze the convergence of various methods on the Bratu problem. All truncated methods use $k = 10$, nested \textit{nlKrylov} methods employ $20$ GMRES steps, and adaptive methods switch with a threshold $\theta=10^{-3}$. Anderson Acceleration is damped with $\beta=0.1$, and MINRES is used as the inner solver in Newton's method due to the symmetry of the Jacobian. Iterations are stopped once a relative tolerance of $10^{-14}$ is reached. The results are presented in \Cref{fig:bratu}. Here, we observe the impact of refined local linear models for moderately nonlinear problems. With a relatively large truncation window of $k = 10$ and $m=20$ within $\mathcal{SR}_j$, all nested methods converged within $30$ iterations, whereas nlGCR required approximately $500$ iterations. In particular, nlGCRO benefits from the slowly varying Jacobian, requiring only about one-third of the function evaluations of nlGCR. 
This reduction is likewise reflected in the computational timings, where nlGCRO proves to be the most efficient method among all approaches considered. Both nested methods and their adaptive variants, as well as nlGCR-A, outperform Newton-MINRES in terms of function evaluations. Remarkably, nlOrthomin achieves a similar number of iterations as nlGCR but requires more than ten times as many function evaluations.

\begin{figure}[H]
\centering
\includegraphics[width=.8\textwidth]{newplots/bound_legend.pdf}
\begin{minipage}{0.2445\textwidth}
\centering
\includegraphics[width=\textwidth,height=.2\textheight]{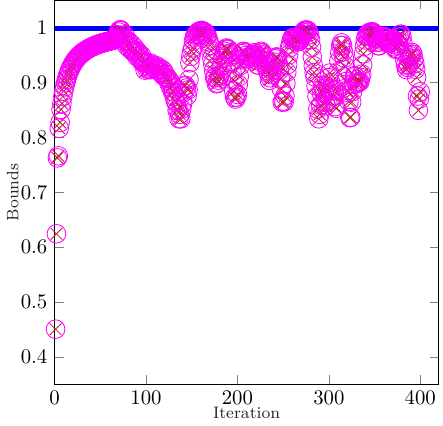}
\end{minipage}
\begin{minipage}{0.2445\textwidth}
\centering
\includegraphics[width=\textwidth,height=.2\textheight]{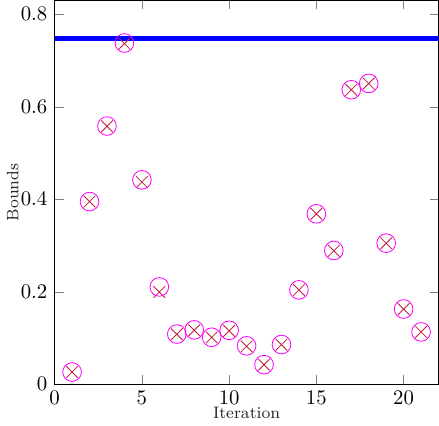}
\end{minipage}
\begin{minipage}{0.2445\textwidth}
\centering
\includegraphics[width=\textwidth,height=.2\textheight]{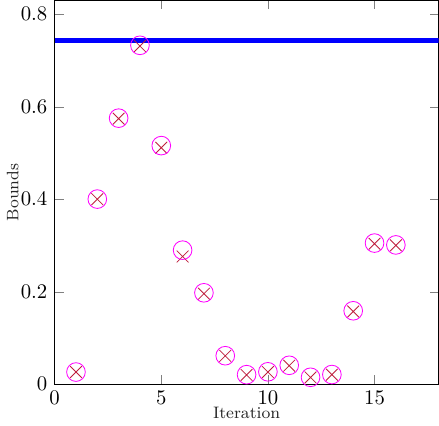}
\end{minipage}
\begin{minipage}{0.2445\textwidth}
\centering
\includegraphics[width=\textwidth,height=.2\textheight]{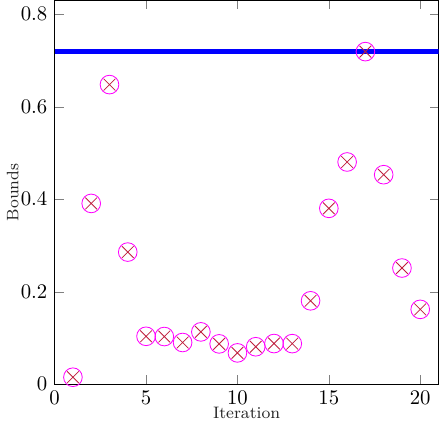}
\end{minipage}
\caption{Theoretical and observed convergence bounds for the Bratu problem, shown for nlGCR, nlGMRESR, nlGCRO, and nlLGMRES (left to right).}
\label{fig:bratubounds}
\end{figure}

\smallskip
\noindent\emph{Experiment 2:} In the second experiment, we analyze the theoretical bounds from \Cref{thm:nonsingular}. Using the same setup as before, we compare the uniform bound $c$ from \cite{He2024}, the sequence $c_j=\mu_j+\eta_j$ from \eqref{eq:upperbound}, and the actual ratio $\vartheta_j=\frac{\|t_j\|}{\|r_j\|}$. As shown in \Cref{fig:bratubounds}, both the uniform bound $c$ as well as $c_j$ are smaller than one, thereby guaranteeing convergence of all algorithms. We also observe that for nlGCR and nlGCRO, $c$ provides a sharp upper bound for $c_j$ and $\vartheta_j$, whereas for nlGMRESR and nlLGMRES, $c_j$ is consistently slightly smaller than $c$.
\noindent
On the other hand, for most iterations $j$, we observe that $c_j\ll c$, indicating that the uniform bound $c$ significantly overestimates the actual error, whereas the $c_j$'s remain highly accurate. Furthermore, the smaller values of $c_j$ and $\vartheta_j$ for the nested methods compared to nlGCR demonstrate that convergence, measured by the number of outer GCR steps, is substantially faster for these methods, consistent with the results shown in \Cref{fig:bratu}. Once again, $c_j$ provides an asymptotically sharp bound for $\vartheta_j$ in all cases.

\subsection{A nonlinear algebraic Riccati equation}
In our final example, we consider a non-symmetric algebraic Riccati equation (NARE) of the form 
\begin{equation}
    \mathcal{R}(X) = FG^T + AX + XB - XPQ^TX =0, \label{eq:nare}
\end{equation}
where $A\in\RR^{n,n}$, $B\in\RR^{p,p}$, $F\in\RR^{n,r}$, $G\in\RR^{p,r}$, $P\in\RR^{p,s}$, $Q\in\RR^{n,s}$ and $X\in\RR^{n,p}$ is the solution of interest. In our application, $p\ll n$ and all matrices are sparse. In particular, we adopt a setup similar to \cite[Example 1]{BenKS2016}, where 
\begin{align*}
    & A= \begin{bmatrix}3 & -1 & & \\ & \ddots & \ddots & \\ & & 3 & -1\\ -1 & & & 1.9
    \end{bmatrix},  &&B=\begin{bmatrix}2 & -1 & & \\ & 3 & \ddots & \\ & & \ddots & -1\\ -1 & & & 3\end{bmatrix}, &&G=I_{p,r},\quad Q=I_{n,s} \\
    &F_1=\begin{bmatrix}-1 & -1 & & \\ & \ddots & \ddots & \\ & & -1 & -1 \\ & & & -0.9
    \end{bmatrix}\in\RR^{r,r}, &&P_1=\begin{bmatrix}1 & & &\\1 & 1 & & \\ & \ddots & \ddots & \\ & & 1 & 1 \end{bmatrix} \in\RR^{s,s}, && F=\begin{bmatrix}F_1\\0\end{bmatrix},\quad P=\begin{bmatrix}P_1\\0\end{bmatrix},
\end{align*}
with $n=30,000$, $p=200$, $r=3$ and $s=5$. Equation \eqref{eq:nare} yields the Fr\'echet derivative
\[
L_\mathcal{R}(X,\Delta{X})=A\Delta{X} + \Delta{X}B - \left(\Delta{X}PQ^TX + XPQ^T\Delta{X}\right).
\]
A highly efficient low-rank Newton-ADI algorithm for NAREs is proposed in \cite{BenKS2016}. For this illustrative example, we instead compare our methods to standard JFNK applied to \eqref{eq:nare}, which requires reducing $p$ by two orders of magnitude compared to \cite{BenKS2016} to satisfy memory limitations.

\begin{figure}[H]
    \myplotsvert{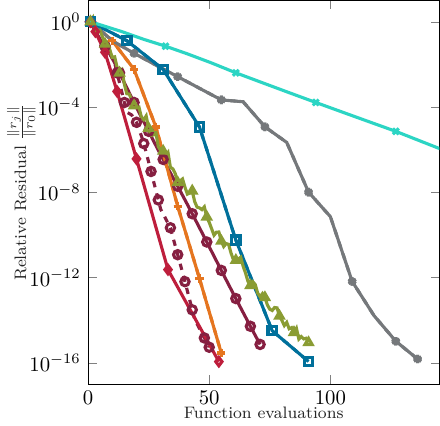}
    \caption{Convergence results for NARE problem}
    \label{fig:nare}
\end{figure}

\smallskip
\noindent\emph{Experiment 1:} In this experiment, we use a truncation window of $k=6$ for the \textit{nlKrylov} methods and nlOrthomin. Nested methods use $m=8$, adaptive methods switch with $\theta=10^{-4}$, and Anderson Acceleration uses a truncation window $k_{AA}=12$ with damping parameter $\beta=0.1$. Newton's method performs up to $50$ steps of global GMRES for the inner solve. Iterations are stopped after $100$ steps or when the relative residual falls below $\tau=10^{-15}$. \Cref{fig:nare} presents the results obtained with the initial guess of all-zeros $X_0=0_{n,p}$, showing that \textit{nlKrylov} methods can successfully solve matrix-valued root finding problems. 

\noindent As in previous experiments, adaptive versions of nlGMRESR, nlGCRO, and nlLGMRES did not switch to the linear update, while nlGCR-A switched after about $10$ iterations. Both nlGMRESR and nlLGMRES converged in the same number of iterations as Newton-GMRES; nlGMRESR required a comparable number of function evaluations, while nlLGMRES needed significantly more. Consistent with other experiments, nlOrthomin shows similar iteration counts to nlGCR, but its function evaluations grow too quickly to remain competitive. The timing results indicate that nlGMRESR and its adaptive variant were the top-performing solvers in this test. Notably, Newton-Krylov methods required $p$ to be reduced by two orders of magnitude to stay within memory constraints, while \textit{nlKrylov} methods handled larger values of $p$ without issue, owing to their comparatively modest storage requirements.
\section{Conclusion and future work}
In this work, we have developed a unified framework for nonlinear Krylov methods (\textit{nlKrylov}) derived from nested GCR-type linear solvers. This framework systematically generalizes classical methods such as GMRESR, GCRO, and LGMRES to the nonlinear setting, yielding algorithms including nlGMRESR, nlGCRO, and nlLGMRES. All methods in this class preserve the essential structural properties of their linear counterparts while admitting an interpretation as inexact Newton schemes, thereby providing a unified viewpoint linking Krylov recycling strategies and Newton-type iterations.

Across the numerical experiments, nlGMRESR emerges as the most consistently robust and competitive variant among the proposed methods. In contrast, nlGCRO and nlLGMRES exhibit more problem-dependent behavior, and their performance is particularly favorable in regimes where the nonlinearities evolve slowly and subspace recycling is most effective.
When compared against Newton--Krylov methods, we emphasize that the proposed approaches do not constitute a uniform replacement but rather a complementary class of solvers. In terms of the total number of function evaluations, Newton--Krylov methods remain a strong and highly competitive baseline across most test problems. The advantage of \textit{nlKrylov} methods is therefore more nuanced: they can achieve comparable performance in several cases (notably Lennard--Jones and Chandrasekhar H-equation) and show clear benefits in problems with slowly varying Jacobians (such as the Bratu problem), where recycling and subspace reuse are most effective.

We have established convergence results under relaxed assumptions and clarified the connections between the proposed methods and existing nonlinear solvers, including nlOrthomin and projection-based Newton frameworks. Moreover, the framework extends naturally to nonlinear matrix equations, which are increasingly relevant in high-dimensional and data-driven applications.

Several directions for future work remain open. These include adaptive strategies for selecting the parameters $m$ and $k$, which strongly influence performance, as well as more systematic sensitivity analyses to better characterize their impact. Further promising avenues include extensions to stochastic or derivative-free settings, adaptive memory truncation for scalability, and deeper integration with Newton--Krylov strategies through dynamic selection of Krylov subspace dimensions.
\section*{Code availability}
The authors are applying this article for the ''SIAM Reproducibility Badge: Code and data available''. Codes to reproduce the numerical experiments are publicly available via the GitHub repository \url{https://github.com/amiedlar/nlKrylov}.
\section*{Acknowledgements}
The first author sincerely acknowledges the hospitality of the Department of Mathematics at Virginia Tech during his visits in Spring 2024 and Fall 2025. Ning Wan and Agnieszka Mi\k{e}dlar gratefully acknowledge support from the National Science Foundation under grants DMS \#2144181 and \#2324958. The authors are obliged to Mark Embree, Eric de Sturler, Heike Fa{\ss}bender and Yousef Saad for comments on earlier versions of this draft. Finally, the authors extend their sincere gratitude to the anonymous referees, whose perceptive suggestions and close reading contributed meaningfully to the final quality of this work.
\section*{Declaration of competing interests}
The authors declare that they have no known competing financial interests or personal
relationships that could have appeared to influence the work reported in this paper.
\section*{CRediT author statement} \noindent\textbf{Tom Werner:} Conceptualization, Investigation, Methodology, Software, Writing - Original Draft, Review and Editing. \textbf{Ning Wan:} Conceptualization, Investigation, Methodology, Software, Writing - Original Draft, Review and Editing. \textbf{Agnieszka Mi\k{e}dlar:} Conceptualization, Methodology, Supervision, Writing - Original Draft, Review and Editing.
\appendix
\vspace{1cm}
\begin{center}
    \MakeUppercase{\textbf{Supplementary materials: \textit{nlKrylov}: A unified framework for nonlinear GCR-type Krylov subspace methods}}\\[.1in]
    \small{\textsc{Tom Werner, Ning Wan, and Agnieszka Mi\k{e}dlar}}
\end{center}

\section{Extended discussion of algorithms}
In this first section, we want to provide an extended discussion of the main algorithms referenced in the paper. Those include the particular algorithms GCR, GMRESR, GCRO and LGMRES that fall into the nested Krylov class, as well as their nonlinear companions fitting into the \textit{nlKrylov} class (Section \ref{sec:basenlkrylov}). We also present global versions of these methods to solve linear and nonlinear operator equations in Section \ref{sec:globalnlkrylov}. To cap off the section, we state additional algorithms used in the paper that are related to the \textit{nlKrylov} family in Section \ref{sec:addalg}.  
\subsection{Base versions of nlKrylov methods}\label{sec:basenlkrylov}
In this Section, we provide a one-to-one comparison and detailed statement of the linear and nonlinear Krylov subspace methods presented in this paper. Those include (nl)GCR($k$), (nl)GMRESR($m,k$), (nl)GCRO($m,k$) and (nl)LGMRES($m,k$). The linear version, designed to solve the linear system $Ax=b$, is displayed on the left, while its nonlinear extension for the root finding problem $f(x)=0$ can be found on the right. nlGCR($k$) (\Cref{alg:nlGCR}) has been already proposed in \cite{He2024} as "nlTGCR($k$)" while the other nonlinear methods (\Cref{alg:nlGMRESR}, \Cref{alg:nlGCRO} and \Cref{alg:nlLGMRESO}) have been proposed in this paper as notable members of the \textit{nlKrylov} family.\\
\begin{minipage}{.495\textwidth}
\begin{algorithm}[H]
    \caption{GCR($k$) for $Ax=b$  \ \cite{DeS1996,Van1994}}    
    \begin{algorithmic}[1]
    \setstretch{1.05}
    \REQUIRE $A\in\mathbb{R}^{n,n}$, $k\in\mathbb{N}$, $x_0,b\in\mathbb{R}^{n}$
    \ENSURE $x^*$ approximate solution to $Ax=b$
    \STATE $\widehat{p}=r_0=b-Ax_0$,\quad $\widehat{v}=A\widehat{p}$
    \STATE $p_0=\frac{\widehat{p}}{\lVert \widehat{v}\rVert}$,\quad $v_0=\frac{\widehat{v}}{\lVert \widehat{v}\rVert}$
    \STATE $j=0$
    \WHILE{\nc}
        \STATE $\alpha_j=\langle v_j,r_j\rangle$ 
        \STATE $x_{j+1} = x_j + \alpha_j p_j$,\quad $r_{j+1} = r_j - \alpha_j v_j$
        \STATE $\widehat{p} = r_{j+1}$, $\widehat{v}=A\widehat{p}$
        \FOR{$i=j_k:j$}
            \STATE $\beta_i=\langle \widehat{v},v_i\rangle$ 
            \STATE $\widehat{p} = \widehat{p} - \beta_ip_i$,\quad $\widehat{v}=\widehat{v} - \beta_iv_i$
        \ENDFOR
        \STATE $p_{j+1}=\frac{\widehat{p}}{\lVert \widehat{v}\rVert}$,\quad $v_{j+1}=\frac{\widehat{v}}{\lVert \widehat{v}\rVert}$
        \STATE $j=j+1$
    \ENDWHILE
    \RETURN $x^*=x_j$
    \end{algorithmic}
\end{algorithm}
\end{minipage}
\begin{minipage}{.495\textwidth}
\begin{algorithm}[H]
    \caption{nlGCR($k$) for $f(x)=0$  \ \cite{He2024}}   
    \begin{algorithmic}[1]
    \setstretch{1.05}
    \REQUIRE $x_0\in\mathbb{R}^{n}$, $k\in\mathbb{N}$, $f,J_f$
    \ENSURE{$x^*$ approximate solution to $f(x)=0$}
    \STATE $\widehat{p}=r_0=-f(x_0)$,\quad $\widehat{v}=J_f(x_0)\widehat{p}$
    \STATE $p_0=\frac{\widehat{p}}{\lVert \widehat{v}\rVert}$,\quad $v_0=\frac{\widehat{v}}{\lVert \widehat{v}\rVert}$
    \STATE $j=0$
    \WHILE{\nc}
        \STATE $y_j={V_j}^Tr_j$ 
        \STATE $x_{j+1} = x_j +  P_jy_j$,\quad $r_{j+1} = -f(x_{j+1})$ 
        \STATE $\widehat{p} = r_{j+1}$, $\widehat{v}=J_f(x_{j+1})\widehat{p}$ 
        \FOR{$i=j_k:j$}
            \STATE $\beta_i=\langle \widehat{v},v_i\rangle$
            \STATE $\widehat{p} = \widehat{p} - \beta_ip_i$,\quad $\widehat{v}=\widehat{v} - \beta_iv_i$
        \ENDFOR
        \STATE $p_{j+1}=\frac{\widehat{p}}{\lVert \widehat{v}\rVert}$,\quad $v_{j+1}=\frac{\widehat{v}}{\lVert \widehat{v}\rVert}$ 
        \STATE $j=j+1$
    \ENDWHILE
    \RETURN $x^*=x_j$
\end{algorithmic}
\end{algorithm}
\end{minipage}\\
\begin{minipage}{.495\textwidth}
\begin{algorithm}[H]
    \caption{GMRESR($m,k$) \cite{DeS1996,Van1994}}  
    \label{alg:GMRESR}
    \begin{algorithmic}[1]
    \setstretch{1.1}
    \REQUIRE{$A\in\mathbb{R}^{n,n}$, $m,k\in\mathbb{N}$, $x_0,b\in\mathbb{R}^{n}$}
    \ENSURE{$x^*$ approximate solution to $Ax=b$}
    \STATE $r_0=b-Ax_0$
    \STATE $\widehat{p}=\texttt{GMRES}(A,r_0,m)$,\quad $\widehat{v}=A\widehat{p}$
    \STATE $p_0=\frac{\widehat{p}}{\lVert \widehat{v}\rVert}$,\quad $v_0=\frac{\widehat{v}}{\lVert \widehat{v}\rVert}$
    \STATE $j=0$
    \WHILE{\nc}
        \STATE $\alpha_j=\langle v_j,r_j\rangle$
        \STATE $x_{j+1} = x_j + \alpha_j p_j$,\quad $r_{j+1} = r_j - \alpha_j v_j$
        \STATE $\widehat{p} = \texttt{GMRES}(A,r_{j+1},m)$
        \STATE $\widehat{v}=A\widehat{p}$
        \FOR{$i=j_k:j$}
            \STATE $\beta_i=\langle \widehat{v},v_i\rangle$
            \STATE $\widehat{p} = \widehat{p} - \beta_ip_i$,\quad $\widehat{v}=\widehat{v} - \beta_iv_i$
        \ENDFOR
        \STATE $p_{j+1}=\frac{\widehat{p}}{\lVert \widehat{v}\rVert}$,\quad $v_{j+1}=\frac{\widehat{v}}{\lVert \widehat{v}\rVert}$
        \STATE $j=j+1$\\
    \ENDWHILE
    \RETURN $x^*=x_j$    
    \end{algorithmic}
\end{algorithm}
\end{minipage}
\begin{minipage}{.495\textwidth}
\begin{algorithm}[H]
    \caption{nlGMRESR($m,k$)}
    \label{alg:nlGMRESR}
    \begin{algorithmic}[1]
    \setstretch{1.1}
    \REQUIRE{$x_0\in\mathbb{R}^{n}$, $m,k\in\mathbb{N}$, $f,J_f$}
    \ENSURE{$x^*$ approximate solution to $f(x)=0$}
    \STATE $r_0=-f(x_0)$
    \STATE $\widehat{p}=\texttt{GMRES}(J_f(x_0),r_0,m)$,\quad $\widehat{v}=J_f(x_0)\widehat{p}$
    \STATE $p_0=\frac{\widehat{p}}{\lVert \widehat{v}\rVert}$,\quad $v_0=\frac{\widehat{v}}{\lVert \widehat{v}\rVert}$
    \STATE $j=0$
    \WHILE{\nc}
        \STATE $y_j= {V_j}^Tr_j$
        \STATE $x_{j+1} = x_j +  P_jy_j$,\quad $r_{j+1} = -f(x_{j+1})$
        \STATE $\widehat{p} = \texttt{GMRES}(J_f(x_{j+1}),r_{j+1},m)$ \label{line:gmressolve}
        \STATE $\widehat{v}=J_f(x_{j+1})\widehat{p}$ \label{line:gmresupdate}
        \FOR{$i=j_k:j$}
            \STATE $\beta_i=\langle \widehat{v},v_i\rangle$
            \STATE $\widehat{p} = \widehat{p} - \beta_ip_i$,\quad $\widehat{v}=\widehat{v} - \beta_iv_i$
        \ENDFOR
        \STATE $p_{j+1}=\frac{\widehat{p}}{\lVert \widehat{v}\rVert}$,\quad $v_{j+1}=\frac{\widehat{v}}{\lVert \widehat{v}\rVert}$
        \STATE $j=j+1$
    \ENDWHILE
    \RETURN $x^*=x_j$
\end{algorithmic}
\end{algorithm}
\end{minipage}\\
\begin{minipage}{.495\textwidth}
\begin{algorithm}[H]
    \caption{GCRO($m,k$) \cite{DeS1996,Des1999}} 
    \label{alg:GCRO}
\begin{algorithmic}[1]
\setstretch{1.1}
    \REQUIRE{$A\in\mathbb{R}^{n,n}$, $m,k\in\mathbb{N}$, $x_0,b\in\mathbb{R}^{n}$}
    \ENSURE{$x^*$ approximate solution to $Ax=b$}
    \STATE $r_0=b-Ax_0$
    \STATE $\widehat{p}=\texttt{GMRES}(A,r_0,m)$,\quad $\widehat{v}=A\widehat{p}$
    \STATE $p_0=\frac{\widehat{p}}{\lVert \widehat{v}\rVert}$,\quad $v_0=\frac{\widehat{v}}{\lVert \widehat{v}\rVert}$
    \STATE $j=0$
    \WHILE{\nc}
        \STATE $\alpha_j=\langle v_j,r_j\rangle$
        \STATE $x_{j+1} = x_j + \alpha_j p_j$,\quad $r_{j+1} = r_j - \alpha_j v_j$
        \STATE $\widehat{p} = \texttt{GMRES}(\orth{A}{V_j},r_{j+1},m)$
        \STATE $\widehat{v}=A\widehat{p}$
        \FOR[\textcolor{blue}{Optional}]{$i=j_k:j$} 
            \STATE $\beta_i=\langle \widehat{v},v_i\rangle$
            \STATE $\widehat{p} = \widehat{p} - \beta_ip_i$,\quad $\widehat{v}=\widehat{v} - \beta_iv_i$
        \ENDFOR
        \STATE $p_{j+1}=\frac{\widehat{p}}{\lVert \widehat{v}\rVert}$,\quad $v_{j+1}=\frac{\widehat{v}}{\lVert \widehat{v}\rVert}$
        \STATE $j=j+1$
    \ENDWHILE
    \RETURN $x^*=x_j$     
\end{algorithmic}
\end{algorithm}
\end{minipage}
\begin{minipage}{.495\textwidth}
\begin{algorithm}[H]
    \caption{nlGCRO($m,k$)}
    \label{alg:nlGCRO}
    \begin{algorithmic}[1]
    \setstretch{1.1}
    \REQUIRE{$x_0\in\mathbb{R}^{n}$, $m,k\in\mathbb{N}$, $f,J_f$}
    \ENSURE{$x^*$ approximate solution to $f(x)=0$}
    \STATE $r_0=-f(x_0)$
    \STATE  $\widehat{p}=\texttt{GMRES}(J_f(x_0),r_0,m)$,\quad $\widehat{v}=J_f(x_0)\widehat{p}$
    \STATE $p_0=\frac{\widehat{p}}{\lVert \widehat{v}\rVert}$,\quad $v_0=\frac{\widehat{v}}{\lVert \widehat{v}\rVert}$
    \STATE $j=0$
    \WHILE{\nc}
        \STATE $y_j= {V_j}^Tr_j$
        \STATE $x_{j+1} = x_j +  P_jy_j$,\quad $r_{j+1} = -f(x_{j+1})$
        \STATE $\widehat{p} = \texttt{GMRES}(\orth{J_f(x_{j+1})}{V_j},\widetilde{r}_{j+1},m)$ \label{line:gcrosolve}
        \STATE $\widehat{v}=J_f(x_{j+1})\widehat{p}$
        \FOR[\textcolor{blue}{Optional}]{$i=j_k:j$}
            \STATE $\beta_i=\langle \widehat{v},v_i\rangle$
            \STATE $\widehat{p} = \widehat{p} - \beta_ip_i$,\quad $\widehat{v}=\widehat{v} - \beta_iv_i$
        \ENDFOR
        \STATE $p_{j+1}=\frac{\widehat{p}}{\lVert \widehat{v}\rVert}$,\quad $v_{j+1}=\frac{\widehat{v}}{\lVert \widehat{v}\rVert}$
        \STATE $j=j+1$
    \ENDWHILE
    \RETURN $x^*=x_j$ 
    \end{algorithmic}
\end{algorithm}
\end{minipage}\\
\begin{minipage}{.495\textwidth}
\begin{algorithm}[H]
    \caption{LGMRES($m,k$) \cite{Baker2005,HicZ2010} }  
    \label{alg:LGMRES}       
\begin{algorithmic}[1]
    \setstretch{1.04}
    \REQUIRE{$A\in\mathbb{R}^{n,n}$, $m,k\in\mathbb{N}$, $x_0,b\in\mathbb{R}^{n}$}
    \ENSURE{$x^*$ approximate solution to $Ax=b$}
    \STATE $r_0=b-Ax_0$
    \STATE $p_0=\texttt{GMRES}(A,r_0,m+k)$
    \STATE $v_0=Ap_0$
    \STATE $j=0$
    \WHILE{\nc}
        \STATE \phantom{$y_j = V_j^+r_j$}
        \STATE $x_{j+1} = x_j + p_j$,\quad $r_{j+1} = r_j - v_j$
        \STATE $p_{j+1} = \texttt{AGMRES}(A,r_{j+1},m,k,P_j)$
        \STATE $v_{j+1}=Ap_{j+1}$
        \STATE $j=j+1$
    \ENDWHILE
    \RETURN $x^*=x_j$
\end{algorithmic}
\end{algorithm}
\end{minipage}
\begin{minipage}{.495\textwidth}
\begin{algorithm}[H]
\begin{algorithmic}[1]
    \setstretch{1.05}
    \REQUIRE{$x_0\in\mathbb{R}^{n}$, $m,k\in\mathbb{N}$, $f,J_f$}
    \ENSURE{$x^*$ approximate to $f(x)=0$}
    \STATE $r_0=-f(x_0)$
    \STATE $p_0=\texttt{GMRES}(J_f(x_0),r_0,m+k)$
    \STATE $v_0=J_f(x_0)p_0$
    \STATE $j=0$
    \WHILE{\nc}
        \STATE $y_j = V_j^+r_j$  \label{line:linesearchlgmres}
        \STATE $x_{j+1} = x_j + P_jy_j$,\quad $r_{j+1} = -f(x_{j+1})$
        \STATE $p_{j+1}\!=\!\texttt{AGMRES}(J_f(x_{j+1}),r_{j+1},m,k,P_j)$
        \STATE $v_{j+1}=J_f(x_{j+1})p_{j+1}$ \label{line:updatev}
        \STATE $j=j+1$
    \ENDWHILE
    \RETURN $x^*=x_j$
    \caption{nlLGMRES($m,k$)}
    \label{alg:nlLGMRES}       
\end{algorithmic}
\end{algorithm}
\end{minipage}\medskip\\
Since the two algorithms presented above do not directly fit the \textit{nlKrylov} framework, we suggest to modify the nonlinear algorithm to include outer orthogonalization. The resulting algorithm nlLGMRESO($m,k$) does not simplify to \Cref{alg:LGMRES} for linear problems but has significant stability improvements. It also allows for a comprehensive convergence analysis within the \textit{nlKrylov} framework. Whenever nlLGMRES($m,k$) is referenced in the paper, we refer to \Cref{alg:nlLGMRESO}.
\begin{algorithm}[H]    
    \caption{nlLGMRESO($m,k$) for $f(x)=0$}
    \label{alg:nlLGMRESO}
    \begin{algorithmic}[1]
    \setstretch{0.98}
    \REQUIRE{$x_0\in\mathbb{R}^{n}$, $m,k\in\mathbb{N}$, $f,J_f$}
    \ENSURE{$x^*$ approximate solution to $f(x)=0$}
    \STATE $r_0=-f(x_0)$
    \STATE $\widehat{p}=\texttt{GMRES}(J_f(x_0),r_0,m+k)$,\quad $\widehat{v}=J_f(x_0)\widehat{p}$
    \STATE $p_0=\frac{\widehat{p}}{\lVert \widehat{v}\rVert}$,\quad $v_0=\frac{\widehat{v}}{\lVert \widehat{v}\rVert}$
    \STATE $j=0$
    \WHILE{\nc}
        \STATE $y_j = V_j^Tr_j$ 
        \STATE $x_{j+1} = x_j + P_jy_j$,\quad $r_{j+1} = -f(x_{j+1})$
        \STATE $\widehat{p} = \texttt{AGMRES}(J_f(x_{j+1}),r_{j+1},m,k,P_j)$,\quad $\widehat{v}=J_f(x_{j+1})\widehat{p}$
        \FOR{$i=j_k:j$}
            \STATE $\beta_i=\langle \widehat{v},v_i\rangle$,\quad $\widehat{p} = \widehat{p} - \beta_ip_i$,\quad $\widehat{v}=\widehat{v} - \beta_iv_i$
        \ENDFOR
        \STATE $p_{j+1}=\frac{\widehat{p}}{\|\widehat{v}\|}$,\quad $v_{j+1}=\frac{\widehat{v}}{\|\widehat{v}\|}$
        \STATE $j=j+1$
    \ENDWHILE
    \RETURN $x^*=x_j$    
    \end{algorithmic}
\end{algorithm}
\subsection{Global nlKrylov methods} \label{sec:globalnlkrylov}
All strategies introduced in \Cref{sec:nlKrylov} of the main text can be easily modified to solve nonlinear matrix equations of the form 
\begin{equation}
F(X)=0, \quad F:\RR^{n,p}\rightarrow\RR^{n,p},~p\leq n. \label{eq:matnonlinear}   
\end{equation}
For problems of this form, the standard multivariate Newton iteration can be generalized by replacing the system \eqref{eq:Newton} involving the Jacobian $J_f(x_j)$ with a linear matrix equation 
\begin{equation}
    L_F(X_j,\Delta{X}_j)=-F(X_j),  \label{eq:matnewton}  
\end{equation}
where $L_F(X_j):\RR^{n,p}\rightarrow\RR^{n,p}$ is the Fréchet derivative of $F$ at $X_j$. In Section \ref{sec:glkrylov}, we briefly recall the idea of global Krylov solvers for linear operator equations and afterwards, in Section \ref{sec:glnlkrylov}, we explain how to combine these ideas and the results on \textit{nlKrylov} methods introduced in this paper to acquire global \textit{nlKrylov} methods for nonlinear matrix equations. For the remainder of this section, we use capital letters to denote matrices and calligraphic letters the denote the column wise concatenation of multiple matrices into a block matrix. Specifically, if $U_1,\dots,U_j$ are matrices of size $n\times p$, we define the corresponding block matrix $\mathcal{U}_j=\left[U_1,U_2,\dots,U_j\right]\in\RR^{n,p\cdot j}$. On such block matrices, we define the matrix-vector product
\begin{equation}
 \mathcal{U}_j\dprod \gamma = \mathcal{U}_j(\gamma\otimes I_p)=\sum_{i=1}^j\gamma_iU_i,\quad\gamma=(\gamma_i)_{i=1,\dots,j}\in\RR^j. \label{eq:blockmatvec}  
\end{equation}
Finally, we denote by $\langle\cdot,\cdot\rangle$ and $\lVert\cdot\rVert$ the Frobenius inner product and norm, respectively. For a matrix $W\in\RR^{n,p}$ and a block matrix $\mathcal{U}_j\in\RR^{n,p\cdot j}$, we define
$$\langle\mathcal{U}_j,W\rangle = \begin{bmatrix}\langle U_1,W\rangle\\\vdots\\\langle U_j,W\rangle\end{bmatrix}\in\RR^j \quad \text{and}\quad \langle \mathcal{U}_j,\mathcal{U}_j\rangle=\begin{bmatrix}\langle U_1,U_1\rangle & \cdots & \langle U_1,U_j\rangle\\ \vdots & \ddots & \vdots \\ \langle U_j,U_1\rangle & \cdots & \langle U_j,U_j\rangle
\end{bmatrix}\in\RR^{j,j}.$$
We call the block matrix $\mathcal{U}_j$ \emph{F-orthogonal}, if $\langle\mathcal{U}_j,\mathcal{U}_j\rangle = I_j$. Additionally, for a linear operator $\mathcal{A}:\RR^{n,p}\rightarrow\RR^{n,p}$, and an F-orthogonal block matrix $\mathcal{U}_j$, we denote by $\orth{\mathcal{A}}{\mathcal{U}_j}:=(I-\mathcal{U}_j\dprod \langle\mathcal{U}_j,\cdot\rangle)\mathcal{A}$ the operator such that
$$\orth{\mathcal{A}}{\mathcal{U}_j}(X):=((I-\mathcal{U}_j\dprod\langle\mathcal{U}_j,\cdot\rangle)\mathcal{A})(X) = \mathcal{A}(X)-\mathcal{U}_j\dprod\langle\mathcal{U}_j,\mathcal{A}(X)\rangle,$$
i.e., the operator that orthogonalizes $\mathcal{A}(X)$ against $\mathcal{U}_j$ with respect to the Frobenius inner product. It is easy to see that this is the matrix-valued equivalent to the GCRO-operator $\orth{A}{V_j}=(I-V_jV_j^T)A$ introduced in equation \eqref{eq:gcrosolve}, hence, we use the same notation.
\subsubsection{Global linear methods} \label{sec:glkrylov}
Global Krylov subspace methods \cite{JbiMS1999} have been introduced in recent years as an efficient iterative tool for solving linear matrix equations 
\begin{equation}
\mathcal{A}(X)=B,\quad \mathcal{A}:\RR^{n,p}\rightarrow\RR^{n,p},\quad B\in\RR^{n,p},   \label{eq:linmatrix} 
\end{equation}
where usually $p\ll n$ and $\mathcal{A}$ is a linear operator that is inexpensive to evaluate. Classical examples include Sylvester equations 
$$A_1X+XA_2=B,$$
where $A_1\in\RR^{n,n}$ and $A_2\in\RR^{p,p}$ \cite{ZadTW2019,WanG2007} are large and sparse, or linear systems with multiple right hand sides \cite{TouK2006}.
Let us recall that, given an initial residual $R_0=B-\mathcal{A}(X_0)$, $X_0\in\RR^{n,p}$, global GMRES \cite{JbiMS1999} seeks to find a solution of the linear equation \eqref{eq:linmatrix} in the block Krylov subspace 
$$\mathcal{K}^\Box_j(\mathcal{A},R_0)=\mathrm{blockspan}\left\lbrace R_0,\mathcal{A}(R_0),\dots,\mathcal{A}^{j-1}(R_0)\right\rbrace.$$
After $m$ steps of global GMRES, we obtain the global Arnoldi relation 
$$\mathcal{A}\mathcal{Q}_m=\mathcal{Q}_{m+1}\dprod \underline{H}_m,$$
where $\mathcal{Q}_{m+1}\in\RR^{n,(m+1)\cdot p}$ is F-orthogonal and $\underline{H}_m\in\RR^{m+1,m}$ is upper Hessenberg with an extra row. Note that the size of $\underline{H}_m$ depends only on the size $m$ of the Krylov subspace and is not affected by the block size $p$ of the problem. The global Krylov methods rely on the fact that every linear operator $\mathcal{A}:\RR^{n,p}\rightarrow\RR^{n,p}$ has a matrix representation $\mathcal{M}_\mathcal{A}\in\RR^{n\cdot p,n\cdot p}$ such that
$$\mathrm{vec}(\mathcal{A}(X))=\mathcal{M}_\mathcal{A}\cdot\mathrm{vec}(X).$$
Hence, solving \eqref{eq:linmatrix} by $m$ steps of global GMRES is mathematically equivalent to solving the vectorized equation 
$$\mathcal{M}_\mathcal{A}x=b,\quad \mbox{where} \quad x=\mathrm{vec}(X),~b=\mathrm{vec}(B)\in\RR^{n\cdot p},$$
by standard GMRES. Although, to the best of our knowledge it has not been done, GCR-based nested Krylov subspace methods can also be extended to linear matrix equations by using the Frobenius inner product and Frobenius norm in the modified Gram-Schmidt process (lines \ref{line:inner1} and \ref{line:inner2} of \Cref{alg:GCR}). Analogously, GMRESR, GCRO and LGMRES can be adapted by adding global GMRES as the inner solver. As shown in
\cite{Wer2024}, global Krylov methods can serve as efficient inner solvers for \eqref{eq:matnewton} in the context of matrix-valued Newton methods. In the remainder of this section, we extend the nonlinear Krylov methods discussed in \Cref{sec:nlKrylov} to matrix-valued problems.
\subsubsection{Moving to global nonlinear methods} \label{sec:glnlkrylov}
Using the notation introduced in \Cref{sec:nlKrylov} and the ideas of global Krylov subspace methods, we can easily extend \Cref{alg:nlGCR}, \Cref{alg:nlGMRESR}, \Cref{alg:nlGCRO} and \Cref{alg:nlLGMRESO} to the case of nonlinear matrix equations. For this sake, we will denote by $\mathcal{P}_j,\mathcal{V}_j\in\RR^{n,n_j\cdot p}$ the truncated block basis matrices of global GCR, i.e., 
$$\mathcal{P}_j=[P_{j_k},\dots,P_j],\quad \mathcal{V}_j=[V_{j_k},\dots,V_j]=\mathcal{A}\mathcal{P}_j,$$
and require three minor modifications to the vector-valued \textit{nlKrylov} framework: First, we have to replace all inner products and norms in the base algorithms by their Frobenius counterparts. Second, matrix-vector multiplications have to be performed using the $\dprod$-operator \eqref{eq:blockmatvec} instead of the standard matrix-vector-product, e.g., to update $X_{j+1}=X_j + \mathcal{P}_j\dprod y_j$. Note that computing $y_j=V_j^Tr_j$ requires column-wise inner products instead of standard matrix-vector multiplication, i.e., $y_j=\langle\mathcal{V}_j,R_j\rangle\in\RR^{n_j}$. Finally, all matrix-vector multiplications involving the Jacobian, e.g., $J_f(x_j)z$, $z\in\RR^n$, have to be replaced by evaluations of the Fr\'echet derivative $L_F(X_j,Z)$ of $F$ at $X_j\in\RR^{n,p}$ in the direction of $Z\in\RR^{n,p}$. Hence, global Krylov methods can also be used to efficiently solve the inner equation 
\begin{equation}
  L_F(X_j,\widehat{P})=R_j \label{eq:glrefinedlocal}  
\end{equation}
for $\widehat{P}$. Section \ref{sec:changematrix} of the main text contains comments on the efficient implementation of these modifications using {\sc MATLAB}'s high level linear algebra routines.\\
Using the three modifications, the general framework from Section \ref{sec:nlkrylovframework} can in principle be extended to any scheme $\mathcal{SR}_j(R_j,L_F(X_j))$ that employs $R_j$ together with $L_F(X_j,\cdot):\RR^{n,p}\to\RR^{n,p}$ for the refinement of the local linear model \eqref{eq:glrefinedlocal} in order to solve nonlinear matrix equations \eqref{eq:matnonlinear}, thereby encompassing the broader class of \textit{global nlKrylov} methods. A template for \textit{global nlKrylov} methods can be found in \Cref{alg:GLnlKrylov}, global versions of the \textit{nlKrylov} methods covered in Section \ref{sec:basenlkrylov} are stated in Algorithms \ref{alg:GLnlGCR}--\ref{alg:GLnlLGMRESO}. Note that truncation in this case means we drop matrices $P_{j-k}$ from the block matrix $\mathcal{P}_j$. 
\begin{algorithm}[H]    
    \caption{Unifying global nlKrylov($k$) methods for $F(X)=0$}
    \label{alg:GLnlKrylov}
    \begin{algorithmic}[1]
    \REQUIRE{$X_0\in\mathbb{R}^{n,p}$, $m,k\in\mathbb{N}$, $F,L_F,\mathcal{SR}$}
    \ENSURE{$X^*$ approximate solution to $F(X)=0$}
    \STATE $R_0=-F(X_0)$
    \STATE $\widehat{P}=\mathcal{SR}_0(R_0,L_F(X_0))$,\quad $\widehat{V}=L_F(X_0,\widehat{P})$
    \STATE $P_0=\frac{\widehat{P}}{\lVert \widehat{V}\rVert}$,\quad $V_0=\frac{\widehat{V}}{\lVert \widehat{V}\rVert}$
    \STATE $j=0$
    \WHILE{\nc}
        \STATE $y_j = \langle \mathcal{V}_j,R_j\rangle$ 
        \STATE $X_{j+1} = X_j + \mathcal{P}_j\dprod y_j$,\quad $R_{j+1} = -F(X_{j+1})$
        \STATE $\widehat{P} = \mathcal{SR}_{j+1}(R_{j+1},L_F(X_{j+1}))$,\quad $\widehat{V}=L_F(X_{j+1},\widehat{P})$
        \FOR[Optional for GL-nlGCRO($m,k$)]{$i=j_k:j$}
            \STATE $\beta_i=\langle \widehat{V},V_i\rangle$,\quad $\widehat{P} = \widehat{P} - \beta_iP_i$,\quad $\widehat{V}=\widehat{V} - \beta_iV_i$
        \ENDFOR
        \STATE $P_{j+1}=\frac{\widehat{P}}{\|\widehat{V}\|}$,\quad $V_{j+1}=\frac{\widehat{V}}{\|\widehat{V}\|}$
        \STATE $j=j+1$
    \ENDWHILE
    \RETURN $X^*=X_j$    
    \end{algorithmic}
\end{algorithm}
\vspace{-.75cm}
\noindent
\begin{minipage}{.495\textwidth}
\begin{algorithm}[H]
    \caption{GL-nlGCR($k$)}   
    \label{alg:GLnlGCR}
    \begin{algorithmic}[1]
    \setstretch{1.03}
    \REQUIRE{$X_0\in\mathbb{R}^{n,p}$, $k\in\mathbb{N}$, $F,L_F$}
    \ENSURE{$X^*$ approximate solution to $F(X)=0$}
    \STATE $R_0=-F(X_0)$
    \STATE $\widehat{P}=R_0$
    \STATE $\widehat{V}=L_F(X_0,\widehat{P})$ 
    \STATE $P_0=\frac{\widehat{P}}{\lVert \widehat{V}\rVert}$,\quad $V_0=\frac{\widehat{V}}{\lVert \widehat{V}\rVert}$
    \STATE $j=0$
    \WHILE{\nc}
        \STATE $y_j=\langle \mathcal{V}_j,R_j\rangle$
        \STATE $X_{j+1} = X_j +  \mathcal{P}_j\dprod y_j$
        \STATE $R_{j+1} = -F(X_{j+1})$
        \STATE $\widehat{P} = R_{j+1}$ 
        \STATE $\widehat{V}=L_F(X_{j+1},\widehat{P})$
        \FOR{$i=j_k:j$}
            \STATE $\beta_i=\langle \widehat{V},V_i\rangle$
            \STATE $\widehat{P} = \widehat{P} - \beta_iP_i$,\quad $\widehat{V}=\widehat{V} - \beta_iV_i$
        \ENDFOR
        \STATE $P_{j+1}=\frac{\widehat{P}}{\lVert \widehat{V}\rVert}$,\quad $V_{j+1}=\frac{\widehat{V}}{\lVert \widehat{V}\rVert}$
        \STATE $j=j+1$
    \ENDWHILE
    \RETURN $X^*=X_j$      
    \end{algorithmic}
\end{algorithm}
\end{minipage}
\begin{minipage}{.495\textwidth}
\begin{algorithm}[H]    
\caption{GL-nlGMRESR($m,k$)}   
\label{alg:GLnlGMRESR}    
\begin{algorithmic}[1]
\setstretch{1.03}
    \REQUIRE{$X_0\in\mathbb{R}^{n,p}$, $m,k\in\mathbb{N}$, $F,L_F$}
    \ENSURE{$X^*$ approximate solution to $F(X)=0$}
    \STATE $R_0=-F(X_0)$
    \STATE $\widehat{P}=\texttt{GLGMRES}(L_F(X_0),R_0,m)$
    \STATE $\widehat{V}=L_F(X_0,\widehat{P})$
    \STATE $P_0=\frac{\widehat{P}}{\lVert \widehat{V}\rVert}$,\quad $V_0=\frac{\widehat{V}}{\lVert \widehat{V}\rVert}$
    \STATE $j=0$
    \WHILE{\nc}
        \STATE $y_j=\langle \mathcal{V}_j,R_j\rangle$
        \STATE $X_{j+1} = X_j +  \mathcal{P}_j\dprod y_j$
        \STATE $R_{j+1} = -F(X_{j+1})$
        \STATE $\widehat{P} = \texttt{GLGMRES}(L_F(X_{j+1}),R_{j+1},m)$
        \STATE $\widehat{V}=L_F(X_{j+1},\widehat{P})$
        \FOR{$i=j_k:j$}
            \STATE $\beta_i=\langle \widehat{V},V_i\rangle$
            \STATE $\widehat{P} = \widehat{P} - \beta_iP_i$,\quad $\widehat{V}=\widehat{V} - \beta_iV_i$
        \ENDFOR
        \STATE $P_{j+1}=\frac{\widehat{P}}{\lVert \widehat{V}\rVert}$,\quad $V_{j+1}=\frac{\widehat{V}}{\lVert \widehat{V}\rVert}$
        \STATE $j=j+1$
    \ENDWHILE
    \RETURN $X^*=X_j$
\end{algorithmic}
\end{algorithm}
\end{minipage}\\
\begin{minipage}{.495\textwidth}
\begin{algorithm}[H]
    \caption{GL-nlGCRO($m,k$)}   
    \label{alg:GLnlGCRO}
    \begin{algorithmic}[1]
    \setstretch{1.06}
    \REQUIRE{$X_0\in\mathbb{R}^{n,p}$, $m,k\in\mathbb{N}$, $F,L_F$}
    \ENSURE{$X^*$ approximate solution to $F(X)=0$}
    \STATE $R_0=-F(X_0)$
    \STATE $\widehat{P}=\texttt{GLGMRES}(L_F(X_0),R_0,m)$
    \STATE $\widehat{V}=L_F(X_0,\widehat{P})$
    \STATE $P_0=\frac{\widehat{P}}{\lVert \widehat{V}\rVert}$,\quad $V_0=\frac{\widehat{V}}{\lVert \widehat{V}\rVert}$
    \STATE $j=0$
    \WHILE{\nc}
        \STATE $y_j=\langle \mathcal{V}_j,R_j\rangle$
        \STATE $X_{j+1} = X_j +  \mathcal{P}_j\dprod y_j$
        \STATE $R_{j+1} = -F(X_{j+1})$
        \STATE $\widehat{P} = \texttt{GLGMRES}(\orth{L_F(X_{j+1})}{\mathcal{V}_j},\widetilde{R}_{j+1},m)$
        \STATE $\widehat{V}=L_F(X_{j+1},\widehat{P})$
        \FOR{$i=j_k:j$}
            \STATE $\beta_i=\langle \widehat{V},V_i\rangle$
            \STATE $\widehat{P} = \widehat{P} - \beta_iP_i$,\quad $\widehat{V}=\widehat{V} - \beta_iV_i$
        \ENDFOR
        \STATE $P_{j+1}=\frac{\widehat{P}}{\lVert \widehat{V}\rVert}$,\quad $V_{j+1}=\frac{\widehat{V}}{\lVert \widehat{V}\rVert}$
        \STATE $j=j+1$
    \ENDWHILE
    \RETURN $X^*=X_j$       
    \end{algorithmic}
\end{algorithm}
\end{minipage}
\begin{minipage}{.495\textwidth}
\begin{algorithm}[H]
    \caption{GL-nlLGMRESO($m,k$)}   
    \label{alg:GLnlLGMRESO}
\begin{algorithmic}[1]
\setstretch{1.06}
    \REQUIRE{$X_0\in\mathbb{R}^{n,p}$, $m,k\in\mathbb{N}$, $F,L_F$}
    \ENSURE{$X^*$ approximate solution to $F(X)=0$}
    \STATE $R_0=-F(X_0)$
    \STATE $\widehat{P}=\texttt{GLGMRES}(L_F(X_0),R_0,m+k)$
    \STATE $\widehat{V}=L_F(X_0,\widehat{P})$
    \STATE $P_0=\frac{\widehat{P}}{\lVert \widehat{V}\rVert}$,\quad $V_0=\frac{\widehat{V}}{\lVert \widehat{V}\rVert}$
    \STATE $j=0$
    \WHILE{\nc}
        \STATE $y_j=\langle \mathcal{V}_j,R_j\rangle$
        \STATE $X_{j+1} = X_j +  \mathcal{P}_j\dprod y_j$
        \STATE $R_{j+1} = -F(X_{j+1})$
        \STATE \resizebox{.83\linewidth}{!}{$\widehat{P}=\texttt{GLAGMRES}(L_F(X_{j+1}),R_{j+1},m,k,\mathcal{P}_j)$}
        \STATE $\widehat{V}=L_F(X_{j+1},\widehat{P})$
        \FOR{$i=j_k:j$ }
            \STATE $\beta_i=\langle \widehat{V},V_i\rangle$
            \STATE $\widehat{P} = \widehat{P} - \beta_iP_i$,\quad $\widehat{V}=\widehat{V} - \beta_iV_i$
        \ENDFOR
        \STATE $P_{j+1}=\frac{\widehat{P}}{\lVert \widehat{V}\rVert}$,\quad $V_{j+1}=\frac{\widehat{V}}{\lVert \widehat{V}\rVert}$
        \STATE $j=j+1$
    \ENDWHILE
    \RETURN $X^*=X_j$   
\end{algorithmic}
\end{algorithm}
\end{minipage}\\

\medskip
\noindent
Since we usually require relatively low $m$ and $k$ within nested \textit{nlKrylov} methods, they are particularly appealing in larger matrix equations where storage requirements impose a major challenge and provide a competitive alternative to inexact Newton or simple fixed-point schemes.

\subsection{Additional algorithms} \label{sec:addalg}
In this section, we briefly state two additional algorithms referenced in the paper and used for the implementation. The first one is the augmented GMRES solver used within the (nl)LGMRES(O)-algorithm. It aims at finding a solution to the linear system $Ax=b$ over the augmented space $\mathcal{K}_m(A,r_0)\oplus \ran(P)$, where $r_0$ is the initial residual and $P\in\RR^{n,s}$ is an arbitrary augmentation space. The augmented space simplifies to the relation
$$A \left[Q_m, P\right]=\left[Q_{m+1}, V\right]\underline{H}_{m+s}$$
whenever $V=AP$ holds, as is within the linear GCR framework. For nlLGMRES($m,k$), we need to use \Cref{alg:AGMRES} as stated below since generally, we have $J_f(x_j)P_j\neq V_j$. This leads to $k$ additional function evaluations being required compared to nlGMRESR($m,k$). Whenever the linear update version of nlLGMRES-A is used, the additional inexactness introduced by using $v_{i-m_s}$ in line \ref{line:usev} can be justified since the linear model is already assumed to be sufficiently accurate and we use this version to further reduce the necessary function evaluations.\\
The second algorithm we want to state is the Armijo-Line-Search algorithm used to determine damped updates of $x_{j+1}$ that improve the residual norm 
$$\|f(x_{j+1})\|=\|f(x_j+\alpha_jP_jy_j)\|.$$
The estimated angle $\zeta_j=r_j^TJ_f(x_{j})d_j$ is computed using the initial step length $\alpha_j^{(0)}$ and can be inexpensively refined during the iteration by using the reduced step length $\alpha_j^{(\ell)}$. As stated in the paper, we choose $\alpha_j^{(0)}$ depending on the number of line search steps taken in iteration $j-1$.
\begin{algorithm}[H]
\caption{AGMRES($m,k$) for $Ax=b$}
\label{alg:AGMRES}
\begin{algorithmic}[1]
    \REQUIRE{$A\in\mathbb{R}^{n,n}$, $m,k\in\mathbb{N}$, $x_0,b\in\mathbb{R}^{n}$, $P=[p_1,\dots,p_s]\in\RR^{n,s},s\leq k$}
    \ENSURE{$x^*$ approximate solution to $Ax=b$} 
    \STATE $r_0=b-Ax_0$,\quad $\beta=\lVert r_0\rVert$
    \STATE $q_1=\frac{r_0}{\beta}$,\quad $m_s=m+(k-s)$
        \FOR{$i=1:m+k$}
            \IF[Standard GMRES]{$i\leq m_s$}
                \STATE $w=Aq_i$ 
            \ELSE[Augment by P]
                \STATE $w=Ap_{i-m_s}$   \label{line:usev}  
            \ENDIF
            \FOR{$\ell=1:i$}
                \STATE $h_{\ell,i}=\langle w,q_\ell\rangle$,\quad $w=w-h_{\ell,i}q_\ell$ 
            \ENDFOR
            \STATE $h_{i+1,i}=\lVert w\rVert$,\quad $q_{i+1}=\frac{w}{h_{i+1,i}}$ 
        \ENDFOR
        \STATE Define $Z_{m+k}=\left[q_1,\dots,q_{m_s},P\right]\in\mathbb{R}^{n,(m+k)}$
        \STATE $\gamma_{m+k}=\mathrm{argmin}_{\gamma\in\mathbb{R}^{m+k}}\lVert \beta e_1-\underline{H}_{m+k}\gamma\rVert$
        \STATE $x_{m+k}=x_0+Z_{m+k}\gamma_{m+k}$
        \RETURN $x^*=x_{m+k}$    
\end{algorithmic}
\end{algorithm}
\vspace{-.65cm}
\begin{algorithm}[H]
\caption{Armijo-Goldstein-Linesearch}   
\label{alg:linesearch}
\begin{algorithmic}[1]
\setstretch{.98}
\REQUIRE{$f:\RR^n\to\RR^n$, $x_j,d_j,r_j\in\RR^n$, $c_1,\alpha_j^{(0)}>0$}
\ENSURE{$\alpha_j$ damping parameter such that $x_{j+1}=x_j+\alpha_jd_j$}
\STATE Compute $f_{j+1}^{(0)}=f(x_j+\alpha_j^{(0)}d_j)$ and estimate $\zeta_j^{(0)}=\frac{1}{\alpha_j^{(0)}}\langle r_j,f_{j+1}^{(0)}+r_j\rangle$ 
\IF{$\zeta_j^{(0)}<0$}
\STATE Set $d_j=-d_j$ and $\zeta_j^{(0)}=-\zeta_j^{(0)}$, recompute $f_{j+1}^{(0)}=f(x_j+\alpha_j^{(0)}d_j)$
\ENDIF
\STATE Set $\ell=0$
\WHILE{$\|f_{j+1}^{(\ell)}\|^2 > \|r_j\|^2 - c_1\alpha_j^{(\ell)} \zeta_j^{(\ell)}$} 
    \STATE Set $\alpha_j^{(\ell+1)}=\frac{\alpha_j^{(\ell)}}{2}$ and compute damped step $f_{j+1}^{(\ell+1)}=f(x_j+\alpha_j^{(\ell+1)}d_j)$
    \STATE If desired, refine $\zeta_j^{(\ell+1)}=\frac{1}{\alpha_j^{(\ell+1)}}\langle r_j,f_{j+1}^{(\ell+1)}+r_j\rangle$, else, just use $\zeta_j^{(\ell+1)}=\zeta_j^{(\ell)}$
    \STATE Set $\ell=\ell+1$
\ENDWHILE
\RETURN $\alpha_j=\alpha_j^{(\ell)}$
\end{algorithmic}
\end{algorithm}
\newpage

\section{Extended proof of \Cref{thm:singular}}
In this section, we complete the proof of \Cref{thm:singular} by presenting the full induction step that was omitted from the manuscript due to space constraints. The proof proceeds along the same lines as the argument for the initial step.
\begin{proof}[Proof of \Cref{thm:singular} (Continued)]
The base case ($j=0$) was established in the manuscript. Assume now that $j>0$ and that there exist $\rho_j,\gamma_j$ such that 
$x_j\in W(\rho_j,\gamma_j)\subset W(\rho,\gamma)$. Recall that 
$$\widetilde{x}_j:=x_j-x^*,\quad \rho_j:=\|\widetilde{x}_j\|,\quad \gamma_j:=K_1\rho_{j-1}(1+\gamma_{j-1})(\frac{1}{2}-K_0\gamma_{j-1})^{-1/2}.$$
Since $x_j\in W(\rho_j,\gamma_j)$, we may expand $f$ and $J_f$ about $x^*$
exactly as in the case $j=0$ (see \eqref{eq:expandf},\eqref{eq:expandJf} to obtain 
\begin{align}
    f(x_j)&=J_f(x^*)\widetilde{x}_j + \frac{1}{2}\left[A_1(x_j)+B_1(x_j)+C_1(x_j)+D_1(x_j)\right]\widetilde{x}_j + \mathcal{O}(\|\widetilde{x}_j\|^3), \label{eq:expandfj} \\
    J_f(x_j)&=J_f(x^*)+\left[A_1(x_j)+B_1(x_j)+C_1(x_j)+D_1(x_j)\right] + \mathcal{O}(\|\widetilde{x}_j\|^2). \label{eq:expandJfj}
\end{align}
Combining \eqref{eq:expandfj} and \eqref{eq:expandJfj} yields
\begin{equation}
    f(x_j)=J_f(x_j)\widetilde{x}_j-\frac{1}{2}\left[A_1(x_j)+B_1(x_j)+C_1(x_j)+D_1(x_j)\right]\widetilde{x}_j+\mathcal{O}(\|x_j\|^3). \label{eq:star} 
\end{equation}
Using \eqref{eq:expandfj}, the inclusion $x_j\in W(\rho_j,\gamma_j)$ and the identity $J_f(x^*)\widetilde{x}_j=J_f(x^*)P_\cX\widetilde{x}_j$, we obtain
\begin{align*}
    \|f(x_j)\|&\leq K\|P_\cX \widetilde{x}_j\|\leq \gamma_j K\|P_\cN \widetilde{x}_j\| \leq \gamma_j K\|\widetilde{x}_j\|,\\
    \|t_j\|&\leq c\|f(x_j)\|^2\leq c\gamma_j^2K^2\|\widetilde{x}_j\|^2. 
\end{align*}
Hence, $\|t_j\|=\gamma_j^2\mathcal{O}(\|\widetilde{x}_j\|^2)$. 
This estimate together with~\eqref{eq:star}, \Cref{lem:Jinv}, and the definition of $A_1,B_1,C_1$ and $D_1$ (see \eqref{eq:defA}--\eqref{eq:defD} yields
\begin{equation}
    \widetilde{x}_{j+1}=\widetilde{x}_j - J_f(x_j)^{-1}\left(f(x_j)-t_j\right)=\frac{1}{2}P_\cN \widetilde{x}_j + \gamma_j P_\cN\mathcal{O}(\|\widetilde{x}_j\|)+\gamma_j^2P_\cN\mathcal{O}(\|\widetilde{x}_j\|)+\mathcal{O}(\|\widetilde{x}_j\|^2). \label{eq:xjupdate}
\end{equation}
Applying $P_\cX$ and $P_\cN$ to \eqref{eq:xjupdate}, respectively, and taking norms shows that, for $\gamma_j$ sufficiently small, there exist constants $K_1>0$ and $K_0>0$ satisfying $K_0\gamma_j<\frac{1}{2}$ such that 
\begin{align}
    \|P_\cX \widetilde{x}_{j+1}\|&\leq K_1\|\widetilde{x}_j\|^2 \label{eq:showconv}\\
    \left(\frac{1}{2}-K_0\gamma_j\right)\|P_\cN\widetilde{x}_j\|&\leq\|P_\cN \widetilde{x}_{j+1}\|\leq\left(\frac{1}{2}+K_0\gamma_j\right)\|P_\cN\widetilde{x}_j\| \label{eq:boundPn}
\end{align}
Equation \eqref{eq:showconv} establishes \eqref{eq:ConvQuadTerm}. Moreover,
\begin{equation*}
    \|P_\cX \widetilde{x}_{j+1}\|\leq K_1\rho_j^2\leq K_1\rho_j(1+\gamma_j)\|P_\cN\widetilde{x}_j\|\leq K_1\rho_j(1+\gamma_j)(\frac{1}{2}-K_0\gamma_j)^{-1}\|P_\cN\widetilde{x}_{j+1}\|=\gamma_{j+1}\|P_\cN\widetilde{x}_{j+1}\|,
\end{equation*}
which implies that $x_{j+1}\in W(\rho_{j+1},\gamma_{j+1})$. 
Combining the upper bound in~\eqref{eq:boundPn} with \eqref{eq:showconv} gives
\begin{equation*}
    \rho_{j+1}\leq (\frac{1}{2}+K_0\gamma_j)\|P_\cN\widetilde{x}_j\| + K_1\rho_j^2\leq\left((\frac{1}{2}+K_0\gamma_j)(1+\gamma_j)^{-1} + K_1\rho_j\right)\rho_j    
\end{equation*}
Therefore, for $\gamma_j$ and $\rho_j$ sufficiently small, there exists $\tau\in(\frac{1}{2},1)$ such that $\gamma_{j+1}\leq \tau \gamma_j$ and $\rho_{j+1}\leq \tau \rho_j$. By induction, it follows that
$\gamma_j\rightarrow0$ and $\rho_j\rightarrow0$ as $j \rightarrow \infty$.
Consequently, $x_j\to x^*$ with $q$-linear convergence rate bounded by $\tau$. Moreover, from \eqref{eq:boundPn},
\[
\left(\frac12-K_0\gamma_j\right)
\le
\frac{\|P_{\mathcal N}\widetilde{x}_{j+1}\|}
     {\|P_{\mathcal N}\widetilde{x}_j\|}
\le
\left(\frac12+K_0\gamma_j\right),
\]
and, since \(\gamma_j\to0\), the squeeze theorem implies that
\begin{equation*}
    \lim_{j\to\infty} \frac{\|P_\cN \widetilde{x}_{j+1}\|}{\|P_\cN \widetilde{x}_j\|}=\frac{1}{2},  
\end{equation*}
This proves \eqref{eq:limitnull} and thus completes the proof.
\end{proof}
\section{Additional examples}
In this section, we present additional numerical results for examples not included in the main text, including a modified setup of the symmetric Bratu problem with stronger nonlinearity in Section \ref{sec:bratu6}, a nonlinear eigenvalue problem (NEP) in Section \ref{sec:nep} and a nonlinear eigenvector problem (NEPv) in Section \ref{sec:nepv}. Additionally, in Section \ref{sec:sensitivity} we provide a brief sensitivity analysis for selected \textit{nlKrylov} methods on the Bratu problem discussed in the main text to illustrate the dependency on $m$ and $k$. 
\subsection{The symmetric Bratu problem with stronger nonlinearity} \label{sec:bratu6}
In this section, we revisit the symmetric Bratu problem from Section \ref{sec:bratu} of the main text,
\begin{eqnarray*}
    \Delta u + \lambda e^u &=&0, \quad(y,z)\in\Omega,\\
    u &=&0, \quad(y,z)\in\partial\Omega,
\end{eqnarray*}
on the domain $\Omega=(0,1)\times(0,1)$, where now we choose $\lambda=6$ instead of $\lambda=0.5$ to have a stronger nonlinearity in the problem. This setup is also considered in \cite{WalN2011}, where $\Omega$ is discretized using $N=128$ grid points in all directions to obtain the nonlinear problem 
\begin{equation*}
    f(x)=Lx - h^2 \lambda \exp(x)=0,
\end{equation*}
displayed in equation \eqref{eq:bratu} of the main text, where $x\in\mathbb{R}^{N^2}$ with $N^2=16,384$.

\begin{figure}[H]
    \myplotsvert{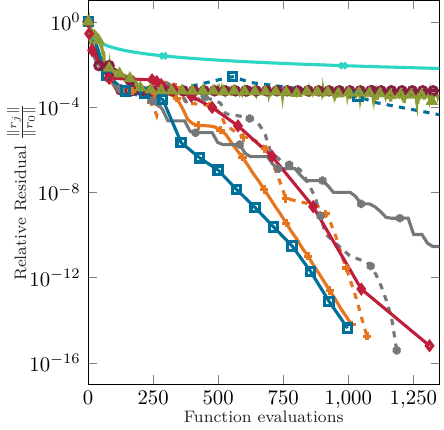}
    \caption{Convergence results for Bratu problem with $\lambda=6$}
    \label{fig:bratu6}
\end{figure}

\medskip
\noindent
\emph{Experiment 1.} We adopt the same parameter choices as in the main text, i.e., nested \textit{nlKrylov} methods use $m=20$, adaptive methods use $\theta=10^{-3}$, Anderson Acceleration is damped with $\beta=0.1$, MINRES is employed as the inner solver within the Newton--Krylov framework, and all methods are terminated once a relative tolerance of $10^{-14}$ is reached. However, we use a relatively large truncation parameter $k=k_{AA}=50$, as suggested in \cite{WalN2011}, for all truncated methods. The results are reported in \Cref{fig:bratu6}.
\noindent We observe that both variants of nlGCR, as well as Anderson Acceleration and nlOrthomin, fail to converge for this value of $\lambda$. The nlLGMRES method performs comparably to Newton--Krylov in terms of outer iterations; however, due to the large truncation window, its inner solves become significantly more expensive, making it less efficient overall (see \Cref{tab:runtime} in the main text). In contrast, nlGMRESR is the most efficient method: while it requires a similar number of function evaluations as nlLGMRES, it is considerably faster since the cost of $\mathcal{SR}_j$ is independent of $k$. Finally, the strong performance of nlGCRO observed for $\lambda=0.5$ in the original experiments is not reproducible in this more strongly nonlinear regime.
\subsection{A nonlinear eigenvalue problem} \label{sec:nep}
\begin{figure}[H]
    \myplotsvert{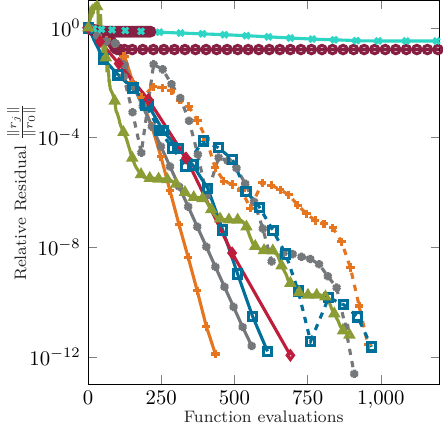}
    \caption{Convergence results for PDDE-NEP}
    \label{fig:nleig}
\end{figure}

In this vector valued example, we consider a nonlinear eigenvalue problem (NEP) originating from the discretization of the partial delay-differential equation (PDDE)
\begin{align*}
    \Delta u(x,t) + a(x)u(x,t) + b(x)u(x,t-2)-u_t(x,t) &=0,&&x\in\Omega,\\
    u(x,t)&=0, &&x\in\partial\Omega,
\end{align*}
where $\Omega=[0,\pi]^2$, $t\geq0$,\;$a(x)=\sin^2(x_1)\sin^2(x_2)$\;and\;$b(x)=\sin(x_1+x_2)+1.31.$
The example is called the \texttt{pdde\_symmetric} problem in the NLEVP database \cite{BetHMST2013,HigPT2019}. Upon discretization, we obtain the NEP
\begin{equation}
  T(\lambda)v=\lambda v, \label{eq:NEP}  
\end{equation}
with a nonlinear function $T(\lambda):\RR\rightarrow\RR^{n,n},~\lambda\mapsto (M+A)+e^{-2\lambda}B,$ where $n=(N-1)^2$, $N$ is the number of discretization points, and $M$, $A$ and $B$ are sparse symmetric $n\times n$ matrices that discretize $\Delta(\cdot)$, $a(\cdot)$ and $b(\cdot)$, respectively. The eigenvalue $\lambda^*$ nearest to zero is of interest in this example.
Using an arbitrary normalization vector $c\in\RR^n$, the eigenvalue problem \eqref{eq:NEP} can be rewritten in terms of a root finding problem \cite{GutT2017}, i.e.,
\begin{equation}
    0=f(x)=\begin{bmatrix} T(\lambda)v-\lambda v\\1-c^Tv       
    \end{bmatrix},\quad x=\begin{bmatrix}v\\\lambda\end{bmatrix}\in\RR^{n+1}. \label{eq:neproot} 
\end{equation}
It is straightforward to see that $T'(\lambda)=-2e^{-2\lambda}B$ and
$$J_f(x)\Delta{x}=\begin{bmatrix} T(\lambda)-\lambda I & (T'(\lambda)-I)v\\-c^T & 0\end{bmatrix}\begin{bmatrix}\Delta{v}\\\Delta\lambda\end{bmatrix} = \begin{bmatrix} (T(\lambda)-\lambda I)\Delta{v}+\Delta\lambda(T'(\lambda)-I)v\\-c^T\Delta{v} 
\end{bmatrix}.$$
\emph{Experiment 1:} In this experiment, we use $N=64$ grid points in both directions, which  results in a problem of size $n+1=(N-1)^2+1=3,970$. We set $c=\frac{1}{\sqrt{n}}\mathrm{1}_{n}$ and initial guess $x_0=[v_0^T,\lambda_0]^T$, where $(v_0,\lambda_0)$ is the smallest magnitude eigenpair of $T(0)=M+A+B$. We use a truncation window of $k=20$ for \textit{nlKrylov} methods and $k_{AA}=30$ for Anderson Acceleration. Nested \textit{nlKrylov} methods use $m=30$, their adaptive versions use $\theta=5\times10^{-3}$. Anderson Acceleration uses the damping parameter $\beta=6\times10^{-4}$ and the Newton--Krylov method allows for a maximum of $200$ GMRES steps. The iterations are stopped once the tolerance of $\tau=\sqrt{n}\cdot10^{-13}$ is reached.\\
In \Cref{fig:nleig}, all nested Krylov methods converge within $31$ iterations. Among them, nlLGMRES required fewer iterations than nlGMRESR and nlGCRO. This behavior is mirrored by their adaptive version, which converged after around ten additional iterations. The methods nlGCR, nlGCR-A, and nlOrthomin did not converge in this experiment. The Newton--Krylov solver converged after just six steps while AA required around 900 steps. The eigenvalue $\lambda^*\approx-0.23937\times10^{-3}$ of interest is found by all converged methods.\\
Looking at function evaluations, it is evident that nlGMRESR, nlLGMRES and nlGCRO all outperform the Newton--Krylov solver, with nlGMRESR requiring around half the function evaluations of Newton--GMRES, while their adaptive counterparts require a number of function evaluations comparable to AA. In terms of total computation time, nlGMRESR significantly outperforms the other methods, requiring only around 70\% of the time needed by its closest competitor, nlGCRO. The adaptive version turns out to be inferior in this experiment. 
\subsection{A nonlinear eigenvector problem (NEPv)} \label{sec:nepv}
\begin{figure}[H]
    \myplotsvert{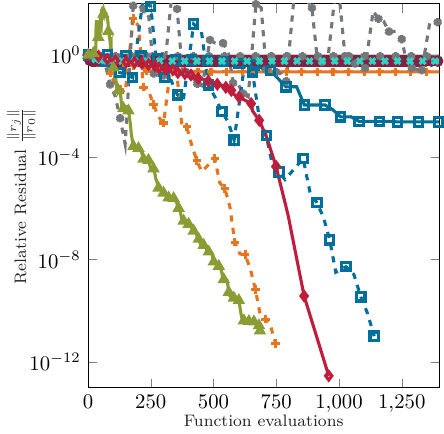}
    \caption{Convergence results for singular Kohn-Sham NEPv}
    \label{fig:nepv}
\end{figure}

\begin{figure}[H]
    \myplotsvert{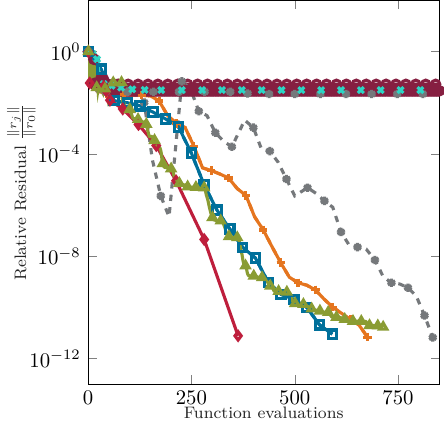}
    \caption{Convergence results for nonsingular Kohn-Sham NEPv}
    \label{fig:nepvnormal}
\end{figure}

For our last experiment, we consider a nonlinear eigenvector problem of the form 
\begin{equation}
    H(V)V=V\Lambda,\quad V^TV=I_p,\quad \Lambda=\Lambda^T, \label{eq:nepv}
\end{equation}
where $H:\RR^{n,p}\rightarrow\RR^{n,n}$ is a symmetric matrix function, i.e., $H(V)=H(V)^T$, and $p\ll n$. Problems of this type arise frequently in quantum physics, most notably discretized Kohn-Sham and Gross-Pitaevskii equations \cite{CaiZBL2018,JarKM2014}, in data science applications such as Robust-Rayleigh-Quotient or Trace-Ratio optimization \cite{BaiLV2018,CaiZBL2018} as well as the modeling of dissipative Hamiltonian DAEs \cite{BaiL2024}. Recently, the use of matrix-valued Newton methods to solve \eqref{eq:nepv} via the root finding problem
\begin{equation}
    F(X)=\begin{bmatrix}H(V)V-V\Lambda\\I_p-V^TV
    \end{bmatrix},\quad X=\begin{bmatrix}V\\\Lambda \end{bmatrix} \in\RR^{(n+p),p}, \label{eq:nepvroot}
\end{equation}
has been studied in \cite{Wer2024}. We consider a simple discrete 3D-Kohn-Sham model 
$$ H(V)=L + \mathrm{Diag}(L^{-1}\rho(V)-\gamma\rho(V)^{1/3}), \quad \gamma\geq 0, $$
where $L\in\RR^{n,n}$ is a 3D-Laplacian on a cube discretized using $N$ equally spaced grid points in every direction, leading to $n=N^3$, and $\rho(V)=\mathrm{diag}(VV^T)$ is the charge density of electrons. Note that the Fr\'echet derivative of \eqref{eq:nepvroot} is singular close to the solution due to the orthogonal invariance of $H(V)$, as was pointed out in \cite{Wer2024}. However, if one chooses a different normalization criterion similar to \eqref{eq:neproot}, i.e., using a constant full rank matrix $C\in\RR^{n,p}$, we can rewrite \eqref{eq:nepvroot} to
\begin{equation}
    \widetilde{F}(X)\begin{bmatrix}H(V)V-V\Lambda\\I_p-C^TV
    \end{bmatrix},\quad X=\begin{bmatrix}V\\\Lambda \end{bmatrix} \in\RR^{(n+p),p}, \label{eq:nepvrootnormal}
\end{equation}
to obtain a nonsingular problem.\\
\emph{Experiment 1:} We first want to examine the singular case, i.e., the root finding problem \eqref{eq:nepvroot}. In our experiment, we use $N=16$ and $p=2$, leading to $X\in\RR^{4098,2}$. The iteration is started with an initial value generated by two steps of SCF initialized using the two smallest eigenpairs of $L$. We use a truncation window of $k=10$ for the \textit{nlKrylov} methods, the inner solve in nested \textit{nlKrylov} methods is carried out by $m=20$ steps of global GMRES, Newton's method uses at most $200$ steps of global GMRES. The iteration is terminated when the relative residual is below $\tau=(N+p)10^{-12}\approx 1.8\cdot10^{-11}$ or after $50$ iterations. The adaptive methods use an angle of $\theta=10^{-3}$. For NEPv, acceleration schemes such as Anderson, Pulay or secant acceleration have proven to improve and stabilize SCF convergence when used as a mixing scheme \cite{ChuDLS2021,ClaM2023,GarS2012}. As a consequence, here, we are using Anderson Acceleration as an accelerator rather than a solver, with SCF as the underlying fixed-point scheme. The window size is $k=30$ and the update is damped using $\beta=-0.1$. The convergence results are displayed in \Cref{fig:nepv}. 
\noindent Here we can see that only four methods managed to converge to the desired accuracy while all the others stagnated earlier. The four convergent methods are nlGMRESR--A, nlLGMRES--A, Newton--Krylov and Anderson Acceleration. Among these methods, the adaptive \textit{nlKrylov} methods required the fewest iterations to convergence: $36$ for nlGMRESR-A and $39$ for nlLGMRES--A, while Newton's method required $49$ steps. Looking at the function evaluations, Anderson acceleration performed slightly better than nlGMRESR-A, while Newton--Krylov and nlLGMRES required over $250$ function evaluations more to achieve convergence. This is also reflected in the total runtime, i.e., AA achieves the lowest computation time at approximately 6.4 seconds, marginally faster than nlGMRESR-A at around 7 seconds, while the remaining methods each exceeded ten seconds. These results underline the observation from the H-equation example that adaptive versions of \textit{nlKrylov} methods exhibit improved convergence compared to their full nonlinear counterparts and can achieve convergence in cases where the standard nonlinear methods stagnate.\\
\noindent\emph{Experiment 2:} In the second experiment, we use the same setup as before to solve the nonsingular problem \eqref{eq:nepvrootnormal}. Here, we adopt the somewhat academic choice of normalization matrix $C=V^*$, where $(V^*,\Lambda^*)$ is an approximate solution to \eqref{eq:nepv} obtained via SCF. The convergence results are displayed in \Cref{fig:nepvnormal}.\\
Compared to the singular case, more methods achieve the desired accuracy, including the full nonlinear versions of nlGMRESR and nlLGMRES, as well as nlGCRO-A. In this setup, Newton--Krylov converges fastest, reaching the solution in just ten steps with nearly quadratic residual reduction with a total computation time of approximately 3.3 seconds. nlLGMRES and nlLGMRES-A follow, both converging after approximately 20 iterations; notably, the adaptive variant did not switch to the linear update. The same holds for nlGMRESR and nlGMRESR-A, both converging in 32 iterations, while nlGCRO-A reached the desired accuracy after 39 steps. A similar trend is observed in the function evaluations, with the exception of Anderson Acceleration requiring roughly the same number of function evaluations as nlGMRESR to achieve the desired accuracy. In terms of computation time, nlLGMRES and nlGMRESR, along with their adaptive counterparts, converge in approximately 5.3 and 6 seconds, respectively, while AA requires around 6.8 seconds and nlGCRO-A around 7.5 seconds. Once again, nlGMRESR and nlLGMRES demonstrate superior convergence behavior relative to nlGCRO and nlGCR on this matrix-valued problem, consistent with the overall trends observed across the extensive numerical experiments.
\subsection{Sensitivity analysis} \label{sec:sensitivity} 
\begin{figure}[H]
    \centering
    \begin{minipage}{0.496\textwidth}
        \vspace{0pt}
        \centering
        \includegraphics[width=\textwidth,height=.3667\textheight]{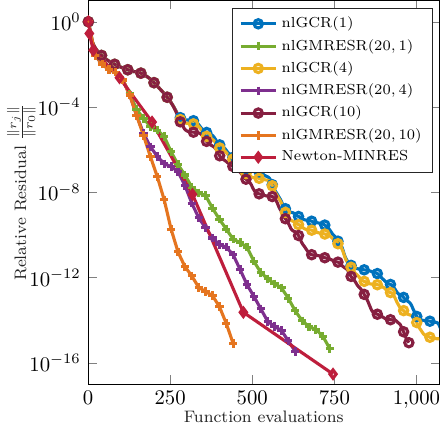}
    \end{minipage}
    \begin{minipage}{0.496\textwidth}
        \vspace{0pt}
        \centering
        \includegraphics[width=\textwidth,height=.3667\textheight]{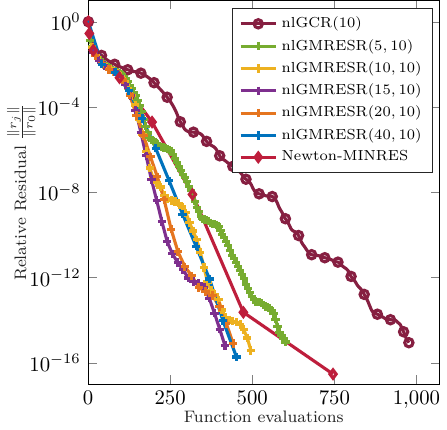}
    \end{minipage}
    \caption{Convergence behavior of nlGCR($k$) and nlGMRESR($m,k$) for different choices of $k$ (left) and $m$ (right) in comparison to Newton--MINRES}
    \label{fig:sensitivity_gcr_gmresr}
\end{figure}

\begin{figure}[H]
    \centering
    \begin{minipage}{0.496\textwidth}
        \vspace{0pt}
        \centering
        \includegraphics[width=\textwidth,height=.3667\textheight]{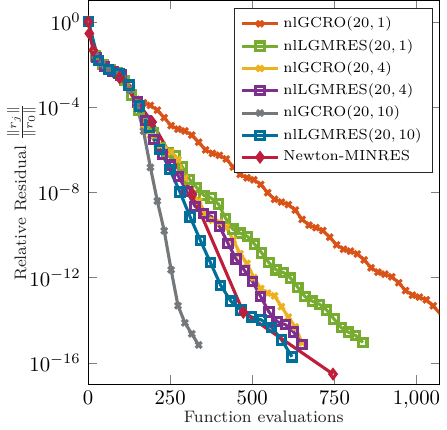}
    \end{minipage}
    \begin{minipage}{0.496\textwidth}
        \vspace{0pt}
        \centering
        \includegraphics[width=\textwidth,height=.3667\textheight]{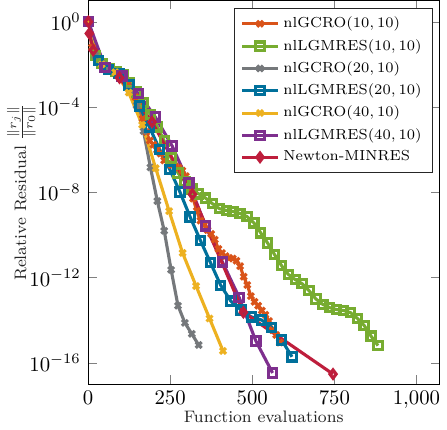}
    \end{minipage}
    \caption{Convergence behavior of nlGCRO($m,k$) and nlLGMRES($m,k$) for different choices of $k$ (left) and $m$ (right) in comparison to Newton--MINRES}
    \label{fig:sensitivity_gcro_lgmres}
\end{figure}

In this section, we present a simplified sensitivity analysis for \textit{nlKrylov} methods applied to the Bratu problem discussed in the main text, where
$$f(x)=Lx-h^2\lambda\exp(x),$$
with $L\in\RR^{n,n}$ a 2D Laplacian, $\lambda=0.5$, $h=\frac{1}{(N+1)^2}$ and $n=N^2$, using $N=100$ grid points per dimension, and initial guess $x_0=\mathbf{1}_n$. All \textit{nlKrylov} methods are benchmarked against Newton--MINRES across four parameter studies. For nlGCR($k$) and nlGMRESR($m,k$), we fix $m=20$ and vary $k\in\lbrace 1,4,10\rbrace$, then we fix $k=10$
and compare nlGCR($10$) to nlGMRESR($m,10$), for $m\in\lbrace 5,10,15,20,40\rbrace$. We repeat the same tests for nlGCRO and nlLGMRES, fixing $m=20$ and compare nlGCRO($20,k$) to nlLGMRES($20,k$) for $k\in\lbrace 1,4,10\rbrace$, then fixing $k=10$ and compare nlGCRO($m,10$) to nlLGMRES($m,10$) for  varying $m\in\lbrace10,20,40\rbrace$. The results are displayed in Figures \ref{fig:sensitivity_gcr_gmresr} and \ref{fig:sensitivity_gcro_lgmres}.\\
The left part of \Cref{fig:sensitivity_gcr_gmresr} shows that increasing $k$ from $1$ to $10$ reduces function evaluations for both nlGCR($k$) and nlGMRESR($20,k$). For nlGCR, the reduction is modest (from $1100$ for $k=1$ to around $1000$ for $k=10$), while for nlGMRESR the effect is more pronounced, with evaluations dropping from roughly $750$ for $k=1$ to $450$ for $k=10$. The right plot, fixing $k=10$, shows that nlGMRESR($m$,$10$) outperforms Newton--MINRES in function evaluations for all tested values of $m$. Performance peaks around $m=15-20$, where the two curves nearly coincide. Increasing to $m=40$ improves performance near convergence but proves counterproductive early in the iteration, where the additional GMRES steps refine the local model but do not translate into meaningfully better nonlinear steps at this stage.
This reflects a general principle, i.e., the closer the local model gets to being linear (the closer the iterates get towards the solution), the more the algorithm benefits from a more accurate linear model, i.e., a larger choice of $m$.
\noindent 
Turning to Figure~\ref{fig:sensitivity_gcro_lgmres}, methods that incorporate the nonlinear basis exhibit a stronger dependence on $m$ and $k$ than nlGMRESR($m,k$). The left plot reveals a striking sensitivity to $k$ for nlGCRO. With $k=1$, nlGCRO($20,1$) behaves similarly to nlGCR($1$), requiring over $1100$ function evaluations; with $k=10$, it converges after around $300$ function evaluations, outperforming all other methods on this problem. nlLGMRES($20,k$) shows a milder but consistent trend, going from roughly $850$ function evaluations for $k=1$ to around $650$ for $k=10$. Notably, at $k=4$, nlLGMRES($20,4$) and nlGCRO($20,4$) perform similarly, both surpassing Newton--MINRES and further improving with an increase in $k$. Fixing $k=10$ and varying $m$, nlGCRO($m,10$) again exhibits a sweet spot near $m=20$, 
while $m=40$ leads to an increase in function evaluations. With $m=10$, nlGCRO($10,10$)  requires a similar number of evaluations as nlLGMRES($20,10$), whereas nlLGMRES($10,10$) needs around $900$ function evaluations to reach the desired tolerance. Increasing to $m=40$ improves nlLGMRES($m,10$) over nlLGMRES($20,10$), particularly near convergence.\\
In summary, $m=20$ and $k=10$ yield good performance across all four algorithms and were therefore adopted in the main experiments. Minor adjustments to these parameters can improve one algorithm while degrading another, making it challenging to select values that ensure a fair comparison; tuning for a single algorithm is, of course, considerably easier. Two key takeaways emerge from this analysis. First of all, nlGMRESR($m,k$) seems to be the most robust and well-behaving algorithm with respect to parameter changes, matching or outperforming Newton--MINRES across all tested values of $m$ and $k$.  Both nlGCRO($m,k$) and nlLGMRES($m,k$) have at least one outlier configuration, where convergence degrades significantly, yet nlGCRO can excel when well-tuned as seen with $m=20$ and $k=10$ for a given problem. The non-nested method nlGCR($k$) is consistent across choices of $k$, but generally inferior to nlGMRESR. Second, nested methods, especially those incorporating the nonlinear basis, benefit from a moderately larger $k$ (e.g., $k=10$) to fully exploit the information provided by the outer model. As noted in the main text, this benefit is most pronounced for moderately nonlinear problems; for highly nonlinear ones, outdated Jacobian information in the inner solve can hinder convergence. For nlGCR($k$), choosing larger $k$ yields little improvement, and small values such as  $k\in\lbrace1,2\rbrace$ are preferable to reduce memory requirements and take advantage of short recurrences in the orthogonalization.   
\bibliographystyle{siamplain}

\end{document}